\documentclass[]{interact}

\usepackage{epstopdf}
\usepackage[caption=false]{subfig}
\usepackage{epstopdf}

\usepackage{graphicx,epsfig}
\usepackage{enumerate}
\usepackage{amssymb}
\usepackage{amsthm}
\usepackage{amsmath}
\usepackage{centernot}
\usepackage{xcolor}
\usepackage{todonotes}
\usepackage{mathrsfs}
\usepackage[utf8]{inputenc}
\usepackage{multirow}
\usepackage{xcolor}
\usepackage{tikz}
\usetikzlibrary{positioning}
\usepackage{pgfplots}
\usepackage{pst-plot}
\usepackage{amsmath}
\usepackage{amsfonts}
\usepackage{amssymb}
\usepackage{multirow}
\usepackage{multicol}
\usepackage{framed}
\usetikzlibrary{angles,quotes}  
 
\usepackage{mathtools}
\usepackage{appendix}
\usepackage{wrapfig}
\usepackage{amsmath}
\usepackage{setspace}
\usepackage{booktabs}
\usepackage{lscape}
\usepackage{algorithm,algorithmicx}
\usepackage{algcompatible}
\usepackage{bm}
\usepackage{graphicx,epsfig}
\usepackage{xfrac}
\usepackage{url}
\usepackage{multirow}
\usepackage{longtable}
\usepackage[linkcolor=blue, urlcolor=blue, citecolor=blue,
colorlinks, bookmarks]{hyperref}
\usepackage[numbers,sort&compress]{natbib}

\hypersetup{colorlinks,%
citecolor=blue,%
filecolor=blue,%
linkcolor=blue,%
urlcolor=blue
}

\bibpunct[, ]{[}{]}{,}{n}{,}{,}
\renewcommand\bibfont{\fontsize{10}{12}\selectfont}
\makeatletter
\def\NAT@def@citea{\def\@citea{\NAT@separator}}
\makeatother

\theoremstyle{plain}
\newtheorem{theorem}{Theorem}[section]
\newtheorem{lemma}[theorem]{Lemma}

\newtheorem{proposition}[theorem]{Proposition}
\newtheorem{assumption}[theorem]{Assumption}

\theoremstyle{definition}
\newtheorem{definition}[theorem]{Definition}
\newtheorem{example}[theorem]{Example}

\theoremstyle{remark}
\newtheorem{remark}{Remark}

\begin{document}

\articletype{ARTICLE TEMPLATE}

\title{Newton Method for Set Optimization Problems with Set-Valued Mapping of Finitely Many Vector-Valued Functions}

\author{
\name{Debdas Ghosh\textsuperscript{a}\thanks{Corresponding Author: D. Ghosh (debdas.mat@iitbhu.ac.in)}, Anshika\textsuperscript{a}, Qamrul Hasan Ansari\textsuperscript{b,}\textsuperscript{c}, Xiaopeng Zhao\textsuperscript{d}}
\affil{\textsuperscript{a}Department of Mathematical Sciences, Indian Institute of Technology (BHU), Varanasi, Uttar Pradesh---221005, India}
\affil{\textsuperscript{b}Department of Mathematics, Aligarh Muslim University, Aligarh 202 002, India}
\affil{\textsuperscript{c}College of Sciences, Chongqing University of Technology, Chongqing 400054, China}
\affil{\textsuperscript{d}School of Mathematical Sciences, Tiangong University, Tianjin, 300387, China}
}

\maketitle

\begin{abstract}
In this paper, we propose a Newton method for unconstrained set optimization problems to find its weakly minimal solutions with respect to lower set-less ordering. The objective function of the problem under consideration is given by finitely many strongly convex twice continuously differentiable vector-valued functions. At first, with the help of a family of vector optimization problems and the Gerstewitz scalarizing function, we identify a necessary optimality condition for weakly minimal solutions of the considered problem. In the proposed Newton method, we derive a sequence of iterative points that exhibits local convergence to a point which satisfies the derived necessary optimality condition for weakly minimal points. To find this sequence of iterates, we formulate a family of vector optimization problems with the help of a partition set concept. Then, we find a \textcolor{black}{descent} direction \textcolor{black}{for this obtained family of vector optimization problems to progress from the current iterate to the next iterate.} As the chosen vector optimization problem differed across the iterates, the proposed Newton method for set optimization problems is not a straight extension of that for vector optimization problems. A step-wise algorithm of the entire process is provided. The well-definedness and convergence of the proposed method are analyzed. \textcolor{black}{To establish the convergence of the proposed algorithm under some regularity condition of the stationary points, we derive three key relations: a condition of nonstationarity, the boundedness of the norm of Newton direction, and the existence of step length that satisfies \textcolor{black}{the Armijo condition}.} We obtain the local superlinear convergence of the proposed method under uniform continuity of the Hessian and local quadratic convergence under Lipschitz continuity of the Hessian. We provide some examples discussing the performance of the proposed method. Also, the performance of the proposed method is compared with the existing steepest-descent method.  
\end{abstract}


\section{Introduction}
A general statement of a set-valued optimization problem is given by 
\begin{equation*} 
 \text{Minimize }F(x)~\text{ subject to }~x\in X,
\end{equation*}
where $F$ is a set-valued map from a nonempty subset $X$ of $\mathbb{R}^n$ to $\mathbb{R}^m$. 
For the solution concepts of set optimization problems, there are two main approaches, namely the vector approach \cite{alonso2008optimality,jahn2004some} and the set approach \cite{alonso2008optimality,kuroiwa1997some}. In the vector approach, the decision maker's preference is based on comparing the vectors in the image set $F(X)$. An alternative definition of solutions of set optimization problems was studied \cite{alonso2008optimality,kuroiwa1997some} via comparing the \textcolor{black}{sets $F(x)$} for all $x \in X$. 

Optimization problems with set-valued maps have extensive applications in various areas, for instance, optimal control, game theory, mathematical economics \cite{bao2010set}, finance \cite{feinstein2015comparison}, differential inclusions \cite{aubin2012differential}, and many others \cite{chen1998optimality,pilecka2016set}. A list of detailed references of mathematical and practical applications of set optimization problems can be found in the introduction section of \cite{bouza2021steepest}.

\textcolor{black}{The research using the solution concept of the set approach} was started with the works in  \cite{kuroiwa1997some,kuroiwa2001set,kuroiwa1997cone}, which considers preorder relations for comparing sets. \textcolor{black}{A detailed discussion on this field is given in \cite{khan2016set}.} The existing methods in the literature to solve set optimization problems fall into one of the following groups: 
\begin{enumerate}
    \item Algorithms \emph{based on scalarization} have been discussed in \cite{ehrgott2014minmax,ide2014concepts}. The methods proposed in these papers address a particular class of set-valued mappings characterized by robust counterparts of vector optimization problems. In \cite{ehrgott2014minmax,ide2014concepts}, a linear scalarization technique was utilized to analyse the optimistic solution of set optimization problem and extended the $\epsilon$-constraint method for the ordering cones with nonnegative orthant to deal with the set optimization problem.
    \item Algorithms of \emph{sorting type} have been discussed in \cite{gunther2019computing,gunther2019strictly,kobis2016treatment,kobis2018numerical}. These methods deal with \textcolor{black}{set optimization} problems having a finite feasible set and rely completely on comparing the images of the set-valued objective mappings. \textcolor{black}{In \cite{kobis2016treatment,kobis2018numerical}, K\"obis et al. extended the work of Jahn \cite{jahn2006multiobjective,jahn2011new} for vector optimization problems. They employed the forward and backward reduction procedures to the algorithms of \cite{jahn2006multiobjective,jahn2011new}, which effectively reduce the number of comparisons compared to the native implementations where every pair of sets needs to be compared. After that, an extension of the algorithm developed by G\"unther and Popovici  \cite{gunther2018new} for vector problems has been discussed in \cite{gunther2019computing,gunther2019strictly}. The approach involves first finding an enumeration of the images of the set-valued mapping, ensuring that their values, when scalarized using a strongly monotone functional, are in an increasing order. This is followed by a forward iteration procedure.} 

    \item Eichfelder et al. \cite{eichfelder2020algorithmic} proposed a \emph{branch and bound} technique for solving multiobjective optimization problems with decision uncertainty based on set-order relation to handling the uncertainties. 
    
    \item Methods based on derivative-free strategy were given in \cite{conn2009introduction}. Jahn \cite{jahn2015derivative,jahn2018derivative} reported \emph{derivative free} descent methods for set optimization using the concepts from \cite{conn2009introduction}. \textcolor{black}{These algorithms are designed to handle unconstrained problems without assuming any particular structure of the set-valued mapping. In \cite{jahn2015derivative}, the method focuses on scenarios where both the epigraphical and typographical multifunction of the set-valued objective mapping exhibit convex values. After that, a relaxation of the convexity assumption was explored in \cite{jahn2018derivative} for upper set less relation. In \cite{jahn2018derivative}, several directions are chosen at once instead of only one. These methods follow a tree generation concept, with roots as the initial point and leaves as the possible solutions. This method is called the rooted tree method.}
    
    \item Bouza et al. \cite{bouza2021steepest} introduced the study of conventional \textcolor{black}{\emph{gradient-based classical approaches} (starting with} the steepest descent method) to solve set optimization problems with finite cardinality. 
\end{enumerate}
We do not repeat again in this article the drawbacks of the methods in the above-mentioned methods in the first four types. These drawbacks are neatly declared in \cite{bouza2021steepest}.

\textcolor{black}{Further, it can be observed that interval-valued optimization problems are considered to be a part of set-valued optimization problems. These are several of Newton's methods present in the literature on interval-valued optimization; see \cite{hansen1983interval,ghosh2016newton,ghosh2017newton,ghosh2017quasi,upadhyay2024newton} and references therein. However, in this article, we have considered a special structure of the set-valued problem (motivated by Bouza et al. \cite{bouza2021steepest}), which differs from the general structure of an interval-valued optimization problem.} 

In this article, we aim to derive a Newton method for the set optimization problems studied in \cite{bouza2021steepest} with a strong convexity assumption. The proposed method in this work exhibits a quadratic convergence near the optimal solution and works well for highly nonlinear objective functions (see Example \ref{exam_4}). 

The paper is organized in the following order. Section \ref{section2} consists of preliminaries, notations, basic results, and definitions that will be used throughout the paper. In Section \ref{section3}, we derive optimality conditions for weakly minimal solutions. In Section \ref{section4}, we propose the Newton method for the considered set optimization problems and analyze its superlinear and quadratic convergence. Section \ref{section5} provides the numerical performance of the proposed method in some test examples. Subsequently, we compare the results of the proposed algorithm with the results of the steepest descent method presented in \cite{bouza2021steepest}. Lastly, we conclude the work in Section \ref{section6} by summarizing our results and proposing ideas for further research.

\section{\textbf{Preliminaries and Definitions}}\label{section2}
Throughout the paper, we use the following notations. 
\begin{itemize}
\item $\mathbb{R},~\mathbb{R}_{+}$, and $\mathbb{R}_{++}$ denote the set of real numbers, nonnegative real numbers, and positive real numbers, respectively.
    \item $\mathbb{R}^m=\mathbb{R}\times\mathbb{R}\times\cdots\times\mathbb{R}$ ($m$-times), $\mathbb{R}_+^m=\mathbb{R}_+\times\mathbb{R}_+\times\cdots\times\mathbb{R}_+$, and $\mathbb{R}_{++}^m=\mathbb{R}_{++}\times\mathbb{R}_{++}\times\cdots\times\mathbb{R}_{++}$.
    \item $\mathcal{P}(\mathbb{R}^m)$ denotes the class of all nonempty subsets of $\mathbb{R}^m$.
    \item For any nonempty set $A\in\mathcal{P}(\mathbb{R}^m)$, the notations int$(A)$, cl($A$),  $|A|$, bd($A$), and conv($A$) denote the interior, closure, cardinality, boundary, and convex hull, respectively, of the set $A$.
    \item $\top$ denotes the transpose operator, and all the elements in $\mathbb{R}^n$ are column vectors.
    \item $\lVert\cdot\rVert$ stands for the standard Euclidean norm of a vector or spectral norm of a matrix. 
    \item For a given $k\in\mathbb{N}$, $[k]$ represents the set $\{1,2,\ldots,k\}$.
    \item A cone $K \subseteq \mathbb{R}^m$ is said to be convex if $K+K=K$, solid if int$(K)\not=\emptyset$, and pointed if $K\cap(-K)=\{0\}$.
    \item Throughout, the notation $K\in\mathcal{P}(\mathbb{R}^m)$ represents a closed, convex, pointed, and solid cone. 
    \item The set $K^*=\{w\in\mathbb{R}^m:~w^\top z\geq0 \text{ for all }z\in K\}$ represents the dual cone of $K$.
    \item \textcolor{black}{$\nabla f^i(x)$ denotes the Jacobian of a vector-valued function $f^i: \mathbb{R}^n \to \mathbb{R}^m$ at $x$.} 
    \item \textcolor{black}{For a vector-valued function $f^i: \mathbb{R}^n \to \mathbb{R}^m$, given by $$f^i(x) : = \left(f^{i, 1}(x), f^{i, 2}(x), \ldots, f^{i, m}(x)\right)^\top,$$ the notation $\nabla^2 f^i(x)$ denotes the follow matrix: 
    $$\nabla^2 f^i(x) := \left(\nabla^2 f^{i, 1}(x), \nabla^2 f^{i, 2}(x), \ldots, \nabla^2 f^{i, m}(x)\right)^\top.$$}
\end{itemize}

Below, we discuss the required fundamental results and basic definitions of set optimization.
\begin{definition}(Partial ordering on $\mathbb{R}^m$ \cite{gopfert2003variational}).\label{partial} The cone $K$ generates a partial order $\preceq$ and a strict order $\prec$ on $\mathbb{R}^m$ defined as follows: for any $y,z\in\mathbb{R}^m$, 
\[
y\preceq z\iff z-y\in K,\text{ and }y\prec z\iff z-y\in \text{ int}(K).
\]
If $y\preceq z$, we often present it by $z \succeq y$. Similarly, $y\prec z$ is often presented by $z \succ y$. 
\end{definition}

\begin{definition}(Lipschitz continuity on $\mathbb{R}^n$).
Let $S$ be a nonempty subset of $\mathbb{R}^n$. A function $F:S\to \mathbb{R}^m$ is said to be Lipschitz continuous on $S$ if there exists a constant $L>0$ such that 
\[
\lVert F(y)-F(x)\rVert\leq L\lVert y-x\rVert \text{ for all }x, y\in S.
\]
\end{definition}

\begin{definition}(Minimal and weakly minimal elements of a set \cite{jahn2009vector}). Let $A\in\mathcal{P}(\mathbb{R}^m)$. The set of minimal and weakly minimal elements of $A$ with respect to $K$ is defined by 
    \begin{eqnarray*} 
    &&\text{Min}(A,K)=\{z\in A:(z-K)\cap A=\{z\}\}\text{ and }\\
    &&\text{WMin}(A,K)=\{z\in A:(z-\text{ int}(K))\cap A=\emptyset\},\text{ respectively}.
     \end{eqnarray*}
\end{definition}

\begin{proposition}\label{min_min} \emph{\cite{jahn2009vector}} Any compact set $A\in \mathcal{P}(\mathbb{R}^m)$ satisfies the domination property with respect to $K$, i.e., $A+K=\text{Min}(A, K) + K$. 
\end{proposition}

Next, we recall the concept of the Gerstewitz scalarizing function, which performs a significant role in the main results of the paper.
\begin{definition}(Gerstewitz function \cite{gerstewitz1983nichtkonvexe}). For an element $e\in$ int($K)$ and $z\in\mathbb{R}^m$, the Gerstewitz function $\Psi_e:\mathbb{R}^m\to \mathbb{R}$ associated with $e$ and $K$ is  defined by 
\[
\Psi_e(z)= \min\{t\in\mathbb{R}~:~te\in z+K\}.
\]
\end{definition}

\begin{proposition}\label{gerstewitz}\emph{(See \cite{khan2016set})}.
For a given element $e\in$ int($K$), the function $\Psi_e:\mathbb{R}^m\to \mathbb{R}$ has the following properties:
 \begin{enumerate}
     \item[(i)] $\Psi_e$ is sublinear on $\mathbb{R}^m$.\label{ge_1} 
     
     \item[(ii)] $\Psi_e$ is positive homogenous of degree 1 on $\mathbb{R}^m$.\label{ge_2} 

     \item[(iii)] $\Psi_e$ is Lipschitz continuous on $\mathbb{R}^m$.\label{ge_3}
     \item[(iv)] $\Psi_e$ is monotone, i.e., for all $p,q\in \mathbb{R}^m$,
     \[
      p\preceq q\implies \Psi_e(p)\leq\Psi_e(q)
     ~~~\text{ and }~~ p \prec q\implies \Psi_e(p)<\Psi_e(q).
     \]  
    
    \item[(v)]\label{ge_5} $\Psi_e$ satisfies the representability property, i.e., 
     \[
     -K=\{z\in\mathbb{R}^m:\Psi_e(z)\leq0\}\text{ and } -\text{int}(K)=\{z\in\mathbb{R}^m:\Psi_e(z)<0\}.
     \]
     \item[(vi)] \label{ge_6} $\Psi_e$ has the translativity property, i.e., $\Psi_e(y + te) = \Psi_e(y) + t ~\forall~ y \in \mathbb{R}^m.$ 
\end{enumerate}
\end{proposition}

Next, we \textcolor{black}{recall the set order relations between the nonempty subsets} of $\mathbb{R}^m$. 

\begin{definition}(Lower set less and strict lower set less relations \cite{kuroiwa1997cone}).\label{set_less}
Let $A$ and $B$ be two sets in $\mathcal{P}(\mathbb{R})^m$. The lower set less $(\preceq^l)$ and strict lower set less $(\prec^l)$ relations for a given \textcolor{black}{cone $K$ are defined} by
\[
A\preceq^l B\iff B\subseteq A+K ~\text{ and }~ A\prec^l B\iff B\subseteq A+\text{ int}(K),\text{ respectively}.
\]
\end{definition}

\begin{framed}
\noindent
In this paper, we aim to derive \textcolor{black}{a} Newton method to identify weakly minimal solutions of the following unconstrained set optimization problem. Let $F:\mathbb{R}^n \rightrightarrows \mathbb{R}^m$ be a nonempty set-valued mapping. The unconstrained set optimization problem that we study is as follows: 
\[
\preceq^l\text{--}\underset{x\in\mathbb{R}^n}{\min}~F(x), \label{sp_equation} \tag{SOP}
\] 
where the solution concept is the notion of weakly minimal solutions with respect to a given ordering cone $K$ in $\mathbb{R}^m$ as given below.

\begin{definition}\label{weakly_minimal_solution}(Weakly minimal solution of (\ref{sp_equation}) \cite{bouza2021steepest}). A point $\bar{x}\in\mathbb{R}^n$ is a local weakly minimal solution of (\ref{sp_equation}) if there exists a neighbourhood $U\subset \mathbb{R}^n$ of $\bar{x}$ such that there does not exist any $x\in U$ with $F(x)\prec^l F(\bar{x}).$ The point $\bar{x}$ is called a weakly minimal solution of (\ref{sp_equation}) if $U=\mathbb{R}^n$.
\end{definition}
 
\begin{assumption} \label{assumption} 
For the sake of deriving \textcolor{black}{our} Newton method for \eqref{sp_equation}, we employ throughout the following assumption for dealing with  \eqref{sp_equation}. 
The function $F:\mathbb{R}^n\rightrightarrows\mathbb{R}^m$ in \eqref{sp_equation} is given by finitely many vector-valued functions: 
\[
F(x)=\left\{f^1(x),f^2(x),\ldots,f^p(x)\right\}\text{ for all }x\in\mathbb{R}^n, 
\]
where $f^1,f^2,\ldots,f^p$ are twice continuously differentiable and $K$-strongly convex functions from $\mathbb{R}^n$ to $\mathbb{R}^m$ with respect to a common $e \in \text{int} (K)$, i.e., there exist positive constants \textcolor{black}{$\rho_1, \rho_2, \ldots, \rho_p$ and $e \in \text{int} (K)$ such that for each $i \in [p]$,}  
\begin{equation}\label{strong_convex_assumption} 
f^{i}(\lambda x_1 + (1- \lambda) x_2) \preceq \lambda f^i(x_1) + (1 - \lambda) f^i(x_2) - \tfrac{1}{2} \rho_i \lambda (1 - \lambda) \|x_1 - x_2\|^2 e  
\end{equation}
for all $x_1, x_2 \in \mathbb{R}^n$ and $\lambda \in [0, 1]$. 
\end{assumption}
\end{framed}

\begin{remark}\label{remark_eigen_value_bigger_than_rho}
It is straightforward to observe that the strong convexity condition \eqref{strong_convex_assumption} in Assumption \ref{assumption} is equivalent to the \textcolor{black}{existence} of positive constants \textcolor{black}{$\rho_1, \rho_2, \ldots, \rho_p$ such that for any $x, u \in \mathbb{R}^n$ and $i \in [p]$, 
\begin{equation}\label{assumption_eigen_value}
u^\top \nabla^2 f^{i}(x) u \succeq \rho_i \|u\|^2 e. 
\end{equation}}
\end{remark}

\section{Optimality Conditions}\label{section3}
In this section, we recall some results on optimality conditions for weakly minimal solutions of (\ref{sp_equation}) under Assumption \ref{assumption}. These results are the foundation for constructing the proposed Newton method to capture weakly minimal solutions of \eqref{sp_equation}. We start with some index-related set-valued mappings analogous to those given in \cite{bouza2021steepest}.

\begin{definition}(Active indices for set-valued maps \cite{bouza2021steepest}).
\begin{enumerate}
    \item[(i)] The set-valued function of active indices of minimal elements associated with the set-valued mapping $F$ of (\ref{sp_equation}) is $I:\mathbb{R}^n\rightrightarrows[p]$, which is defined by 
    \[
    I(x)=\{i\in[p]: f^i(x)\in\text{ Min}(F(x),K)\}.
    \]
    \item[(ii)] The set-valued function of active indices of weakly minimal elements associated with the set-valued mapping $F$ is $I_W:\mathbb{R}^n\rightrightarrows[p]$, which is defined by
    \[
    I_W(x)=\{i\in[p]: f^i(x)\in\text{ WMin}(F(x),K)\}.
    \]
    \item[(iii)] For a vector $r\in\mathbb{R}^m$, the set-valued mapping $I_r:\mathbb{R}^n\rightrightarrows[p]$ is given by 
    \[ 
    I_r(x)=\{i\in I(x): f^i(x)=r\}.
    \] 
    It is to notice that for any $u\in\mathbb{R}^m$, $I_u(x)= \emptyset $ for $u \notin \text{Min}(F(x),K)$;~ $I_u(x) \cap I_v(x) = \emptyset$ for any $u \neq v \in \mathbb{R}^m$;~  $I(x)=\bigcup\limits_{u\in\text{Min}(F(x),K)} I_u(x)$.
\end{enumerate}   
\end{definition}

\begin{definition}\label{cardinality_w}(Cardinality of a set of minimal elements \cite{bouza2021steepest}).
The map $w:\mathbb{R}^n\to \mathbb{N}\cup\{0\}$, which is defined by 
\[w(x)= |\text{Min}(F(x),K)|\]
is called the cardinality of the set of minimal elements of $F(x)$ with respect to $K$. \textcolor{black}{Further, at any $\bar{x}\in\mathbb{R}^n$, we denote $w(\bar{x})=\bar{w}$ for simplicity.} 
\end{definition}

\begin{definition}(Partition set at a point \cite{bouza2021steepest}).
Let us consider an element $x\in\mathbb{R}^n$ and an enumeration 
$\{r_1^x,r_2^x,\ldots,r^x_{w(x)}\}$ of the set $\text{Min}(F(x),K)$. The partition set at $x$ is defined by 
\textcolor{black}{\[
P_x=\prod_{j=1}^{w(x)} I_{r_j^x}(x), \text{ where }  \prod_{j=1}^{w(x)} I_{r_j^x}(x)=I_{r_j^x}(x)\times I_{r_j^x}(x)\times \cdots\times I_{r_j^x}(x).
\]
}
\end{definition}

Throughout the paper, for an iterative point $x_k \in \mathbb{R}^n$, a generic element of the partition set $P_{x_k}$ is denoted by $a^i$, and the components of $a^i$ are denoted by $a^i_j$, where  $j\in[w(x_k)]$. Specifically, if $\lvert P_{x_k}\rvert=p_k$ and Min$(F(x_k),K) = \{r_1^{x_k},r_2^{x_k},\ldots,r^{x_k}_{w(x_{k})}\}$, then 
\[ P_{x_k}=\{a^1,a^2,\ldots,a^{p_k}\}, \]
where for each $i\in[p_k]$, 
\[ a^i=\left(a^i_1,a^i_2,\ldots,a^i_{w(x_k)}\right),~a_i^j \in I_{r_{j}^{x_k}},~j\in [w(x_k)].
\]

Now, we describe a family of vector optimization problems that will help us to find weakly minimal solutions of (\ref{sp_equation}). At a given iterate, one of the vector optimization problems from this family is solved, and the solution is checked for optimality. If the optimality condition is met, then we stop. Otherwise, we solve another vector optimization problem from the family to progress the iteration. Towards finding a stopping condition, we study the following result from \cite{bouza2021steepest}.

\begin{theorem}\emph{(See \cite{bouza2021steepest}).}\label{equivalence_relation}
Let $P_{\bar{x}}$ be the partition set at $\bar{x}$ and $\bar{w} = w(\bar{x})$. For every $a=(a_1,a_2,\ldots,a_{\bar{w}})\in P_{\bar{x}}$, define a vector-valued function  \textcolor{black}{$\widetilde{f}^a:\mathbb{R}^n\to\prod_{j=1}^{\bar{w}} \mathbb{R}^m=\mathbb{R}^{m\bar{w}}$} by 
\[ 
\widetilde{f}^a(x)= \left(
f^{a_1}(x), f^{a_2}(x), \ldots, f^{a_{\bar{w}}}(x)\right)^\top.  
\] 
Let $\widetilde{K}\in\mathcal{P}(\mathbb{R}^{m\bar{w}})$ be the cone given by $\widetilde{K}=\prod_{j=1}^{\bar{w}} K,$ and $\preceq_{\widetilde{K}}$ denote  the partial order in $\mathbb{R}^{m\bar{w}}$ induced by $\widetilde{K}$. Then, $\bar{x}$ is a local weakly minimal solution of (\ref{sp_equation}) if and only if for every $a\in P_{\bar{x}}$, $\bar{x}$ is a local weakly minimal solution of the vector optimization problem 
 \[
\preceq_{\widetilde{K}}\text{--}\underset{x\in\mathbb{R}^n}{\min}~\widetilde{f}^a(x).\tag{VOP}\label{vp_equation}
 \]
\end{theorem}

Next, to find a necessary optimality condition for weakly minimal solutions of (\ref{sp_equation}),  we provide the concept of a stationary point for  (\ref{sp_equation}).

\begin{definition}(Stationary points of \eqref{sp_equation}  \cite{bouza2021steepest}).\label{stationary_dfn} 
A point $\bar{x}$ is called a stationary point of (\ref{sp_equation}) if for every $a=(a_1,a_2,\ldots,a_{\bar{w}})\in P_{\bar{x}}$ and $u\in\mathbb{R}^n$, there exists $j\in[\bar{w}]$ such that 
\begin{eqnarray}\label{stationary_1}  
\nabla f^{a_j}(\bar{x})^\top u \not\in -\text{ int}(K), \text{ i.e., }\Psi_e(\nabla f^{a_j}(\bar{x})^\top u) \geq 0.\label{stationary_2}
\end{eqnarray}   
\end{definition}

\begin{lemma}\label{rtyrrsv}
A point $\bar{x} \in \mathbb{R}^n$ is a stationary point of \eqref{sp_equation}  if and only if  for every $a \in {P}_{\bar x}$, $\bar x$ is a stationary point of \eqref{vp_equation} for every $a \in {P}_{\bar{x}}$. 
\end{lemma}

\begin{proof}
Let $\bar{x}$ be a stationary of  (\ref{sp_equation}). We prove that $\bar{x}$ is a stationary point of \eqref{vp_equation}. On the contrary, suppose $\bar{x}$ is not a stationary point of (\ref{vp_equation}) for some $a \in {P}_{\bar{x}}$. Then, there exists $\bar u \in \mathbb R^{n}$ such that
\begin{align}\label{rsexi}
       & \nabla f^{a_{j}}(\bar{x})^{\top} \bar u \in -\text{int}(K)  ~\forall~ j \in [\bar{w}] \nonumber\\
    \text{ i.e., } & 0 \in \nabla f^{a_{j}}(\bar{x})^{\top} \bar u + \text{int}(K) ~\forall~ j \in [\bar{w}], 
    \end{align}
which implies that 
\[\{ 0\} \subseteq \{ 0\}+ K \overset{(\ref{rsexi})} \subseteq   \{\nabla f^{a_{j}}(\bar{x})^{\top} \bar u\}_{j \in [\bar{w}]}+ \text{int}(K)+K  \subseteq   \{\nabla f^{i}(\bar{x})^{\top} \bar u\}_{i \in [p]}+ \text{int}(K).
\]
Therefore, $ \{\nabla f^{i}(\bar{x})^{\top} \bar u\}_{i \in [p]} \prec^{l} \{ 0\}$, which is contradictory to $\bar x$ \textcolor{black}{being} a stationary point of (\ref{sp_equation}). Hence, $\bar{x}$ must be a stationary point of (\ref{vp_equation}) for every $a \in {P}_{\bar{x}}$. \\

\textcolor{black}{Conversely, suppose $\bar{x}$ is a stationary point of (\ref{vp_equation}) for every $a \in {P}_{\bar{x}}$. We prove that $\bar{x}$ is a stationary point of \eqref{sp_equation}. To prove this, let us assume that $\bar x$ is not a stationary point of (\ref{sp_equation}). Thus, there exists $ \bar u \in \mathbb R^{n}$ such that  
\begin{align*} 
& \nabla f^{i}(\bar{x})^{\top} \bar u \prec 0 ~\forall~ i \in [p]  \\
\Longrightarrow ~ & 0 \in \nabla f^{i}(\bar{x})^{\top} \bar u + \text{int}(K) ~\forall~i \in [p] \\
\Longrightarrow ~ & \nabla f^{a_{j}}(\bar x)^{\top} \bar u \in -\text{int}(K) ~\forall ~j \in [\bar{w}] \subseteq [p], 
\end{align*}
which is contradictory to $\bar{x}$ being a stationary point of (\ref{vp_equation}) for every $a \in P_{\bar x}$. Hence, the result follows.}
\end{proof}

\begin{lemma}\label{weak_min_iff_stationary}
A point $\bar x$ is a stationary point of \eqref{sp_equation} if and only if $\bar x$ is a local weakly minimal solution of \eqref{sp_equation}.     
\end{lemma}

\begin{proof}
Let $\bar x$ be a local weakly minimal solution of \eqref{sp_equation}. On the contrary, suppose $\bar x$ is not a stationary point of (\ref{sp_equation}). Then, by Lemma \ref{rtyrrsv}, $\bar{x}$ is not a stationary point of \eqref{vp_equation} for at least one $\bar{a} = (\bar{a}_{1}, \bar{a}_{2}, \ldots, \bar{a}_{\bar{w}})^{\top} \in P_{\bar x}$. 
Therefore, there exists $\bar u \in \mathbb R^{n}$ such that 
\begin{equation}\label{tt_K}
\nabla f^{\bar a_{j}}(\bar x)^{\top} \bar u \in -\text{int}(K) ~\forall~ j \in [ \bar{w}].
\end{equation}

\noindent
Since $f^{a_{j}}: \mathbb R^{n} \rightarrow \mathbb R^{m}$ is continuously differentiable (Assumption \ref{assumption}) for all  $j \in [\bar{w}]$, we have 
\begin{align}\label{tylor_series}
    f^{\bar a_{j}}(x) = f^{\bar a_{j}}(\bar x) + \nabla f^{\bar a_{j}}(\bar x)^{\top}(x-\bar x)+ o(\lVert x-\bar x\rVert),
\end{align}
where $\lim\limits_{x \to \bar x} \frac{o(\lVert x -\bar x\rVert)}{\lVert x-\bar{x}\rVert}=0.$ As $\bar x$ is a local weakly minimal solution of \eqref{sp_equation}, by Lemma \ref{rtyrrsv}, $\bar x$ is a local weakly minimal point of \eqref{vp_equation} for all $a \in P_{\bar x}$. So, $\bar x$ is a local weakly minimal point of \eqref{vp_equation} for $\bar a \in P_{\bar x}$. Thus, there exists a neighborhood $U$ of $\bar x$ such that 
$$\nexists ~x \in U \text{ with } {{f}^{\bar a_j}}(x) - {{f}^{\bar a_{j}}}(\bar x) \notin -\text{int}(K) \text{ for all } j \in [\bar w ].$$ 
From (\ref{tylor_series}), there exists a neighborhood $B \subseteq U$ of $\bar x$  such that for all $j \in [\bar w]$, 
     \begin{align}\label{tyu}
     & \nabla f^{\bar a_{j}}(\bar x)^{\top}(x-\bar x) \notin -\text{int} (K) ~~\forall~ x \in B
     \end{align}
As $B$ is a neighborhood of $\bar x$, there exists $\bar t > 0$ such that $x' = \bar x + \bar t u \in B$. The relation \eqref{tyu} with $x = x'$ yields 
\[\nabla f^{\bar a_{j}}(\bar x)^{\top} \bar u \notin -\text{int} (K), 
\] 
which is contradictory to \eqref{tt_K}. Therefore, $\bar x$ is a stationary point for (\ref{sp_equation}).  \\

Conversely, suppose that $\bar x$ is a stationary point of \eqref{sp_equation}. As ${F}$ is strongly $K$-convex (Assumption \ref{assumption}), we have for all $ j \in [\bar{w}], x \in \mathbb R^{n}, \alpha \in (0,1]$ that 
\begin{align}
~&~  \textcolor{black}{\alpha} f^{a_{j}}(x)+(1-\alpha)f^{a_{j}}(\bar x)- f^{a_{j}}(\alpha x+ (1-\alpha)\bar x)\in { K} \notag \\
\textcolor{black}{\overset{\text{Def. \ref{partial}}}{\implies}} ~&~ \textcolor{black}{(f^{a_{j}}(x)-f^{a_{j}}(\bar x))\preceq\tfrac{1}{\alpha}(f^{a_{j}}( \bar x+\alpha(x-\bar x))-f^{a_{j}}(\bar x))}\notag\\
\textcolor{black}{\overset{\text{Def. \ref{partial}}}{\implies}}  ~&~\textcolor{black}{(f^{a_{j}}(x)-f^{a_{j}}(\bar x))-\tfrac{1}{\alpha}(f^{a_{j}}( \bar x+\alpha(x-\bar x))-f^{a_{j}}(\bar x)) \in {K} }\notag \\ 
\overset{\alpha \to 0+}{\implies}  ~&~  (f^{a_{j}}(x)-f^{a_{j}}(\bar x))- \nabla f^{a_{j}}(\bar x)(x-\bar x)\in {K} \notag \\
\implies ~&~ \nabla f^{a_{j}}(\bar x)(x-\bar x) \preceq (f^{a_{j}}(x)-f^{a_{j}}(\bar x)) \notag \\ 
\implies ~&~ \Psi_e(\nabla f^{a_{j}}(\bar x)(x-\bar x) ) \leq \Psi_e( f^{a_{j}}(x)-f^{a_{j}}(\bar x)).
\label{cri} 
\end{align}
Since $\bar x$ is a stationary point for \eqref{sp_equation}, by Lemma \ref{rtyrrsv}, $\bar{x}$ is stationary point for (\ref{vp_equation}) for all $a \in P_{\bar x}$. Thus, for any $u = x - \bar x \in \mathbb R^{n}$, there exists $j_{u} \in [\bar w]$ such that $\Psi_e(\nabla f^{a_{j_{u}}}(\bar x)^{\top}u) \geq 0$. 
So, from \eqref{cri} we get for any $x \in \mathbb R^{n}$ that 
\begin{align*}
  ~&~ 0 \le \Psi_e(\nabla f^{a_{j_{u}}}(\bar x)( x- \bar x) ) \leq \Psi_e( f^{a_{j_{u}}}(x) - f^{a_{j_{u}}}(\bar x)) \\
  \overset{\text{Lemma \ref{gerstewitz}}}{\implies} &  f^{a_{j_{u}}}( x)-f^{a_{j_{u}}}(\bar x) \notin - \text{int }{K} \\
  \implies ~&~ f^{a_{j_{0}}}(x) \nprec f^{a_{j_{0}}}(\bar x).
\end{align*}
Therefore, there exists $j_{u}\in [\bar w]$ for which there does not exist any $x \in \mathbb R^{n}$ such that $f^{a_{j_{u}}}(x) \prec f^{a_{j_{0}}}(\bar x)$. This implies $\bar x$ is a weakly minimal element of (\ref{vp_equation}) for all $a \in P_{\bar x}$. Therefore, by Lemma \ref{rtyrrsv}, $\bar{x}$ is a weakly minimal solution of (\ref{sp_equation}). 
 \end{proof}

\begin{lemma}\label{g_is_strong_convex}
For any given $x \in \mathbb{R}^n$, $a \in P_x$ and $j \in [w(x)]$, the function $g: \mathbb{R}^n \to \mathbb{R}$ given by 
\[g(u) = \Psi_{e}\left(\nabla f^{a_j}(x)^{\top}u + \tfrac{1}{2}u^\top\nabla^2 f^{a_j}(x)u\right)\]
is strongly convex on $\mathbb{R}^n$. 
\end{lemma}

\begin{proof}
Consider the function 
\[
h_j(u) = \nabla f^{a_j}(x)^{\top}u + \tfrac{1}{2}u^\top\nabla^2 f^{a_j}(x)u,~ u \in \mathbb{R}^n. 
\]
According to Remark \ref{remark_eigen_value_bigger_than_rho}, \textcolor{black}{we have $u^\top \nabla^2 h_j(x) u \succeq \rho_j \|u\|^2 e$.} Therefore, by Corollary 2.2 in \cite{drummond2014quadratically}, the function $h_j: \mathbb{R}^n \rightarrow \mathbb{R}^m$ is strongly convex. Hence, there exists \textcolor{black}{$\rho > 0$} such that for any $u_1, u_2 \in \mathbb{R}^n$ and $\lambda \in [0, 1]$,  
\begin{equation}\label{hj_ineq}
h_j(\lambda u_1 + (1 - \lambda)u_2) \preceq \lambda h_j(u_1) + (1 - \lambda)h_j(u_2) - \tfrac{\textcolor{black}{\rho}}{2}\lambda (1 - \lambda) \|u_1 - u_2\|^2e. 
\end{equation}
Therefore, for any $u_1, u_2 \in \mathbb{R}^n$ and $\lambda \in [0, 1]$, we have 
\begin{align*}
~&~ g(\lambda u_1 + (1 - \lambda)u_2) \\ 
= ~&~ \Psi_e(h_j(\lambda u_1 + (1 - \lambda)u_2)) \\ 
\textcolor{black}{\overset{\text{\ref{gerstewitz}(iv)\&\eqref{hj_ineq}}}{\le}} ~&~\textcolor{black}{ \Psi_e(\lambda h_j(u_1) + (1 - \lambda)h_j(u_2) - \tfrac{\rho}{2}\lambda (1 - \lambda) \|u_1 - u_2\|^2e)}\\ 
\textcolor{black}{\overset{\text{\ref{gerstewitz}(i)\&(ii)}}{\le}} ~&~ \textcolor{black}{\lambda \Psi_e(h_j(u_1)) + (1 - \lambda) \Psi_e(h_j(u_2)) - \tfrac{\rho}{2}\lambda (1 - \lambda) \|u_1 - u_2\|^2}  \\
= ~&~ \lambda g(u_1) + (1 - \lambda) g(u_2) - \tfrac{\rho}{2}\lambda (1 - \lambda) \|u_1 - u_2\|^2. 
\end{align*} 
Hence, $g$ is strongly convex on $\mathbb{R}^n$.  
\end{proof}

\begin{framed}
\noindent 
Next, we discuss a necessary condition for weakly minimal solutions of (\ref{sp_equation}). From Lemma \ref{weak_min_iff_stationary}, we note that a weakly minimal solution of \eqref{sp_equation} is a stationary point of \eqref{sp_equation} and vice-versa. From Definition \ref{stationary_dfn},    
\begin{align*}
&~~ \text{A point } \bar{x} \text{ is a stationary point of (\ref{sp_equation})} \\ 
\Longleftrightarrow &~~ \textcolor{black}{\forall} a \in P_{\bar x} \text{ and } u \in \mathbb{R}^n,  ~\exists a_j \text{ with } \Psi_e(\nabla f^{a_j}(\bar x)^\top u) \ge 0. 
\end{align*}
So, by Proposition \ref{ge_6} and Remark \ref{remark_eigen_value_bigger_than_rho}, for a stationary point $\bar x$, for any $a \in P_{\bar x}$ and $u \in \mathbb{R}^n$, there exists $a_j$ with 
\textcolor{black}{\begin{align}\label{psi_e_ge_0}
&u^\top \nabla^2 f^{a_{j}}(\bar x) u\succeq \rho_{a_j} \|u\|^2 e\notag\\
\text{or, } &\nabla f^{a_{j}}(\bar x)^\top u+u^\top \nabla^2 f^{a_{j}}(\bar x) u\succeq \nabla f^{a_{j}}(\bar x)^\top u+\rho_{a_j} \|u\|^2 e\notag\\
\text{or, } &\Psi_e\left(\nabla f^{a_{j}}(\bar x)^\top u+u^\top \nabla^2 f^{a_{j}}(\bar x) u\right)\ge\Psi_e\left(\nabla f^{a_{j}}(\bar x)^\top u+\rho_{a_j} \|u\|^2 e\right)\text{ by \ref{ge_6}(iv)}\notag\\
\text{or, }&\Psi_e\left(\nabla f^{a_{j}}(\bar x)^\top u + u^\top \nabla^2 f^{a_{j}}(\bar x) u \right)   
\ge \Psi_e\left(\nabla f^{a_{j}}(\bar x)^\top u\right) + \tfrac{1}{2}\rho_{a_j} \|u\|^2 \ge 0 \text{ by \ref{ge_6}(vi)},
\end{align}}
where $\rho = \min\{\rho_1, \rho_2, \ldots, \rho_p\}$. So, for any $x \in \mathbb{R}^n$, if we define a function $\xi_x:P_x \times \mathbb{R}^n \to\mathbb{R}$ by 
\begin{align}
\xi_x(a,u) = \underset{j\in [w(x)]}{\max}\left\{\Psi_{e}(\nabla f^{a_j}(x)^{\top}u + \tfrac{1}{2}u^\top\nabla^2 f^{a_j}(x)u)\right\}, a\in P_x,~u\in\mathbb{R}^n, \label{function}   
\end{align}
then by \eqref{psi_e_ge_0}, at a stationary point $\bar x$, we have 
\begin{align}\label{xi_bar_x_ge_0} 
~&~ \xi_{\bar{x}}(a, u) \ge 0 \text{ for all } a \in P_{\bar x} \text{ and } u \in \mathbb{R}^n \notag \\ 
\Longrightarrow ~&~ \min_{u \in \mathbb{R}^n} \xi_{\bar{x}}(a, u) \ge 0 \text{ for all }~ a \in P_{\bar x} \notag \\ 
\Longrightarrow ~&~ \forall a \in P_{\bar x}: 0 \le \min_{u \in \mathbb{R}^n} \xi_{\bar{x}}(a, u) \le \xi_{\bar{x}}(a, 0) = 0 \notag \\ 
\Longrightarrow ~&~ \forall a \in P_{\bar x}: \min_{u \in \mathbb{R}^n} \xi_{\bar{x}}(a, u) = 0. 
\end{align}

\noindent
Moreover, as for any $x \in \mathbb{R}^n$, $P_x$ is finite, we note from Lemma \ref{g_is_strong_convex} that for any $a\in P_x$, the function $\xi_x(a,\cdot)$ is strongly convex in $\mathbb{R}^n$. Hence, the function $\xi_{\bar x}(a,\cdot)$ has a unique minimum over $\mathbb{R}^n$. If for $a\in P_{\bar x}$, $\bar{u}_{a, \bar{x}}\in\mathbb{R}^n$ be such that $\xi_{\bar x}(a,\bar{u}_{a, \bar{x}})=\underset{u\in\mathbb{R}^n}{\min}~\xi_{\bar x}(a,u)$, then from \eqref{xi_bar_x_ge_0}, we have 
\begin{eqnarray}
\xi_{\bar x}(a,\bar{u}_{a, \bar{x}}) = 0 \text{ if and only if }\bar{u}_{a, \bar{x}} = 0. \label{function_4}
\end{eqnarray}
As for any $x \in \mathbb{R}^n$, the partition set $P_{x}$ is finite, $\xi_{x}$ attains its minimum over the set $P_{x}\times\mathbb{R}^n$. Let us define a function $\Phi:\mathbb{R}^n\to\mathbb{R}$ by	
\begin{eqnarray}
&&\Phi(x)=\underset{(a,u)\in P_{x} \times \mathbb{R}^n}{\min}~\xi_x(a,u)\label{functionn}.
\end{eqnarray}
Then, for any $x \in \mathbb{R}^n$, 
\begin{equation}\label{function_3}
\Phi(x) =\underset{(a,u)\in P_{x} \times \mathbb{R}^n}{\min}~\xi_x(a,u) \le \xi_x(a,0) = 0.
\end{equation}
Also, in view of (\ref{function_4}) and \eqref{xi_bar_x_ge_0}, if for $(a,\bar{u}) \in P_{\bar x}\times\mathbb{R}^n$ we have $\Phi(\bar x)=\xi_{\bar x}(\bar{a},\bar{u})$, then   
\begin{eqnarray}\label{function_2}
\Phi(\bar x)= 0 \text{ and } \bar{u}=0. 
\end{eqnarray}
Accumulating all, we obtain the following result.

\begin{proposition}\emph{(Necessary condition for weakly minimal points)}.\label{prop_minimal}
Let $\bar{x}$ be a weakly minimal point of (\ref{sp_equation}) and  $\bar{a}\in P_{\bar{x}}$ and $\bar{u}\in\mathbb{R}^n$ be such that $\Phi(\bar{x})=\xi_{\bar{x}}(\bar{a},\bar{u})$, where $\xi_{\bar{x}}$ and $\Phi$ are as defined in (\ref{function}) and (\ref{functionn}), respectively. Then, $\bar{u}=0$.  
\end{proposition}
\end{framed}


Next, we derive a few properties of $\Phi$, which play an important role in the convergence analysis of the proposed Newton method for \eqref{sp_equation}.

\begin{proposition}\label{continuity}
The function $\Phi$ as given in (\ref{functionn})  is continuous at any  $\bar{x}\in \mathbb{R}^n$.    
\end{proposition}

\begin{proof}
Let $\{x_k\}$ be a sequence in $\mathbb{R}^n$ that converges to $\bar{x}\in\mathbb{R}^n$. We show that 
\[\lim_{k\to\infty}\Phi(x_k) = \Phi(\bar{x}).
\] 
Since the set $P_{\bar x}$ is finite and $\xi_{\bar x}$ attains its minimum over the set $P_{\bar x}\times\mathbb{R}^n$, there exists $(\bar a, \bar u) \in P_{\bar x} \times \mathbb{R}^n$ such that $\Phi(\bar x) = \xi_{\bar x} (\bar a, \bar u)$. \\ \\ 
Let $(a^k, u_k)$ be an element in $P_{x_k} \times \mathbb{R}^n$ such that $\Phi(x_k) = \xi_{x_k}(a^k,u_k)$. Such an element $(a^k, u_k)$ exists since the set $P_{x_k}$ is finite and $\xi_{x_k}$ attains its minimum over the set $P_{x_k}\times\mathbb{R}^n$. 
Since $\Psi_e$ is Lipschitz continuous on $\mathbb{R}^n$ (Proposition \ref{gerstewitz} (iii)) and \textcolor{black}{$f^{a^k_j}$ is twice continuously differentiable for each $j\in[w(x_k)]$, therefore the function $\xi_{x_k}$ is continuous on $P_{x_k} \times \mathbb{R}^n$. Thus, we get 
\begin{eqnarray}\label{continuity_1} 
\underset{k\to \infty}{\limsup}~\Phi(x_k) = \underset{k\to \infty}{\limsup}~\xi_{x_k}(a^k,u_k) \le \underset{k\to \infty}{\limsup}~\xi_{x_k}(\bar a, \bar u) =\xi_{\bar{x}}(\bar{a},\bar{u})
=\Phi(\bar{x}).
\allowdisplaybreaks
\end{eqnarray}}
Let $L>0$ be the Lipschitz constant of $\Psi_e$. Then, from the Definition (\ref{functionn}) of $\Phi$ at $\bar{x}$, we observe that 
\allowdisplaybreaks
\begin{align}\label{cont_dg}
&~ \Phi(\bar{x}) \nonumber \\ 
=&~  \underset{(a,u)\in P_{\bar{x}}\times\mathbb{R}^n}{\min}\xi_{\bar{x}}(a,u)\nonumber\\
\leq &~\xi_{\bar{x}}(a^k,u_k)\nonumber\\
=&~\underset{k\to \infty}{\liminf}~\xi_{\bar{x}}(a^k,u_k)\text{ since }\xi_{\bar{x}}\text{ is continuous}\nonumber \\
=&~\underset{k\to \infty}{\liminf}~\left\{\underset{j\in[w(\bar{x})]}{\max}\left(\Psi_e(\nabla f^{a^k_j}(\bar{x})^\top u_k+\tfrac{1}{2}u_k^\top\nabla^2f^{a^k_j}(\bar{x})u_k)\right)\right\}\nonumber\\
=&~\underset{k\to \infty}{\liminf}~\left\{\underset{j\in[w(\bar{x})]}{\max}\left(\Psi_e(\nabla f^{a^k_j}(\bar{x})^\top u_k+\tfrac{1}{2}u_k^\top\nabla^2f^{a^k_j}(\bar{x})u_k+\nabla f^{a^k_j}(x_k)^\top u_k\right.\right.\nonumber\\
&~\quad \quad \quad \quad \left.\left. +\tfrac{1}{2}u_k^\top\nabla^2f^{a^k_j}(x_k)u_k-\nabla f^{a^k_j}(x_k)^\top u_k-\tfrac{1}{2}u_k^\top\nabla^2f^{a^k_j}(x_k)u_k)
\right) \right\} \nonumber \\ 
\overset{\ref{gerstewitz}\text{(i)}}{\le}&~\underset{k\to \infty}{\liminf}~\left\{\underset{j\in[w(\bar{x})]}{\max}\left\{\Psi_e\left(\nabla f^{a^k_j}(x_k)^\top u_k + \tfrac{1}{2} u_k^\top\nabla^2f^{a^k_j}(x_k)u_k\right) \right. \right.  \nonumber\\
&~\left. \left. + \Psi_e\left(\nabla f^{a^k_j}(\bar{x})^\top u_k +\tfrac{1}{2}u_k^\top\nabla^2f^{a^k_j}(\bar{x})u_k-\nabla f^{a^k_j}(x_k)^\top u_k-\tfrac{1}{2}u_k^\top\nabla^2f^{a^k_j}(x_k)u_k
\right)\right\}\right\}\nonumber\\
=&~\underset{k\to \infty}{\liminf}~\left\{\xi_{x_k}(a^k,u_k) +\underset{j\in[w(\bar{x})]}{\max}\left\{\Psi_e\left(\nabla f^{a^k_j}(\bar{x})^\top u_k \right.\right.\right.\nonumber\\
&~ \left. \left. \left. \qquad~~ + \tfrac{1}{2}u_k^\top\nabla^2f^{a^k_j}(\bar{x})u_k-\nabla f^{a^k_j}(x_k)^\top u_k - \tfrac{1}{2}u_k^\top\nabla^2f^{a^k_j}(x_k)u_k\right)\right\}\right\}\nonumber\\
\overset{\ref{gerstewitz}\text{(iii)}}{\le} & \underset{k\to \infty}{\liminf}~\left\{\xi_{x_k}(a^k,u_k) + L\underset{j\in[w(\bar{x})]}{\max}\left\{\lVert\nabla f^{a^k_j}(\bar{x})^\top u_k+\tfrac{1}{2}u_k^\top\nabla^2f^{a^k_j}(\bar{x})u_k \right. \right.\nonumber\\
&~\quad \quad \quad \left.\left.-\nabla f^{a^k_j}(x_k)^\top u_k -\tfrac{1}{2}u_k^\top\nabla^2f^{a^k_j}(x_k)u_k\rVert\right\}\right\} \nonumber\\
\leq&~\underset{k\to \infty}{\liminf}\left\{~\xi_{x_k}(a^k,u_k) + L~\underset{j\in[w(\bar{x})]}{\max}\{\lVert\nabla f^{a^k_j}(\bar{x})-\nabla f^{a^k_j}(x_k)\rVert\lVert u_k\rVert\}\right. \nonumber \\
&~\left.+\tfrac{L}{2}\underset{j\in[w(\bar{x})]}{\max}\left\{\lVert u_k^\top\left(\nabla^2f^{a^k_j}(\bar{x})-\nabla^2f^{a^k_j}(x_k)\right)u_k\rVert\right\}\right\}. 
\end{align}
Note that for \textcolor{black}{$j\in[w_k]$,} each $f^{a^k_j}$ is a twice continuously differentiable and the sequence $\{x_k\}$ converges to $\bar{x}$. Also, note that there is no loss of generality if $\{u_k\}$ is assumed to be in $\{u \in \mathbb{R}^n: \|u\| \le 1\}$. Thus, we obtain from \eqref{cont_dg} that 
\begin{eqnarray}\label{continuity_3}
\Phi(\bar{x})\leq \underset{k\to \infty}{\liminf}~\xi_{x_k}(a^k,u_k)=\underset{k\to \infty}{\liminf}~\Phi(x_k).
\end{eqnarray}
Finally, in view of (\ref{continuity_1}) and (\ref{continuity_3}), we conclude that
\[
\underset{k\to\infty}{\lim}\Phi(x_k)=\Phi(\bar{x}).
\]
Thus, the function $\Phi$ is continuous at $\bar{x}$.
\end{proof}

\begin{proposition}\label{bounded}
Let $U$ be a nonempty subset of $\mathbb{R}^n$. Suppose there exists $\gamma\in\mathbb{R}_{++}$ such that for any $x\in U$ and $a\in P_x$, $\nabla^2 f^{a_j}(x)\leq \gamma I$ for all $j\in[w(x)]$. Then, for any $a \in P_{x}$, there exist $\lambda_j\geq0$, $j\in [w(x)]$, with $\sum_{j=1}^{[w(x)]}\lambda_j=1$ such that  
\[
\lvert\Phi(x)\rvert\leq\tfrac{3L}{2\gamma }\left\lVert\sum_{j=1}^{[w(x)]}\lambda_j\nabla f^{a_j}(x)\right\rVert^2, 
\]
where $L$ is the Lipschitz constant of $\Psi_e$.
\end{proposition}

\begin{proof}
Let $x\in U$ and $P_x$ be the partition set of (\ref{sp_equation}) at $x$. Note that for any $b_1,b_2,\ldots,b_{w(x)}\in\mathbb{R}$, the identity $\max\{b_1,b_2,\ldots,b_{w(x)}\}=\underset{\lambda\in \Delta_{w(x)}}{\max}\sum_{j=1}^{w(x)} \lambda_j b_j $ holds, where $\Delta_{w(x)} = \{(\lambda_1, \lambda_2, \ldots, \lambda_{w(x)}) \in \mathbb{R}^{w(x)}_+: \sum_{j = 1}^{w(x)} \lambda_j = 1\}$. Thus, in view of definition (\ref{functionn}) of $\Phi$, we have for any $a \in P_{x}$ that 
\allowdisplaybreaks 
\begin{eqnarray} 
|\Phi(x)|&=& \left|\underset{(a,u)\in P_x\times\mathbb{R}^n}{\min}\xi_{x}(a,u) \right|\nonumber \\
&=&\left| \underset{(a,u)\in P_x\times\mathbb{R}^n}{\min}~\left\{\underset{j\in[w(x)]}{\max}\Psi_e\left(\nabla f^{a_j}(x)^\top u+\tfrac{1}{2}u^\top\nabla^2f^{a_j}(x)u\right)\right\}\right|\nonumber\\ 
&=&\underset{(a,u)\in P_x\times\mathbb{R}^n}{\min}~\left|\underset{\lambda\in \Delta_{w(x)}}{\max}\sum_{j=1}^{[w(x)]}\lambda_j\Psi_e\left(\nabla f^{a_j}(x)^\top u+\tfrac{1}{2}u^\top\nabla^2f^{a_j}(x)u\right)\right|\nonumber\\ 
&\leq&\underset{(a,u)\in P_x\times\mathbb{R}^n}{\min}~\sum_{j=1}^{w(x)} \left|\lambda_j\Psi_e\left(\nabla f^{a_j}(x)^\top u+\tfrac{1}{2}u^\top\nabla^2f^{a_j}(x)u\right)\right| \notag \\ 
&& \text{ for some } \lambda \in \Delta_{w(x)} \nonumber \\ 
&\overset{\ref{gerstewitz}\text{(iii)}}{\le} & \underset{(a,u)\in P_x\times\mathbb{R}^n}{\min}~\sum_{j=1}^{w(x)} \lambda_j L \left\lVert\nabla f^{a_j}(x)^\top u+\tfrac{1}{2}u^\top\nabla^2f^{a_j}(x)u\right\rVert \notag \\ 
&\leq& L \underset{(a,u)\in P_x\times\mathbb{R}^n}{\min}~\left\{\sum_{j=1}^{w(x)} \lambda_j\left\lVert\nabla f^{a_j}(x)^\top u\right\rVert + \sum_{j=1}^{w(x)} \lambda_j \left\lVert\tfrac{1}{2}u^\top\nabla^2f^{a_j}(x)u\right\rVert \right\} \nonumber \\ 
&\leq& L\underset{(a,u)\in P_x\times\mathbb{R}^n}{\min}\left\{\sum_{j=1}^{w(x)} \lambda_j \lVert\nabla f^{a_j}(x)^\top u\rVert + \tfrac{\gamma }{2}\lVert u\rVert^2\right\} \label{dg_aux1} \text{as }  \nabla^2 f^{a_j}(x) \leq \gamma I. 
\end{eqnarray} 
As the function $u \mapsto \sum_{j=1}^{w(x)}\lambda_j\lVert\nabla f^{a_j}(x)^\top u\rVert+\tfrac{\gamma }{2}\lVert u\rVert^2$ is a strongly convex function on $\mathbb{R}^n$, the first-order optimality condition implies that its minimum is obtained at \textcolor{black}{$u=-\tfrac{1}{\gamma }\sum_{j=1}^{w(x)}\lambda_j\nabla f^{a_j}(x)$.} Thus, \eqref{dg_aux1} gives for any $a \in P_{x}$ that there exist  $\lambda_j\geq0$, $j\in w(x)$ with $\sum_{j=1}^{w(x)}\lambda_j=1$ such that 
\[
\lvert\Phi(x)\rvert\leq\tfrac{3L}{2\gamma }\left\lVert\sum_{j=1}^{[w(x)]}\lambda_j\nabla f^{a_j}(x)\right\rVert^2. 
\]
\end{proof}

Below, we define the notion of the regularity of a point with an essential property for a set-valued mapping, which has a significant role in the convergence of the proposed algorithm.

\begin{definition}\label{regular_condition}
(Regular point \cite{bouza2021steepest}).
Let $U$ be a nonempty subset of $\mathbb{R}^n$. A point $\bar{x}\in U$ is said to be a regular point of $F$ if it satisfies the following conditions:
\begin{enumerate}
    \item[(i)] Min $(F(\bar{x}),K)=$ WMin $(F(\bar{x}),K)$, and 
    \item[(ii)] the cardinality function $w$ in Definition \ref{cardinality_w} is constant in the neighbourhood of $\bar{x}$.
\end{enumerate}
\end{definition}

\begin{lemma}\emph{(See \cite{bouza2021steepest})}.\label{regular}
Let us assume that $\bar{x}\in \mathbb{R}^n$ is a regular point of $F$. Then, there exists a neighbourhood $U$ of $\bar{x}$ such that for every $x\in U$, $w(x)=\bar{w}$, and $P_x\subseteq P_{\bar{x}}.$
\end{lemma}


\section{Newton Method and Its Convergence Analysis}\label{section4}
In this section, we propose a Newton method (Algorithm \ref{algo1}) for set optimization problems (\ref{sp_equation}) with an $F$ as given in Assumption \ref{assumption}. We start the algorithm by selecting an arbitrary initial point. If this point does not satisfy the necessary condition for a weakly minimal point as stated in Proposition \ref{prop_minimal}, then we proceed to update this point as discussed in Algorithm \ref{algo1}. At each iteration, we select an element from the partition set of the current point, and \textcolor{black}{then we evaluate a descent direction for (\ref{vp_equation}) by following the ideas of \cite{drummond2014quadratically,fliege2009newton}.} Once a descent direction is found, we employ a backtracking procedure similar to the classical Armijo-type method to find an appropriate step size and then update the iterate. We keep updating the iterate until the necessary condition in Proposition \ref{prop_minimal} for a weakly minimal point is met. The entire method is given in Algorithm \ref{algo1}.  

\begin{algorithm}[!htb]
\caption{Newton Method for Set Optimization Problem (\ref{sp_equation})  \label{algo1}}
\begin{enumerate}[\textit{Step} 1]
\item \label{step0}  \textbf{Inputs}\\ Provide the objective function $F$ with $f^1,f^2,\ldots,f^p$ being twice continuously differentiable and strongly convex vector-valued functions satisfying Assumption  \ref{assumption}.\\
\item \label{step1}\textbf{Initialization}\\
Choose an initial point $x_0\in\mathbb{R}^n$, a trial step length $\beta\in(0,1)$, and a positive $\nu\in(0,1)$.\\
Set the iteration number $k=0$.\\
Provide a value of the precision level $\epsilon>0$ for termination. \\
\item \label{step2}
 \textbf{Calculate the minimal set and the partition set at the $k$-th iteration}\\
 Compute $M_k=$ Min $(F(x_k),K)=\{r_1,r_2,\ldots,r_{w_k}\}$ and $w_k=\lvert$ Min $F(x_k),K \rvert$.\\
 Find 
 $P_k=P_{x_k}=I_{r_1}\times I_{r_2}\times\cdots\times  I_{r_{w_k}}$,~$p_k=\lvert P_{x_{k}}\rvert$, and $P_{x_k}=\{a_1,a_2,\ldots,a_{p_k}\}$,\\
 and for each $i\in [p_k],~a_i=(a^{1}_i,a^{2}_i,\ldots,a^{w_{k}}_i)\in P_{x_{k}},~a^j_{i}\in I_{r^{x_k}_j},~j\in[w_k].$\\
\item \label{step3}
     \textbf{Computation of a descent direction}\\
     Find $(a^k,u_k)\in \underset{(a,u)\in P_k\times\mathbb{R}^n}{\text{ argmin}} \xi_{x_k}(a,u)$, where $\xi_{x_k}:P_{x_k}\times\mathbb{R}^n\to \mathbb{R}$ is given by 
     $\xi_{x_k}(a,u)=\underset{j\in [w_k]}{\max} ~~ \Psi_e\left(\nabla f^{a_{j}}(x_k)^\top u_k+\tfrac{1}{2}u_k^\top\nabla^2 f^{a_{j}}(x_k) u_k\right).$ \\
\item \label{step4} \textbf{Stopping criterion}\\
  If $\lVert u_k\rVert<\epsilon$, stop. Otherwise, go to Step \ref{step5}.\\
\item \label{step5} \textbf{Compute a step length}\\
\textcolor{black}{Find the step length $t_k$ by the following relation:} 
\[\resizebox{0.9\textwidth}{!}{$
t_k= \underset{q\in\mathbb{N}\cup\{0\}}{\max}\left\{\nu^q~:f^{a_{j}^{k}}(x_k+\nu^q u_k)\preceq f^{a_{j}^{k}}(x_k)+\beta \nu^q \nabla f^{a_{j}^{k}}(x_k)^\top u_k ~\text{for all}~j\in[w_k]\right\}.
$}
\]
\item \label{step6} \textbf{Update the iterate} \\
Update $x_{k+1}\leftarrow x_k+t_k u_k$ and $k \leftarrow k+1$, and go to Step \ref{step2}. 
\end{enumerate}
\end{algorithm}

\begin{remark}
It is to be noted that for $p=1$ in Algorithm \ref{algo1}, that is, for $F(x)=\{f^1(x)\}$, the Step \ref{step3} of Algorithm \ref{algo1} reduces to finding $u_k$ such that
\[
u_k= \underset{u\in\mathbb{R}^n}{\text{ argmin}} ~~\Psi_e\left(\nabla f^{1}(x_k)^\top u_k+\tfrac{1}{2}u_k^\top\nabla^2 f^{1}(x_k) u_k\right).
\] 
\textcolor{black}{In this case, the proposed Algorithm \ref{algo1} reduces to the method given in \cite{drummond2014quadratically} and \cite{fliege2009newton}. In \cite{drummond2014quadratically}, Drummond-Svaiter functional and the support of a generator of the dual cone $K^*$ have been used, and we have used the Gerstewitz function $\Psi_e$.\\   
\\
The proposed Algorithm \ref{algo1} extends the approach proposed in \cite{drummond2014quadratically,fliege2009newton} to the \eqref{sp_equation}, and are found to be equivalent in case of vector optimization problems. It has been proved that the class of Gerstewitz functional $\Psi_e$ is a specific instance of other methods discussed in \cite{bouza2019unified}.}
\end{remark}

\subsection{Convergence analysis}\label{conv}
In this section, first, we show that Algorithm \ref{algo1} is well-defined. After that, we prove the convergence of Algorithm \ref{algo1}.\\

The well-definedness of Algorithm \ref{algo1} essentially depends on the following two points: 
\begin{enumerate}
    \item[(i)] Existence of $(a^k,u_k)$ in Step \ref{step3}, which is assured by the discussion in the paragraph after Definition \ref{stationary_dfn}. 
    \item[(ii)] Existence of step length $t_k$ in Step \ref{step5}, which is assured by the result in Proposition \ref{armijo}. 
\end{enumerate} 
Therefore, Algorithm \ref{algo1} is well-defined.\\

Next, we characterize the stationary points of (\ref{sp_equation}) in terms of the functions $\xi_x$ and $\Phi$ as defined in (\ref{function}) and (\ref{functionn}), respectively.


\begin{theorem}\label{critical}
Let us consider the functions $\xi_x$ and $\Phi$ as given in (\ref{function}) and (\ref{functionn}), respectively. Let $(\bar{a},\bar{u})\in P_x\times\mathbb{R}^n$ be such that $\Phi(\bar{x})=\xi_{\bar{x}}(\bar{a},\bar{u})$. Then, the following conditions are equivalent:
\begin{enumerate}
    \item[(i)] \label{c_0}The point $\bar{x}$ is a nonstationary point of (\ref{sp_equation}).
    \item[(ii)] \label{c_1} $\Phi(\bar{x})<0$. 
    \item[(iii)] \label{c_2} $\bar{u}\not=0.$
\end{enumerate}
\end{theorem}

\begin{proof}
(i)$\implies$(ii). Let us assume that the point $\bar{x}$ is a nonstationary point of (\ref{sp_equation}). Then, in view of (\ref{stationary_1}), there exists an $\widetilde{a}=(\widetilde{a}_1,\widetilde{a}_2,\ldots,\widetilde{a}_{\bar{w}})\in P_{\bar{x}}$ and $\tilde{u} \in \mathbb{R}^n$ such that 
\[
\Psi_e(\nabla f^{\widetilde{a}_j}(\bar{x})^\top \widetilde{u})
<0.
\]
Thus, in view of the above relation, we have  
\allowdisplaybreaks
\begin{eqnarray}\label{critical_0}
\Phi(\bar{x})&=&\underset{(a,u)\in P_{\bar{x}}\times\mathbb{R}^n}{\min}~\xi_{\bar{x}}(a,u)\nonumber\\
&\leq& \xi_{\bar{x}}(\widetilde{a},t\widetilde{u})\text{ for any }t>0, \widetilde{u}\in \mathbb{R}^n, \text{ and } a \in P_{\bar x} \nonumber \\
&=&\underset{j\in[w(\bar{x})]}{\max} \Psi_{e}\left(\nabla f^{\widetilde{a}_j}(\bar{x})^{\top}t\widetilde{u}+\tfrac{1}{2}t\widetilde{u}^\top\nabla^2 f^{\widetilde{a}_j}(\bar{x})t\widetilde{u}\right)  \nonumber\\
&=&t\underset{j\in[w(\bar{x})]}{\max} \Psi_{e}\left(\nabla f^{\widetilde{a}_j}(\bar{x})^{\top}\widetilde{u} + \tfrac{t}{2}\widetilde{u}^\top\nabla^2 f^{\widetilde{a}_j}(\bar{x})\widetilde{u}\right) \text{ by Proposition \ref{gerstewitz} (ii)}\nonumber\\
&\leq&t\underset{j\in[w(\bar{x})]}{\max}\left\{\Psi_{e}\left(\nabla f^{\widetilde{a}_j}(\bar{x})^{\top}\widetilde{u}\right)+\tfrac{t}{2}\Psi_e\left(\widetilde{u}^\top\nabla^2 f^{\widetilde{a}_j}(\bar{x})\widetilde{u}\right)\right\} \text{ by Proposition  \ref{gerstewitz} (i)\&(ii)}\nonumber\\
&\leq&t\left\{\underset{j\in[w(\bar{x})]}{\max}\{\Psi_{e}(\nabla f^{\widetilde{a}_j}(\bar{x})^{\top}\widetilde{u})\}+\tfrac{t}{2}\underset{j\in[w(\bar{x})]}{\max}\{\Psi_e(\widetilde{u}^\top\nabla^2 f^{\widetilde{a}_j}(\bar{x})\widetilde{u})\}\right\}. 
\end{eqnarray} 
Choosing any $t$ such that $0<t<\left(\frac{-2}{\underset{j\in [w(x)]}{\max}\{\Psi_e(\widetilde{u}^\top\nabla^2 f^{\widetilde{a}_j}(\bar{x})\widetilde{u})\}}\right)\left(\underset{j\in [w(x)]}{\max}\{\Psi_e(\nabla f^{\widetilde{a}_j}(\bar{x})^{\top}\widetilde{u})\}\right)$, we obtain from (\ref{critical_0}) that 
\[
\Phi(\bar{x})<t\left\{\underset{j\in[w(\bar{x})]}{\max}\{\Psi_{e}(\nabla f^{\widetilde{a}_j}(\bar{x})^{\top}\widetilde{u})\}-\underset{j\in[w(\bar{x})]}{\max}\{\Psi_{e}(\nabla f^{\widetilde{a}_j}(\bar{x})^{\top}\widetilde{u})\}\right\}=0.
\]
 {}
 \\
    (ii)$\implies$(iii).
   It trivially follows from  (\ref{function_2}). 
   \\
    (iii)$\implies$(i).
   Let us assume contrarily that $\bar{x}$ is a stationary point of (\ref{sp_equation}) and $\bar{u}\not=0$. Then, in view of (\ref{stationary_2}), for $\bar{a}\in P_{\bar{x}}$, there exists  $\widetilde{j}\in[\bar{w}]$  such that
   \begin{eqnarray}\label{critical_1}
    \Psi_e(\nabla f^{{\bar{a}}_{\widetilde{j}}}(\bar{x})^\top \bar{u})\geq0. 
   \end{eqnarray}
   Note from Assumption \ref{assumption} that for any $a\in P_x$ and $x\in\mathbb{R}^n$, we have $\bar{u}^\top\nabla^2 f^{{\bar{a}}_{\widetilde{j}}}(\bar{x})\bar{u}>0$. Therefore, from (\ref{critical_1}) \textcolor{black}{with the help of  Proposition \ref{gerstewitz}(iv)}, we get 
   \allowdisplaybreaks
     \begin{eqnarray*}
     &&\textcolor{black}{\nabla f^{{\bar{a}}_{\widetilde{j}}}(\bar{x})^\top \bar{u} +\tfrac{1}{2}\bar{u} ^\top\nabla^2f^{{\bar{a}}_{\widetilde{j}}}(\bar{x})\bar{u}\ge 0}\\
     &\text{or,}& \Psi_e\left(\nabla f^{{\bar{a}}_{\widetilde{j}}}(\bar{x})^\top \bar{u} +\tfrac{1}{2}\bar{u} ^\top\nabla^2f^{{\bar{a}}_{\widetilde{j}}}(\bar{x})\bar{u} \right)\geq0\\
     &\text{or,}& \underset{j\in[w(\bar{x})]}{\max}\left\{\Psi_e(\nabla f^{{\bar{a}}_{j}}(\bar{x})^\top \bar{u} +\tfrac{1}{2}\bar{u} ^\top\nabla^2f^{{\bar{a}}_{j}}(\bar{x})\bar{u})\right\}\geq0\\
     &\text{or,}&\xi_{\bar{x}}(\bar{a},\bar{u})\geq0\\
     &\text{or,}&\Phi(\bar{x})=0 \text{ from (\ref{function_3})}\\
     &\text{or,}&\bar{u}=0 \text{ from (\ref{function_2})}, 
    \end{eqnarray*}  
which is a contradiction to the considered assumption. Thus, $\bar{x}$ is a nonstationary point of (\ref{sp_equation}). 
\end{proof}
\begin{remark}
    In view of (\ref{function_3}), (\ref{function_2}), and statements (i)--(iii) of Theorem \ref{critical}, we obtain that $\bar{x}$ is a stationary point of (\ref{sp_equation}) if and only if $\Phi(\bar{x})=0$ or $\bar{u}=0.$
\end{remark}

In the next theorem, we estimate an upper bound for the norm of Newton direction $u_k$, generated by Algorithm \ref{algo1}, for (\ref{sp_equation}).


\begin{theorem} \label{uk_bounded}
Let $\{x_k\}$ be a sequence of nonstationary points, $\{u_k\}$ be a sequence of directions generated by Algorithm \ref{algo1}, and  $\{x_k\}$ be convergent. Then, the sequence $\{u_k\}$ is bounded.    
\end{theorem}

\begin{proof}
Let $P_{x_k}$ be the partition set at $x_k$. Then, by Theorem \ref{critical}, there exists $a^k\in P_{x_k}$ such that  
\begin{eqnarray}\label{bounded_eq1} 
& & \underset{j\in [w(x_k)]}{\max}\left\{\Psi_{e}(\nabla f^{a^k_j}(x_k)^{\top}u_k+\tfrac{1}{2}u_k^\top\nabla^2 f^{a^k_j}(x_k)u_k)\right\} <  0\nonumber\\
& \implies & \Psi_{e}(\nabla f^{a^k_j}(x_k)^{\top}u_k + \tfrac{1}{2}u_k^\top\nabla^2 f^{a^k_j}(x_k)u_k) <  0 ~\forall j \in [w(x_k)] \nonumber \\ 
& \overset{\eqref{assumption_eigen_value}}{\implies} & \tfrac{1}{2} \rho_{a_j^k} \|u_k\|^2 < - \Psi_{e}(\nabla f^{a^k_j}(x_k)^{\top}u_k)\nonumber ~\forall j \in [w(x_k)] \\
&\implies& \tfrac{1}{2} \rho_{a_j^k} \|u_k\|^2 < \underset{j\in [w(x_k)]}{\max}\{\lvert\Psi_e(\nabla f^{a^k_j}(x_k)^{\top}u_k)\rvert\} \le L \underset{j\in [w(x_k)]}{\max}\lVert\nabla f^{a^k_j}(x_k)^{\top}u_k\rVert \nonumber\\
&&\text{ from Proposition \ref{gerstewitz}(iii) and }L\text{ is a  Lipschitz constant of } \Psi_e.
\end{eqnarray}
Note that for every $k\in\mathbb{N}$, $f^{a^k_j}$'s are twice continuously differentiable and the sequence $\{x_k\}$ is convergent. Therefore, there exists a positive constants $C$ such that 
\begin{equation}\label{bounded_eq2}
C=\underset{j\in [w(x_k)]}{\max}\lVert\nabla f^{a^k_j}(x_k) \rVert.   
\end{equation}
Let $\rho = \min\{\rho_1, \rho_2, \ldots, \rho_p\}$, where $\rho_1, \rho_2, \ldots, \rho_p$ are as given in \eqref{assumption_eigen_value}. Then, $\rho > 0$ and in view of (\ref{bounded_eq1}) and (\ref{bounded_eq2}), we observe that  
\[\rho \lVert u_k\rVert^2\leq 2CL\lVert u_k\rVert\implies \lVert u_k\rVert\leq \frac{2CL}{\rho}.\]
Thus, the sequence $\{u_k\}$ is bounded. 
\end{proof}

Next, to prove the convergence of the proposed Algorithm \ref{algo1}, we derive a result on the existence of a step size $t_k$ (in Step 6) at every iterate $x_k$ along the chosen (descent) direction $u_k$ of $F$ for the set optimization problem (\ref{sp_equation}) by Algorithm \ref{algo1}. 

\begin{proposition}\label{armijo}
Let $\beta\in(0,1)$ and $(\bar{a},\bar{u})\in P_{\bar{x}}\times\mathbb{R}^n$ be such that $\Phi(\bar{x})=\xi_{\bar{x}}(\bar{a},\bar{u})$ and assume that the point $\bar{x}$ is not a stationary point of (\ref{sp_equation}). Then,  there exists $\widetilde{t}>0$ such that for all $t\in(0,\widetilde{t}]$ and $j\in[\bar{w}]$, 
    \begin{equation}\label{armijo_1}
    f^{\bar{a}_j}(\bar{x}+t\bar{u})\preceq f^{\bar{a}_j}(\bar{x})+\beta t\nabla f^{\bar{a}_j} (\bar{x})^\top\bar{u}.
    \end{equation}
    Additionally, for all $t\in(0,\widetilde{t}]$ and $j\in[\bar{w}]$, we have  
    \begin{equation}\label{armijo_2}
     F(\bar{x}+t\bar{u})\preceq^{l} \{f^{\bar{a}_j}(\bar{x})+\beta t\nabla f^{\bar{a}_j} (\bar{x})^\top\bar{u}\}_{j\in [\bar{w}]}\prec^l F(\bar{x}).   
    \end{equation}
\end{proposition}

\begin{proof} 
If possible, let (\ref{armijo_1}) do not hold. Therefore, there exists a sequence $\{t_k\} \searrow 0$ and ${j'}\in[\bar{w}]$ such that
\begin{eqnarray}\label{armijo_3}
&&f^{\bar{a}_{j'}}(\bar{x}+t_k\bar{u})-f^{\bar{a}_{j'}}(\bar{x})-\beta t_k\nabla f^{\bar{a}_{j'}} (\bar{x})^\top\bar{u}\not\in -K\nonumber\\
&\implies&\underset{k\to0}{\lim}\frac{f^{\bar{a}_{j'}}(\bar{x}+t_k\bar{u})-f^{\bar{a}_{j'}}(\bar{x})}{t_k}-\beta\nabla f^{\bar{a}_{j'}} (\bar{x})^\top\bar{u}\not\in -K\nonumber\\
&\implies&(1-\beta)\nabla f^{\bar{a}_{j'}} (\bar{x})^\top\bar{u}\not\in -\text{int}(K)\nonumber\\
&\implies&\nabla f^{\bar{a}_{j'}} (\bar{x})^\top\bar{u}\not\in -\text{int}(K)\text{ since }\beta\in(0,1).
\end{eqnarray}
Note that $\bar{x}$ is not a stationary point of (\ref{sp_equation}) and $(\bar{a},\bar{u})\in P_{\bar{x}}\times\mathbb{R}^n$ is such that $\Phi(\bar{x})=\xi_{\bar{x}}(\bar{a},\bar{u})$. Therefore, in view of Theorem \ref{critical}, $\bar u \neq 0$ and 
\begin{eqnarray}\label{armijo_4}
&&\xi_{\bar{x}}(\bar{a},\bar{u}) = \Phi(\bar{x}) <0\nonumber\\
&\implies&\underset{j\in[\bar{w}]}{\max}~\Psi_e(\nabla f^{\bar{a}_j}(\bar{x})^\top \bar{u}+\tfrac{1}{2}\bar{u}^\top\nabla^2f^{\bar{a}_j}(\bar{x})\bar{u})<0\text{ since }\bar{u}\not=0\nonumber\\
&\implies&\Psi_e(\nabla f^{\bar{a}_j}(\bar{x})^\top \bar{u}+\tfrac{1}{2}\bar{u}^\top\nabla^2f^{\bar{a}_j}(\bar{x})\bar{u})<0\text{ for all }j\in[\bar{w}] \notag \\ 
&\overset{\eqref{assumption_eigen_value}}{\implies} & \Psi_e(\nabla f^{\bar{a}_j}(\bar{x})^\top \bar{u})+ \tfrac{1}{2} \rho_{\bar{a}_j}\|\bar{u}\|^2 < 0\text{ for all }j\in[\bar{w}] \notag \\ 
&\implies & \Psi_e(\nabla f^{\bar{a}_j}(\bar{x})^\top \bar{u}) < 0 \notag \\ 
& \implies & \nabla f^{\bar{a}_j}(\bar{x})^\top \bar{u}\in -\text{int}(K),  
\end{eqnarray}
which is contradictory to (\ref{armijo_3}). Therefore, for every $j\in[\bar{w}]$, relation (\ref{armijo_1}) holds.\\

Now, from (\ref{armijo_1}) and Proposition \ref{min_min}, we observe that
 for every $t\in(0,\bar{t}]$, 
\begin{eqnarray*}
 F(\bar{x})&\subseteq& \{f^{\bar{a}_j}(x)\}_{j\in[\bar{w}]}+K\\
&\subseteq& \{f^{\bar{a}_j}(\bar{x})+\beta t\nabla f^{\bar{a}_j} (\bar{x})^\top\bar{u}\}_{j\in[\bar{w}]}+K+\text{ int}(K) \\
&\subseteq& \{f^{\bar{a}_j}(\bar{x}+t \bar{u})\}_{j\in[\bar{w}]}+K+K+\text{ int}(K)\text{ from  Definition \ref{partial}}\\
&\subseteq& F(\bar{x}+t \bar{u})+\text{ int}(K),
\end{eqnarray*}
which implies that for every $t\in(0,\bar{t}]$, we have $F(\bar{x}+t\bar{u})\prec^lF(\bar{x})$.
\end{proof}

Now, we present the main theorem of the paper that proves the convergence of the proposed Algorithm \ref{algo1}.

\begin{theorem}\label{convergence}
Let $\{x_k\}$ be a sequence of nonstationary points generated by Algorithm \ref{algo1} and $\bar{x}$ be one of its accumulation points. Additionally, assume that $\bar{x}$ is a regular point of $F$. Then, $\bar{x}$ is a stationary point of (\ref{sp_equation}).     
\end{theorem}

\begin{proof}
\textcolor{black}{Let $\{x_k\}$ be a sequence of nonstationary points and $\bar{x}$ be  an accumulation point of $\{x_k\}$.} We prove that $\bar{x}$ is stationary. Towards this, define a function $\varsigma:\mathcal{P}(\mathbb{R}^m)\to \mathbb{R}\cup\{-\infty\}$ by
 \[
 \varsigma(A)=\underset{z\in A}{\inf}~\Psi_e(z).
 \]
By Proposition \ref{gerstewitz}(iv), the function $\varsigma$ is monotonic with respect to the preorder $\preceq^l$, i.e., for all $A,B\in\mathcal{P}(\mathbb{R}^m)$, we have 
     \begin{equation}\label{pre-order-of-zeta}
     A\preceq^l B\implies \varsigma(A)\leq \varsigma(B).
     \end{equation}
     Now in view of (\ref{armijo_2}) of Proposition \ref{armijo}, for every $k=0,1,2,\ldots$, we obtain
     \allowdisplaybreaks
     \begin{eqnarray}\label{convergence_1}
      &&\varsigma(F(x_{k+1}))\nonumber\\
      =&&\varsigma(F(x_k+t_ku_k))\nonumber\\
       \leq&& \underset{j\in[w_k]}{\min}\left\{ \Psi_e\left(f^{a^k_j}(x_k)+\beta t_k\nabla f^{a^k_j} (x_k)^\top u_k\right)\right\}  \nonumber\\
       \leq&& \underset{j\in[w_k]}{\min}\left\{ \Psi_e\left(f^{a^k_j}(x_k)+\beta t_k\left(\nabla f^{a^k_j} (x_k)^\top u_k+\tfrac{1}{2}u_{k}^\top\nabla^2f^{a_{k,j}}(x_k)u_k\right)\right)\right\} \nonumber\\
       &&\text{ from Proposition \ref{gerstewitz}(iv) and }u_{k}^\top\nabla^2f^{a_{k,j}}(x_k)u_k \succ 0 \nonumber\\
       \leq&& \underset{j\in[w_k]}{\min}\left\{ \Psi_e\left(f^{a^k_j}(x_k)\right)+\beta t_k\Psi_e\left(\nabla f^{a^k_j} (x_k)^\top u_k+\tfrac{1}{2}u_{k}^\top\nabla^2f^{a^k_j}(x_k)u_k\right)\right\}\nonumber\\
       &&\text{ from Proposition \ref{gerstewitz}(i)}\nonumber\\
       \leq&& \underset{j\in[w_k]}{\min}\left\{ \Psi_e\left(f^{a^k_j}(x_k)\right)+\beta t_k\underset{j'\in[w_k]}{\max}\Psi_e\left(\nabla f^{a^k_{j'}} (x_k)^\top u_k+\tfrac{1}{2}u_{k}^\top\nabla^2f^{a^k_{j'}}(x_k)u_k\right)\right\}\nonumber\\
       \leq&& \underset{j\in[w_k]}{\min}\Psi_e\left(f^{a^k_j}(x_k)\right)+\beta t_k\underset{j\in[w_k]}{\max}\left\{\Psi_e\left(\nabla f^{a^k_{j}} (x_k)^\top u_k+\tfrac{1}{2}u_{k}^\top\nabla^2f^{a^k_{j}}(x_k)u_k\right)\right\}\nonumber\\
       =&&\varsigma(F(x_k))+\beta t_k \Phi(x_k).
     \end{eqnarray}
     Therefore, we get 
     \begin{equation}\label{aux_dg_1}
     -\beta t_k\underset{j\in [w_k]}{\max} \left\{\Psi_e\left(\nabla f^{a^k_j} (x_k)^\top u_k+\tfrac{1}{2}u_{k}^\top\nabla^2f^{a^k_j}(x_k)u_k\right)\right\}\leq\varsigma(F(x_k))-\varsigma(F(x_{k+1})).
     \end{equation}
     On adding the above relation for $k=0,1,\ldots,\kappa$, we obtain
     \begin{align}\label{convergence_2}
         -\beta\sum_{k=0}^{\kappa} t_k\underset{j\in [w_k]}{\max} \left\{\Psi_e(\nabla f^{a^k_j}(x_k)^\top u_k+\tfrac{1}{2}u_{k}^\top\nabla^2f^{a^k_j}(x_k)u_k\right\}\leq\varsigma(F(x_0))-\varsigma(F(x_{\kappa+1})) .
     \end{align}
\textcolor{black}{Since $\bar{x}$ is an accumulation point of the sequence $\{x_k\}$, we can find a subsequence $\mathcal{K}\in\mathbb{N}$ such that
\[ 
    \{x_k\}_{k\in\mathcal{K}}\to\bar{x},~\{t_k\}_{k\in\mathcal{K}}\to \bar{t}, ~\text{ and }~\{u_k\}_{k\in\mathcal{K}} \to \bar{u}.
      \]
In view of \eqref{function_3}, \textcolor{black}{the function $\varsigma$ in \eqref{convergence_1}  is monotonic.}} Therefore, from  \eqref{convergence_2}, taking $\kappa \to\infty$ we deduce that 
     \begin{eqnarray}
     -\beta\underset{\kappa\to\infty}{\lim}\sum_{k=0}^{\kappa} t_k\underset{j\in [w_k]}{\max} \left\{\Psi_e\left(\nabla f^{a^k_j}(x_k)^\top u_k+\tfrac{1}{2}u_{k}^\top\nabla^2f^{a^k_j}(x_k)u_k\right)\right\}\leq + \infty.\label{convergence_6}
     \end{eqnarray}
\textcolor{black}{Now, note that if $x_k$ is not a stationary point, then in view of (\ref{stationary_1}),
     for every $a^k\in P_{x_k},u_k\in\mathbb{R}^n,j\in[w(x_k)]$, and $k\in\mathcal{K}$,} we get $\nabla f^{a^k_j}(x_k)^\top u_k\in-\text{ int}(K)$. On proceeding in a similar manner to (\ref{critical_0}), we can find $t_k>0$ such that  
     \[0<t_k<\left(\frac{-2}{\underset{j\in [w_k]}{\max}\{\Psi_e( {u}_{k}^\top \nabla^2 f^{a^k_j}(x_k) {u}_{k})}\right)\left(\underset{j\in [w_k]}{\max}\{\Psi_e(\nabla f^{a^k_j}(x_k)^{\top} {u}_{k})\}\right).\]
     In view of the above chosen \textcolor{black}{$t_k>0,k\in\mathcal{K}$}, we conclude that 
     \begin{eqnarray*}
     &&t_k\underset{j\in [w_k]}{\max} \left\{\Psi_e\left(\nabla f^{a^k_j} (x_k)^\top u_k+\tfrac{1}{2}u_{k}^\top\nabla^2f^{a^k_j}(x_k)u_k\right)\right\}\\
     &\leq & t_k\left\{\underset{j\in[w(\bar{x})]}{\max}\{\Psi_{e}(\nabla f^{a^k_j}(\bar{x})^{\top}\bar{u})\}-\underset{j\in[w(\bar{x})]}{\max}\{\Psi_{e}(\nabla f^{a^k_j}(\bar{x})^{\top}\bar{u})\}\right\}\\
     &=&0.
      \end{eqnarray*}
      Therefore, we get
       \begin{eqnarray}\label{convergence_7}
     &&-t_k\underset{j\in [w_k]}{\max} \left\{\Psi_e(\nabla f^{a^k_j} (x_k)^\top u_k+\tfrac{1}{2}u_{k}^\top\nabla^2f^{a^k_j}(x_k)u_k\right\}\geq0.
     \end{eqnarray}
     On combining (\ref{convergence_6}) and (\ref{convergence_7}), and taking limit $k\to\infty$, we have 
     \begin{eqnarray}\label{convergence_dg}
     && 0 \le - \sum_{k = 0}^\infty t_k\underset{j\in [w_k]}{\max} \left\{\Psi_e(\nabla f^{a^k_j} (x_k)^\top u_k+\tfrac{1}{2}u_{k}^\top\nabla^2f^{a^k_j}(x_k)u_k\right\} \le +\infty.
      \end{eqnarray}
Hence, we obtain     
     \begin{eqnarray}\label{convergence_3}
     &&\underset{\textcolor{black}{k\to\infty,~k\in\mathcal{K}}}{\lim} ~ t_k\underset{j\in [w_k]}{\max} \left\{\Psi_e(\nabla f^{a^k_j} (x_k)^\top u_k+\tfrac{1}{2}u_{k}^\top\nabla^2f^{a^k_j}(x_k)u_k)\right\}=0.
      \end{eqnarray}
      As the number of points in $[p]$ is finite, and $\bar{x}$ is a regular point of $F$, therefore in view of Lemma \ref{regular}, \textcolor{black}{for all $k\in\mathcal{K},u\in\mathbb{R}^n$}, we have $~w_k=\bar{w},P_{x_k}=\bar{P},a^k=\bar{a}\in\bar{P}$ and
      \begin{eqnarray}\label{convergence_4}
        &&\Phi(x_k)=\xi_{x_k}(\bar{a},u_k)\leq \xi_{x_k}(a,u)\nonumber\\
        &\text{and}&\xi_{\bar{x}}(\bar{a},\bar{u})\leq \xi_{\bar{x}}(a,u)\textcolor{black}{\text{ on taking limit }k{\to} \infty,k\in\mathcal{K}}.
      \end{eqnarray}
      Now, we analyze the following two cases:
      \begin{enumerate}
          \item[(i)] Let $\bar{t}>0$. In view of (\ref{convergence_3}) and \textcolor{black}{for all $k\in\mathcal{K}$ such that $w_k=\bar{w},P_{x_k}=\bar{P},a^k=\bar{a}$,} we have
          \begin{eqnarray}\label{convergence_5}
           &&\underset{k\overset{\mathcal{K}}{\to}+\infty}{\lim}\underset{j\in [\bar{w}]}{\max} \left\{\Psi_e(\nabla f^{\bar{a}_{j}} (x_k)^\top u_k+\tfrac{1}{2}u_{k}^\top\nabla^2f^{\bar{a}_{j}}(x_k)u_k\right\}=0\nonumber\\
    &\implies&\underset{k\overset{\mathcal{K}}{\to}+\infty}{\lim}\Phi(x_k)=0, \text{ i.e., } \Phi(\bar x) = 0. 
          \end{eqnarray}
          Thus, by Theorem \ref{critical}, $\bar{u}=0$. Hence, $\bar{x}$ is a stationary point of (\ref{sp_equation}).
          \item[(ii)] Let $\bar{t}=0$. Fix any $\kappa \in \mathbb{N}$. Since $t_k\overset{\mathcal{K}}{\to} \bar t = 0$, large enough $\nu^\kappa$ does not satisfy Armijo condition in Step \ref{step5} of Algorithm \ref{algo1}. Therefore for all $k\in\mathcal{K}$ such that $w_k=\bar{w},~P_{x_k}=\bar{P}$, and $a^k=\bar{a}$, there exists $\bar{j}\in[\bar{w}]$ such that 
          \begin{eqnarray*}
          &&f^{\bar{a}_{\bar{j}}}(x_k+\nu^\kappa u_k)\npreceq f^{\bar{a}_{\bar{j}}}(x_k)+\beta \nu^\kappa \nabla f^{\bar{a}_{\bar{j}}} (x_k)^\top u_k\\
          &\implies&\frac{f^{\bar{a}_{\bar{j}}}(x_k+\nu^\kappa u_k)-f^{\bar{a}_{\bar{j}}}(x_k)}{\nu^\kappa}-\beta \nabla f^{\bar{a}_{\bar{j}}} (x_k)^\top u_k\not\in -K\\
          &\implies&\frac{f^{\bar{a}_{\bar{j}}}(\bar{x}+\nu^\kappa \bar{u})-f^{\bar{a}_{\bar{j}}}(\bar{x})}{\nu^\kappa}-\beta \nabla f^{\bar{a}_{\bar{j}}} (\bar{x})^\top \bar{u}\not\in -\text{ int}(K)~\text{taking }k\overset{\mathcal{K}}{\to}+\infty\\
          &\implies& (1-\beta)\nabla f^{\bar{a}_{\bar{j}}}(\bar{x})^\top\bar{u}\not\in-\text{ int}(K)~\text{taking limit }k\to+\infty\\
          &\implies&\nabla f^{\bar{a}_{\bar{j}}}(\bar{x})^\top\bar{u}\not\in-\text{ int}(K)~\text{since }(1-\beta)\in(0,1).
          \end{eqnarray*}
          Therefore, from (v) of Proposition \ref{gerstewitz}, we have 
          \begin{eqnarray*}
              &&\Psi_e(\nabla f^{\bar{a}_{\bar{j}}}(\bar{x})^\top\bar{u})\geq 0 \\
              &\text{or, }&\Psi_e(\nabla f^{\bar{a}_{\bar{j}}}(\bar{x})^\top\bar{u}+\tfrac{1}{2}\bar{u}^\top\nabla^2 f^{\bar{a}_{\bar{j}}}(\bar{x})\bar{u})\geq0 \\ 
              && \text{ from (iv) of Proposition \ref{gerstewitz} and } \bar{u}^\top\nabla^2 f^{\bar{a}_{\bar{j}}}(\bar{x})\bar{u} \succ 0 \\
              &\text{or, }&0\leq \Psi_e(\nabla f^{\bar{a}_{\bar{j}}}(\bar{x})^\top\bar{u}+\tfrac{1}{2}\bar{u}^\top\nabla^2 f^{\bar{a}_{\bar{j}}}(\bar{x})\bar{u})\\
              &\text{or, }&0\leq \xi_{\bar{x}}(\bar{a},\bar{u})=\underset{(a,u)\in P_{x}\times \mathbb{R}^n}{\min}\xi_x(a,u)=\Phi(\bar{x})\leq0\text{ from (\ref{function_3})}.
          \end{eqnarray*}
          Thus, from Theorems \ref{critical}, we conclude that $\bar{x}$ is a stationary point of (\ref{sp_equation}). 
      \end{enumerate}
\end{proof}

\textcolor{black}{Next, we analyze the rate of convergence of the proposed Algorithm \ref{algo1}. It is found in the following Theorem \ref{superlinear} that under suitable assumptions, the step length $t_k$ is eventually $1$, and a subsequence of the generated sequence $\{x_k\}$ by Algorithm \ref{algo1} converges superlinearly to a locally efficient solution. Towards this, at first, we recall the following result.}
\begin{lemma}\emph{(See \cite{fliege2009newton}).}\label{lipschitz}
Let $V$ be a nonempty convex subset of $\mathbb{R}^n$, and $\epsilon>0$ and $\delta>0$ be such that for any $x,y\in V$ with $\lVert y-x\rVert<\delta$,  
\begin{equation}\label{li_1}
\lVert\nabla^2 f^{a_j}(y)-\nabla^2 f^{a_j}(x)\rVert<\epsilon \text{ for all } j\in [w(x)].
\end{equation} 
Then, for every $j\in [w(x)]$, we have
\begin{equation}\label{li_2} 
\lVert \nabla f^{a_j}(y)-[\nabla f^{a_j}(x)+\nabla^2 f^{a_j}(x)^\top (y-x)]\rVert<\epsilon \lVert y-x\rVert.    
\end{equation}
If $\nabla^2 f^{a_j}$ is Lipschitz continuous on $V$ with constant $\widetilde{L}$ for all $j\in[w(x)]$, then
\begin{equation}\label{li_3} 
\lVert \nabla f^{a_j}(y)-[\nabla f^{a_j}(x)+\nabla^2 f^{a_j}(x)^\top (y-x)]\rVert<\tfrac{\widetilde{L}}{2}\lVert y-x\rVert^2 \text{ for all } j\in w[(x)].    
\end{equation}
\end{lemma}

\begin{theorem}\emph{(Superlinear convergence).}\label{superlinear}
Let $\{x_k\}$ be a sequence of nonstationary points generated by Algorithm \ref{algo1} and $\bar{x}$ be one of its accumulation points. Additionally, assume that $\bar{x}$ is a regular point of $F$, and  there exists a nonempty convex set $V \subseteq \mathbb{R}^n$, and $\gamma >0,\delta>0,$ and $\epsilon>0$ for which the following conditions hold: 
\begin{enumerate}[(i)]
    \item\label{supe_1} $\nabla^2 f^{a_j}(x)\leq \gamma I$ for all $j\in [w(x)]$,
    \item\label{supe_2} $\lVert\nabla^2 f^{a_j}(x)-\nabla^2f^{a_j}(y) \rVert<\epsilon $ for all $x,y\in V$ with $\lVert x-y \rVert<\delta$, and 
    \item \label{supe_3} $\tfrac{\epsilon}{\gamma }\leq 1-2\beta$.
\end{enumerate}
Then, for sufficiently large $k\in\mathbb{N}$, $t_k=1$ holds and there exists a subsequence of $\{x_k\}$ that converges superlinearly to $\bar{x}$.
\end{theorem}

\begin{proof}
From Theorem \ref{convergence}, \textcolor{black}{we obtain that the sequence $\{x_k\}$ converges to $\bar{x}$ and $\bar{x}$ is a stationary point of \eqref{sp_equation}.} \\ 
To prove superlinear convergence, note that each $f^{a_j}$ is twice continuously differentiable. Therefore, for any $\epsilon>0$, there exists $\delta_{\epsilon}>0$ such that
\[
\lVert \nabla^2 f^{a_j}(x)-\nabla^2 f^{a_j}(\bar x)\rVert<\epsilon \text{ for all }x\in \mathcal{B}(\bar{x},\delta_\epsilon)\subseteq V. 
\]
For $x_k\in V$, let $w_k = w(x_k)$. For any $\lambda \in \Delta_k = \{(\lambda_1, \lambda_2, \ldots, \lambda_{w_k}) \in \mathbb{R}^{w_k}_+ : \sum_{i = 1}^{w_k} \lambda_i = 1\},$ we define a function $\Theta_{\lambda}:V\times \mathbb{R}^n\to \mathbb{R}^m$ by 
 \[
 \Theta_{\lambda}(x_k,u)=\sum_{j=1}^{[w_k]}\lambda_j\nabla f^{a^k_j}(x_{k})^\top u + \tfrac{1}{2}\sum_{j=1}^{[w_k]}\lambda_ju^\top \nabla^2 f^{a^k_j}(x_{k}) u. 
 \]
Note that for any $a^k\in P_{x_k}$ and $j\in[w(x_k)]$, each $f^{a_j^k}$ is twice continuously differentiable and strongly convex function. Moreover, the set $P_{x_k}$ is finite. Therefore, the function $\Theta_{\lambda}(a,\cdot)$ is strongly convex in $\mathbb{R}^n$ and hence the function $\Theta_{\lambda}(a,\cdot)$ attains its minimum. Then, using Danskin's theorem (see Proposition 4.5.1, pp. 245--247 in \cite{bertsekas2003convex}) and the first order necessary condition for a minimizer $u_k$ of $\Theta_{\lambda}(x_k,\cdot)$, we have 
\allowdisplaybreaks
\begin{eqnarray}
&&\sum_{j=1}^{[w_k]}\lambda_j\nabla f^{a^k_j}(x_{k})+\sum_{j=1}^{[w_k]}\lambda_j\nabla^2 f^{a^k_j}(x_{k}) u_k=0\label{super_2}\\
\implies&& u_k= -\left[\sum_{j=1}^{[w_k]}\lambda_j\nabla^2 f^{a^k_j}(x_{k})\right]^{-1}\sum_{j=1}^{[w_k]}\lambda_j\nabla f^{a^k_j}(x_{k})\nonumber\\
\implies&& u_k \leq -\tfrac{1}{\gamma } \sum_{j=1}^{[w_k]}\lambda_j\nabla f^{a^k_j}(x_{k})\text{ since }\nabla^2 f^{a^k_j}(x_k)\leq \gamma I\nonumber\\
\implies&& u_k\leq-\tfrac{1}{\gamma } \underset{\lambda\in \Delta_k}{\max}\sum_{j=1}^{[w_k]}\lambda_j\nabla f^{a^k_j}(x_{k})\label{super_1}.
\end{eqnarray}
As the sequence $\{x_k\}$ converges to $\bar{x}$, there exists $k_{\epsilon}\in\mathbb{N}$ such that for all $k\geq k_\epsilon$, we have $
x_k,x_k+u_k\in \mathcal{B}(\bar{x},\delta_{\epsilon})$. Now, using the second-order Taylor expansion at $x_k$ of $f^{a^k_j}$, we have
\begin{eqnarray*}
 f^{a^k_j}(x_k+u_k)&\leq& f^{a^k_j}(x_k) + \nabla f^{a^k_j}(x_k)^\top u_k+\tfrac{1}{2}u_k^\top \nabla^2 f^{a^k_j}(x_k)u_k +\tfrac{\epsilon}{2}\lVert u_k\rVert^2.
 \end{eqnarray*}
As $\max\{b_1,b_2,\ldots,b_{w_k}\}=\underset{\lambda\in \Delta_{w_k}}{\max}\sum_{i=1}^{w_k} \lambda_i b_i $ holds, we get
 \begin{eqnarray}\label{armijo_holds}
&&f^{a^k_j}(x_k+u_k)- f^{a^k_j}(x_k)\nonumber\\
&\leq& \nabla f^{a^k_j}(x_k)^\top u_k+\tfrac{1}{2}u_k^\top \nabla^2 f^{a^k_j}(x_k)u_k +\tfrac{\epsilon}{2}\lVert u_k\rVert^2\nonumber\\
 &\leq& \beta \nabla f^{a^k_j}(x_k)^\top u_k+(1-\beta)\nabla f^{a^k_j}(x_k)^\top u_k+\tfrac{(\gamma +\epsilon)}{2}\lVert u_k\rVert^2\text{ since }\nabla^2 f^{a^k_j}(x)\leq \gamma I\nonumber\\
  &\leq& \beta \nabla f^{a^k_j}(x_k)^\top u_k+(1-\beta)\underset{j\in[w(x_k)]}{\max}\{\nabla f^{a_j}(x_k)^\top u_k\}+\tfrac{(\gamma +\epsilon)}{2}\lVert u_k\rVert^2\nonumber\\
  &\leq& \beta \nabla f^{a^k_j}(x_k)^\top u_k+(1-\beta)\underset{\lambda\in\Delta_k}{\max}\left\{\sum_{j=1}^{[w(x_k)]}\lambda_j\nabla f^{a^k_j}(x_k)^\top u_k\right\}+\tfrac{(\gamma +\epsilon)}{2}\lVert u_k\rVert^2 \nonumber\\
 &\leq& \beta \nabla f^{a^k_j}(x_k)^\top u_k-\gamma (1-\beta)\lVert u_k\rVert^2 + \tfrac{(\gamma +\epsilon)}{2}\lVert u_k\rVert^2\text{ by (\ref{super_1})}\nonumber\\ 
 &\leq& \beta \nabla f^{a^k_j}(x_k)^\top u_k+\frac{\epsilon-\gamma (1-2\beta)}{2}\lVert u_k\rVert^2,
\end{eqnarray}
\textcolor{black}{where from assumption \eqref{supe_3}, we conclude that $\epsilon-\gamma (1-2\beta)\leq 0$, and hence $t_k=1$ holds.} Now, for $k\geq k_\epsilon$, $\lambda\in \Delta_k$, and $j\in[w(x_k)]$, we have
\begin{eqnarray}\label{super_3}
&&\left\lVert \sum_{j=1}^{[w(x_{k+1})]}\lambda_j\nabla f^{a^k_j}(x_{k+1}) \right\rVert\nonumber\\
&=&\left\lVert  \sum_{j=1}^{[w(x_{k+1})]}\lambda_j\nabla f^{a^k_j}(x_{k}+u_k) \right\rVert\nonumber\\
&\overset{\eqref{super_2}}{=}& \left\lVert  \sum_{j=1}^{[w(x_{k+1})]}\lambda_j\nabla f^{a^k_j}(x_{k}+u_k)-\left[\sum_{j=1}^{[w(x_k)]}\lambda_j\nabla f^{a^k_j}(x_{k})+\sum_{j=1}^{[w(x_k)]}\lambda_j\nabla^2 f^{a^k_j}(x_{k})^\top u_k\right]\right\rVert \nonumber\\
&\leq& \epsilon \lVert u_k\rVert \text{ by \eqref{li_2} of Lemma \ref{lipschitz}}. 
\end{eqnarray}
Combining assumption (i) and boundedness of $\{u_{k+1}\}$ (Theorem \ref{uk_bounded}), we observe that
\begin{eqnarray}
\tfrac{1}{2}u_k^\top\nabla^2 f^{a^k_j}(x_k)u_k\leq \tfrac{\gamma }{2}\lVert u_k\rVert^2 \text{ for any }j\in[w(x_k)].    
\end{eqnarray}
Therefore, incorporating the above relation in (\ref{bounded_eq1}), we get
\begin{align*}
\lVert u_{k+1}\rVert 
~\leq~ \tfrac{2L}{\gamma } \underset{j\in [w(x_k)]}{\max}\left\lVert\nabla f^{a^k_j}(x_{k+1})\right\lVert 
~\leq~  \tfrac{2L}{\gamma } \underset{\lambda\in\Delta_k}{\max}\left\lVert\sum_{j=1}^{[w_k]}\lambda_j\nabla f^{a^k_j}(x_{k+1})\right\lVert  
~\overset{\eqref{super_3}}{\leq} ~ \tfrac{2L\epsilon}{\gamma }\lVert u_k\rVert.    
\end{align*}
In view of the above relation, we have
\begin{align*}
\lVert x^{k+1}-x^{k+2}\rVert= \lVert u^{k+1}\rVert  \leq \tfrac{2L\epsilon}{\gamma } \lVert u_k\rVert =\tfrac{2L\epsilon}{\gamma } \lVert x^{k}-x^{k+1}\rVert, 
\end{align*}
and for any $k\geq1$ and $m\geq1$, we obtain
\begin{align}\label{superlinear_1}
\lVert x^{k+m}-x^{k+m+1}\rVert\leq& \left(\tfrac{2L\epsilon}{\gamma }\right)
\lVert x^{k+m-1}-x^{k+m}\rVert\nonumber\\
\leq& \left(\tfrac{2L\epsilon}{\gamma }\right)^2\lVert x^{k+m-1}-x^{k+m}\rVert\nonumber\\
\leq&\cdots\leq \left(\tfrac{2L\epsilon}{\gamma }\right)^m\lVert x^{k}-x^{k+1}\rVert. 
\end{align}
Now, we assume $0<\tau<1$ and define 
\[
\bar{\epsilon}=\min\left\{\gamma (1-2\beta),\tfrac{\tau}{1+2\tau}\left(\tfrac{\gamma }{2L\epsilon}\right)\right\}.
\]
If we take $\epsilon<\bar{\epsilon}$ and $k\geq k_\epsilon$, then by convergence of sequence $\{x_k\}$ and relation (\ref{superlinear_1}), we have
\begin{eqnarray*}
\lVert \bar{x}-x^{k+1}\rVert\leq \sum_{m=1}^{\infty}\lVert x^{k+m}-x^{k+m+1}\rVert &\leq& \sum_{m=1}^{\infty}\left(\tfrac{\tau}{1+2\tau}\right)^m\lVert x^{k}-x^{k+1}\rVert\\
&=&\tfrac{\tau}{1+\tau}\lVert x^{k}-x^{k+1}\rVert. 
\end{eqnarray*}
Therefore, we obtain
\begin{align*}
&\lVert \bar{x}-x^{k}\rVert \geq \lVert x^{k}-x^{k+1}\rVert-\lVert x^{k+1}-\bar{x}\rVert \geq \tfrac{1}{1+\tau}\lVert x^{k}-x^{k+1}\rVert.
\end{align*}
Hence, we can conclude that if $\epsilon<\bar{\epsilon}$ and $k\geq k_\epsilon$, then $\frac{\lVert \bar{x}-x^{k+1}\rVert}{\lVert \bar{x}-x^{k}\rVert}\leq \tau.$
\end{proof}

In the conventional Newton method, it is well known that once the initial point is chosen from a close vicinity of the optimal solution, then the entire sequence of iterates resides in the same vicinity of the optimal solution and converges to it. The next result shows a similar observation in the proposed Newton method for \eqref{sp_equation}.

\begin{proposition}\label{corollary_1}
Let $\bar{x}$ be a stationary point of $F$. Then, there exists $\rho>0$ and a nonempty convex set $V\subseteq\mathbb{R}^n$ such that all the assumptions in Theorem \ref{superlinear} are satisfied for any chosen initial point $x_0\in\mathcal{B}(\bar{x},\rho)$. Moreover, if the initial point $x_0$ belongs to a compact level set of $F$ and $\{x_k\}$ is the sequence generated from $x_0$, then the sequence $\{x_k\}$ converges to some stationary point $\widetilde{x}$ of \eqref{sp_equation}.   
\end{proposition}

\begin{proof}
Given that $\bar{x}$ is a stationary point of $F$. Take $R>0$ such that $\mathcal{B}(\bar{x},R)\subseteq V$. Then, there exist $p> 0$ and $\gamma >0$ such that $ \nabla^2 f^{a_j}(x)\leq \gamma I$ for all $x\in V$ and $j\in [w(x)]$, which is the assumption (i) of Theorem \ref{superlinear}. \\ 
Now, take $\epsilon>0$ such that $\epsilon\leq \gamma (1-2\beta)$. Then, there exists $\delta>0$ such that for all $x,y\in V$ and $j\in[w(\bar{x})]$, 
\[
\lVert\nabla^2 f^{a_j}(x)-\nabla^2f^{a_j}(y) \rVert<\epsilon  \text{ with }\lVert x-y \rVert<\delta,
\]
which is the assumption (ii) of Theorem \ref{superlinear}. Note that $\bar{x}$ is a stationary point of $F$. Then, in view of Theorem \ref{critical}, we have
$\Phi(\bar{x})=0.$ \\ 
From Theorem \ref{continuity}, the function $\Phi$ is continuous. Thus, there exists $\rho \in (0,\tfrac{R}{2})$ such that for any $x\in\mathcal{B}(\bar{x},\rho)$, we have 
\[
\lvert \Phi(x)-\Phi(\bar x)\rvert \leq \tfrac{p}{2}\left[ \min\left\{\delta,\left(\tfrac{R}{2}\right)\left(1-\tfrac{\epsilon}{p}\right)\right\}\right]^2.
\]
Given that $\{x_k\}$ is the sequence generated from $x_0$. Then, in view of Proposition \ref{armijo}, the sequence $\{F(x_k)\}$ is  set-less decreasing. Moreover, the sequence $\{x_k\}$ is bounded and $f^{a_j}$'s are twice continuously differentiable. Therefore, the sequence $\{f^{a_j}(x_k)\}$ is also bounded. Thus, in view of Proposition \ref{armijo} and Step \ref{step5} of Algorithm \ref{algo1}, we get 
\[
\underset{k\to \infty}{\lim}t_k\Phi(x_k)=0.
\]
Hence, we can observe that there exists a subsequence $\{x_{k_j}\}$ of $\{x_k\}$ such that $x_{k_j} \to \widetilde{x}$, for some stationary point $\widetilde{x}\in V$. If this claim is not true, i.e., $\widetilde{x}$ is not stationary, then in view of Proposition \ref{armijo} and Step \ref{step5} of Algorithm \ref{algo1}, we have 
\begin{eqnarray*}
&&\underset{j\to\infty}{\liminf}~t_{k_j}>0 ~\text{ or }~ \underset{j\to\infty}{\lim} \Phi(x_{k_j})=0.
\end{eqnarray*}
Therefore, by Proposition \ref{continuity}, we have $\Phi(\widetilde{x})=0$,
which is contradictory to the result in Theorem \ref{critical}. Thus, the point $\widetilde{x}$ is stationary. Additionally, in view of Proposition \ref{corollary_1}, we conclude that for $j$ large enough, $x_{k_j}$ belongs to the close neighborhood of and converges to stationary point $\widetilde{x}$.
\end{proof} 

Next, we analyze the quadratic convergence of the proposed Algorithm \ref{algo1}.

\begin{theorem}\emph{(Quadratic convergence).} \label{quadratic}
Let $\{x_k\}$ be a sequence of nonstationary points generated by Algorithm \ref{algo1} an initial point $x_0$ \textcolor{black}{belong} to a compact level set of $F$. Let there exist $\gamma >0$ for which $\nabla^2 f^{a_j}(x)\leq \gamma I$ for all $j\in [w(x)]$ and assume that for all $j\in [w(x)]$, $\nabla^2 f^{a_j}$ is Lipschitz continuous with a common Lipschitz constant $\widetilde{L}$. Assume that all the conditions of Theorem \ref{superlinear} hold. Then,  the sequence $\{x_k\}$ converges quadratically to a stationary point 
$\bar{x}$ of \eqref{sp_equation}.
\end{theorem}

\begin{proof}
\textcolor{black}{From Theorem \ref{convergence} and \ref{superlinear}, we observe that $\bar{x}$ is a stationary point of (\ref{sp_equation}) and $t_k=1$ for large enough $k$.} Since $f^{a_j}$'s are twice continuously differentiable, for any $\epsilon>0$, there exists $\delta_{\epsilon}>0$ such that for all $x,y\in \mathcal{B}(\bar{x},\delta_\epsilon)$, we have
\[
\lVert \nabla^2 f^{a_j}(x)-\nabla^2 f^{a_j}(y)\rVert<\epsilon.
\] 
Now proceeding in similar steps as in Theorem \ref{superlinear} and in view of (iii) of Lemma \ref{lipschitz}, we conclude that for $k\geq k_\epsilon$, we have
\allowdisplaybreaks
\begin{eqnarray}\label{quadratic_1}
&&\left\lVert \sum_{j=1}^{[w(x_{k+1})]}\lambda_j\nabla f^{a^k_j}(x_{k+1}) \right\rVert \nonumber\\
&=&\left\lVert  \sum_{j=1}^{[w(x_{k+1})]}\lambda_j\nabla f^{a^k_j}(x_{k}+u_k) \right\rVert\nonumber\\
&\overset{(\ref{super_2})}{=}& \left\lVert  \sum_{j=1}^{[w(x_{k+1})]}\lambda_j\nabla f^{a^k_j}(x_{k}+u_k)-\left[\sum_{j=1}^{[w(x_k)]}\lambda_j\nabla f^{a^k_j}(x_{k})+\sum_{j=1}^{[w(x_k)]}\lambda_j\nabla^2 f^{a^k_j}(x_{k})^\top u_k\right]\right\rVert \nonumber\\
&\leq& \tfrac{\widetilde{L}}{2} \lVert u_k\rVert^2 \text{ by (iii) of Lemma \ref{lipschitz}}. 
\end{eqnarray}
 Note that $\max\{b_1,b_2,\ldots,b_{w(x)}\}=\underset{\lambda\in \Delta_{w(x)}}{\max}\sum_{i=1}^{w(x)} \lambda_i b_i $ holds, where 
 $$\Delta_{w(x)} = \{(\lambda_1, \lambda_2, \ldots, \lambda_{w(x)}) \in \mathbb{R}^{w(x)}_+: \sum_{i = 1}^{w(x)} \lambda_i = 1\}.$$
 On combining the relation (\ref{quadratic_1}) with assumption (i) and boundedness of $\{u_{k+1}\}$ (Theorem \ref{uk_bounded}), we get
\begin{align}\label{quadratic_2}
&\lVert u_{k+1}\rVert\leq\tfrac{2L}{\gamma } \underset{j\in [w(x_{k+1})]}{\max}\left\lVert\nabla f^{a_{k,j}}(x_{k+1})\right\lVert\nonumber\\
\text{or, }&\lVert u_{k+1}\rVert\leq\tfrac{2L}{\gamma }\left \{\underset{\lambda\in\Delta_k}{\max}\left\lVert\sum_{j=1}^{[w(x_k)]}\nabla f^{a^k_j}(x_{k+1})\right\lVert\right\} \nonumber\\
\text{or, }&\lVert u_{k+1}\rVert\leq\tfrac{2L\widetilde{L}}{\gamma }\lVert u_k\rVert^2\text{ from (\ref{quadratic_1})}.  
\end{align}
Since from Theorem \ref{superlinear} the sequence $\{x_k\}$ converges superlinearly to $\bar{x}$ of (\ref{sp_equation}), then in view of (\ref{quadratic_1}) and (\ref{quadratic_2}) there exist $N$ such that for $k\geq N$, we have
\begin{eqnarray}\label{quadratic_eq_a}
\lVert \bar{x}-x^{k+1}\rVert\leq \tau \lVert \bar{x}-x^{k}\rVert\text{ for some }\tau\in(0,1).  
\end{eqnarray}
\textcolor{black}{Further, triangle inequality, for $k=l\geq N$, we obtain 
\begin{eqnarray}
&&\lVert x^{l}-x^{l+1} \rVert \le \lVert x^{l}-\bar{x}\rVert+\lVert \bar{x}-x^{l+1}\rVert \overset{\eqref{quadratic_eq_a}}{\le}   (1+\tau)\lVert \bar{x}-x^{l} \rVert \label{quadratic_3}\\
&\text{and }&\lVert x^{l}-x^{l+1} \rVert \ge \lVert x^{l}-\bar{x}\rVert-\lVert \bar{x}-x^{l+1}\rVert \overset{\eqref{quadratic_eq_a}}{\ge}(1-\tau)\lVert \bar{x}-x^{l} \rVert \label{quadratic_3b}.
\end{eqnarray}}
Finally, from relation (\ref{quadratic_3b}) for $l=k+1$, we conclude that
\begin{eqnarray*}
(1-\tau) \lVert \bar{x}-x^{k+1}\rVert\leq \lVert x^{k+1}-x^{k+2}\rVert&=&\lVert u_{k+1}\rVert\\
&\overset{\eqref{quadratic_2}}{\leq} &\tfrac{2L\widetilde{L}}{\gamma } \lVert u_k\rVert^2 \\
&=& \tfrac{2L\widetilde{L}}{\gamma } \lVert x^{k}-x^{k+1}\rVert^2 \\ 
&\overset{\eqref{quadratic_3}}{=}&\tfrac{2L\widetilde{L}}{\gamma }(1+\tau)^2\lVert x^{k}- \bar x \rVert^2,
\end{eqnarray*}
which proves the quadratic convergence of $\{x_k\}$ to $\bar{x}$.
\end{proof}



\section{Numerical Demonstrations and Execution of Results} \label{section5}
In this section, we implement the proposed Algorithm \ref{algo1} on some numerical experiments. Algorithm \ref{algo1} and its experimentation were performed in MATLAB R2023b software. This MATLAB software is installed in an IOS machine equipped with a 12-core CPU and 8 GB RAM. In the numerical implementation of the algorithm, we have considered the following.  
\begin{itemize}
    \item We take the cone $K$ to be a standard ordering cone, that is, $K=\mathbb{R}^2_+$ for all test instances except for Examples \ref{exam_6} and  \ref{exam_7}. For the scalarizing function $\Psi_e$, we take $e=(1,1,\ldots,1)^\top\in$ int($K$). 
    \item The parameters $\beta$ and $\nu$ in Step \ref{step5} for the line search of the Algorithm \ref{algo1} is chosen as $\beta=0.5$ and $\nu=0.54$.
    \item The employed stopping criterion is $\lVert u_k\rVert <0.001$ or a maximum number of 100 iterations is reached.
    \item To find the set Min($F(x_k),K$) at the $k$-th iteration in Step \ref{step2} of Algorithm \ref{algo1}, \textcolor{black}{we adopted the crude method (pair-wise comparing) of comparing the elements in $F(x_k)$.}
    \item At the $k$-th iteration in Step \ref{step3} of Algorithm \ref{algo1}, we find an $$ (a^k,u_k) \in \underset{(a,u)\in P_k\times\mathbb{R}^n}{\text{argmin}}\xi_{x_k}(a,u)$$ with the help of an inbuilt function \textit{fminsearch} in MATLAB.
    \item We have considered some test problems from the literature subjected to slight modifications, and some are freshly introduced. For each test considered, we generated 100 initial points randomly and ran the algorithm for each of the opted initial points. In the context of each experiment, we have presented a table with four columns: 
    \begin{itemize}
        \item \textbf{Number of initial points:}  This value in the first column is the number of initial points taken to execute the proposed Algorithm \ref{algo1}. 
        \item \textbf{Algorithm:} For the proposed algorithm, we use the abbreviation NM (Newton method), and for the existing steepest descent method \cite{bouza2021steepest}, we use the abbreviation SD. 
        \item \textbf{Iterations:} 
          This value presents the third column with a 6-tuple (Min, Max, Mean, Median, Mode, SD) whose components are the minimum, maximum, mean, median, mode, and standard deviation of the number of iterations until the stopping condition is met. 
          \item \textbf{CPU time:} This value indicates the third column, which is again a 6-tuple (Min, Max, Mean, Median, $\lceil\text{Mode}\rceil$, SD) that shows the minimum, maximum, mean, median, least integer greater or equal to mode, and standard deviation of the CPU time (in seconds) taken by the initial point in reaching the stopping condition.
    \end{itemize}
\end{itemize}

Additionally, the numerical values are presented with precision up to four decimal places to ensure clarity of the difference of function values at different iterations. For every examined problem, the values of the set-valued function $F$ at the initial and final points are marked with black and red colors, respectively. We use shapes $\bullet,~\boldsymbol{\star},$ and $\blacktriangle$ to depict the values $F$ for different initial points. Cyan, magenta, and green colors are used to represent the intermediate iterates for different initial points. Initial points are depicted in black, and the termination point is in red. If the initial point is depicted by black bullet $\bullet$, then the terminating is depicted by the red bullet {\red $\bullet$}, and the intermediate iterates by cyan bullets {\cyan $\bullet$} or magenta bullets {\magenta $\bullet$} or green bullets {\green $\bullet$}. That is, we use the same shape for depicting a complete sequence of iterates generated by Algorithm \ref{algo1}.

Furthermore, we compare the results of the proposed Newton's method (abbreviated as NM) with the existing steepest descent method (abbreviated as SD) for set optimization presented in \cite{bouza2021steepest}. 

Now, we discuss the different test problems on which the proposed algorithm is tested. \textcolor{black}{The first is taken from \cite{kobis2016treatment} with some modifications.}

\begin{example}\label{exam_1}
 Consider the function $F:\mathbb{R}^2\rightrightarrows \mathbb{R}^2$ defined as 
$$F(x)=\{f^1(x),f^2(x)\ldots,f^{20}(x)\},$$
where for each $i\in[20],f^i:\mathbb{R}^2\to\mathbb{R}^2$ is given by 
\[
f^i(x)=\begin{pmatrix}
x_1^2+0.5\sin\left(\tfrac{2\pi(i-1)}{20}\right)+x_2^2\\
2x_1^2+0.5\cos\left(\tfrac{2\pi(i-1)}{20}\right)+2x_2^2
\end{pmatrix}.
\]    

\textcolor{black}{
Figure \ref{figure1} shows the behaviour of Algorithm $\ref{algo1}$ for different initial points within the set $[-4,4]\times[-4,4]$. Firstly, Figure \ref{figure1a} depicts the sequence of iterates $\{F(x_k)\}$ generated by Algorithm \ref{algo1} for the chosen initial point $x_0=(2.5102,0.0000)^\top$. Subsequently, the sequence of iterates $\{x_k\}$ corresponding to $\{F(x_k)\}$ is illustrated in Figure \ref{figure1b}.}

\textcolor{black}{Further, in Figure \ref{figure1c}, the sequence of iterates $\{F(x_k)\}$ generated by Algorithm \ref{algo1} for three randomly selected initial points are depicted with cyan, magenta, and green colors. Additionally, the sequence of iterates $\{x_k\}$ corresponding to $\{F(x_k)\}$ generated by Algorithm \ref{algo1} are shown in Figure \ref{figure1d}}.\\

Next, we have discussed the performance of Algorithm \ref{algo1} for Example \ref{exam_1} in Table \ref{table_1a}. Further, in Table \ref{table_1a}, we have compared the results of NM with SD for set optimization. The values in Table \ref{table_1a} show that the proposed method performs better than the existing SD method.

\begin{figure}[H]
    \centering
    \subfloat[~The value of $F$ at each iteration generated by Algorithm \ref{algo1} for initial point $x_0=(2.5102,0.0000)^\top$ for Example \ref{exam_1}]{\includegraphics[width=0.465\textwidth]{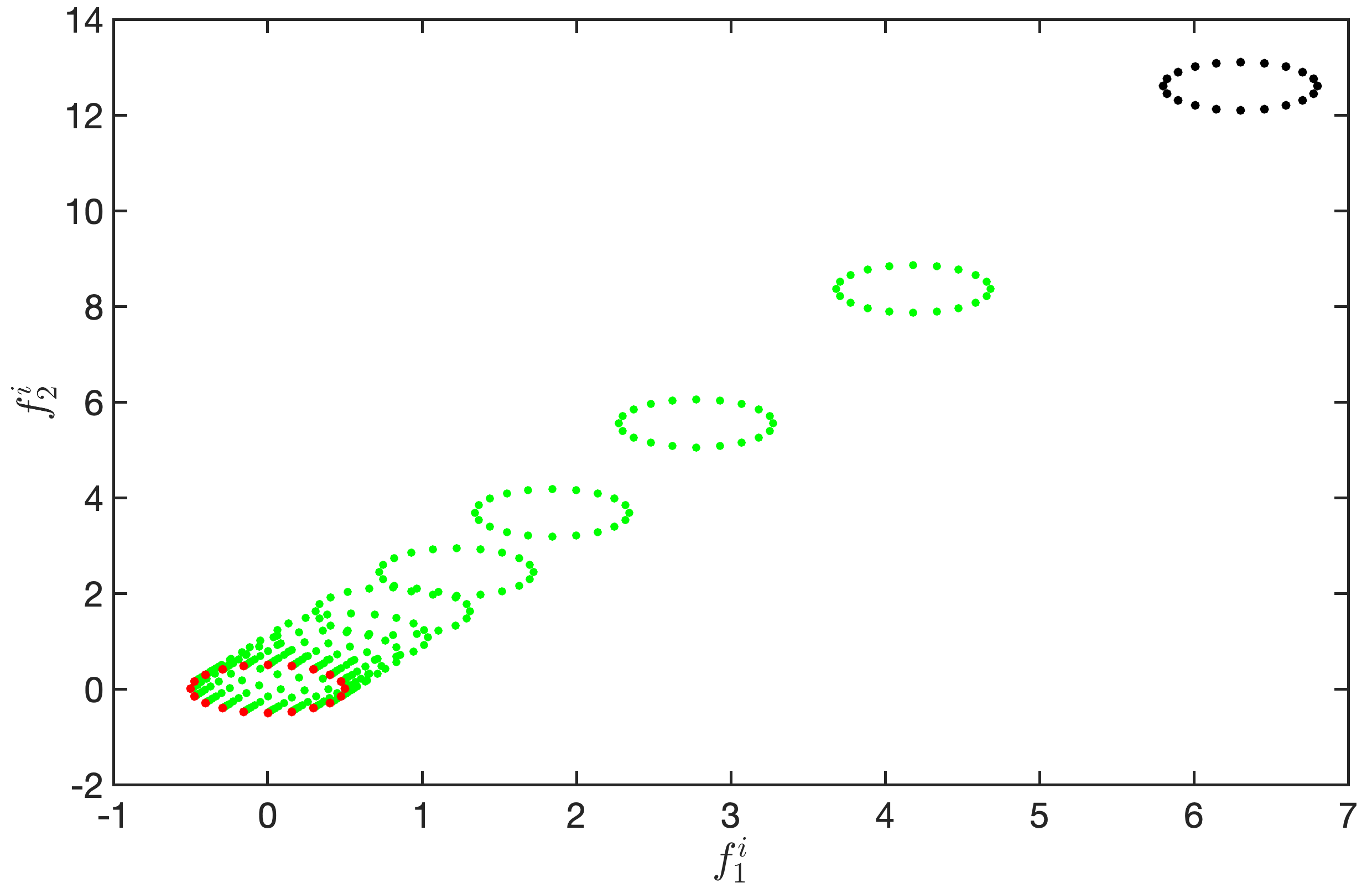}\label{figure1a}} 
    \qquad
    \subfloat[~The value of $x_k$ at each iteration generated by Algorithm \ref{algo1} for initial point $x_0=(2.5102,0.0000)^\top$ for Example \ref{exam_1}]{\includegraphics[width=0.465\textwidth]{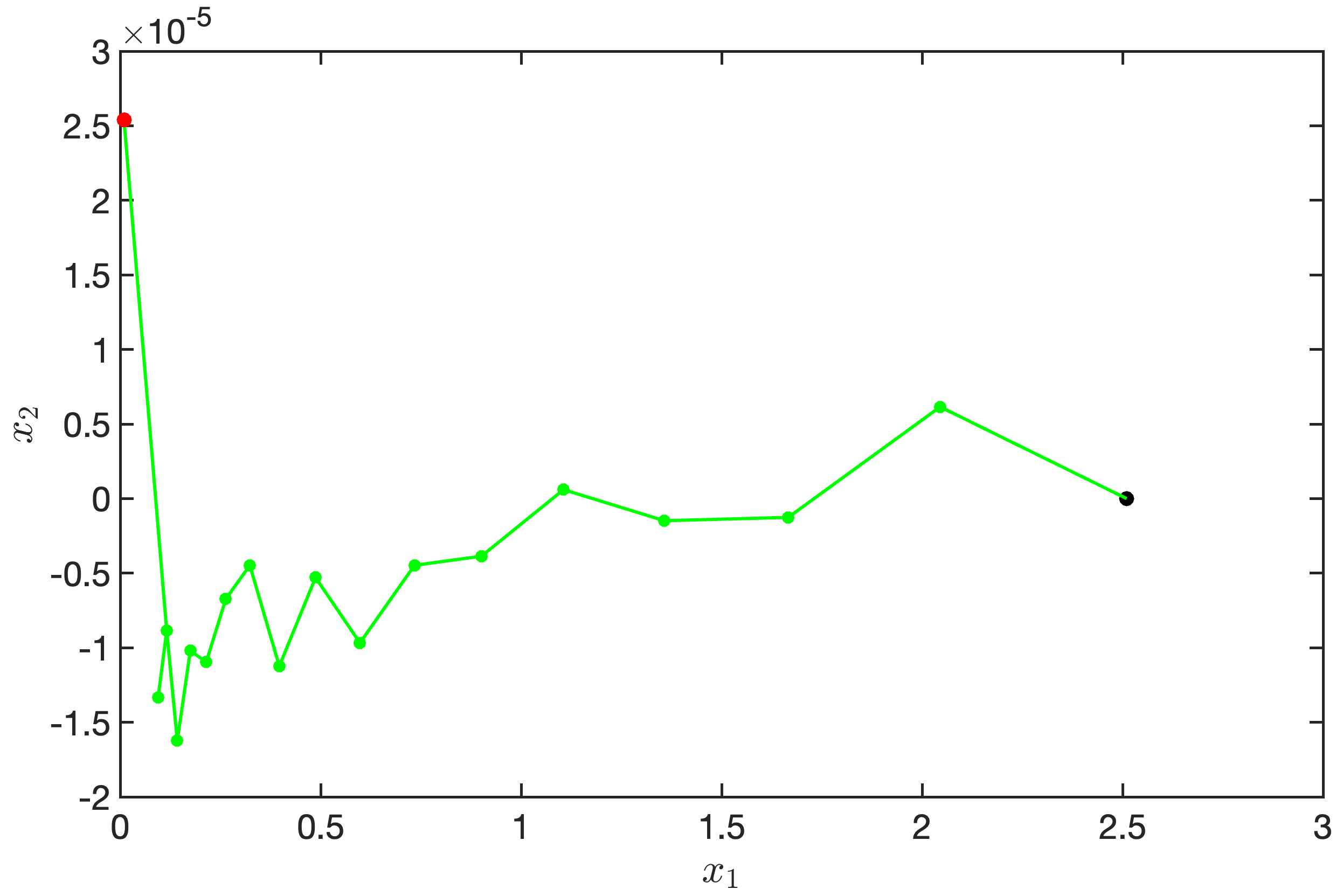}\label{figure1b}} 
    \qquad
    \subfloat[~The value of $F$ at each iteration generated by Algorithm \ref{algo1} for three different randomly chosen initial points for Example \ref{exam_1}]{\includegraphics[width=0.465\textwidth]{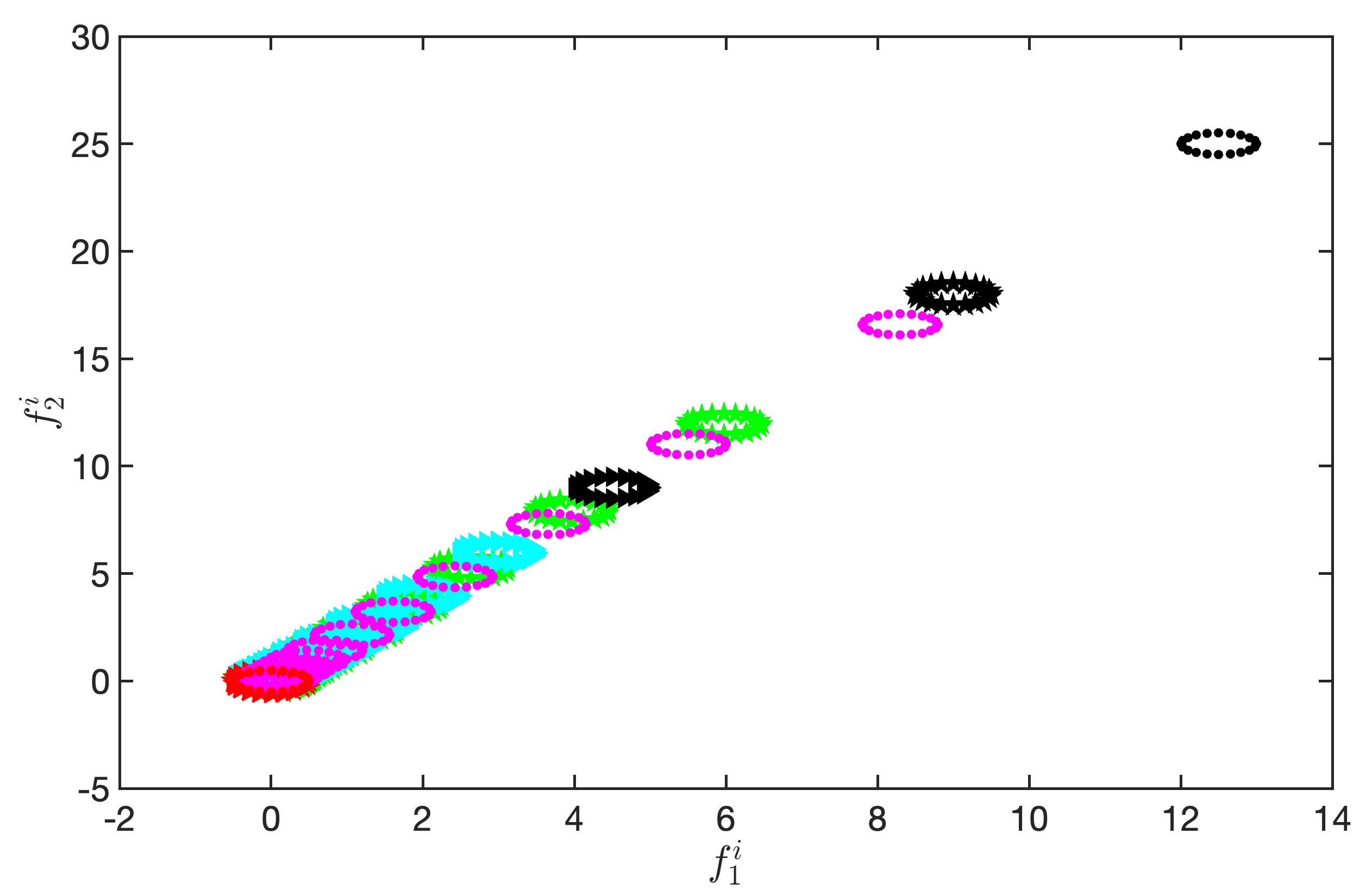}\label{figure1c}}
    \qquad
    \subfloat[~The value of $x_k$ at each iteration generated by Algorithm \ref{algo1} for three different randomly chosen initial points for Example \ref{exam_1}]{\includegraphics[width=0.465\textwidth]{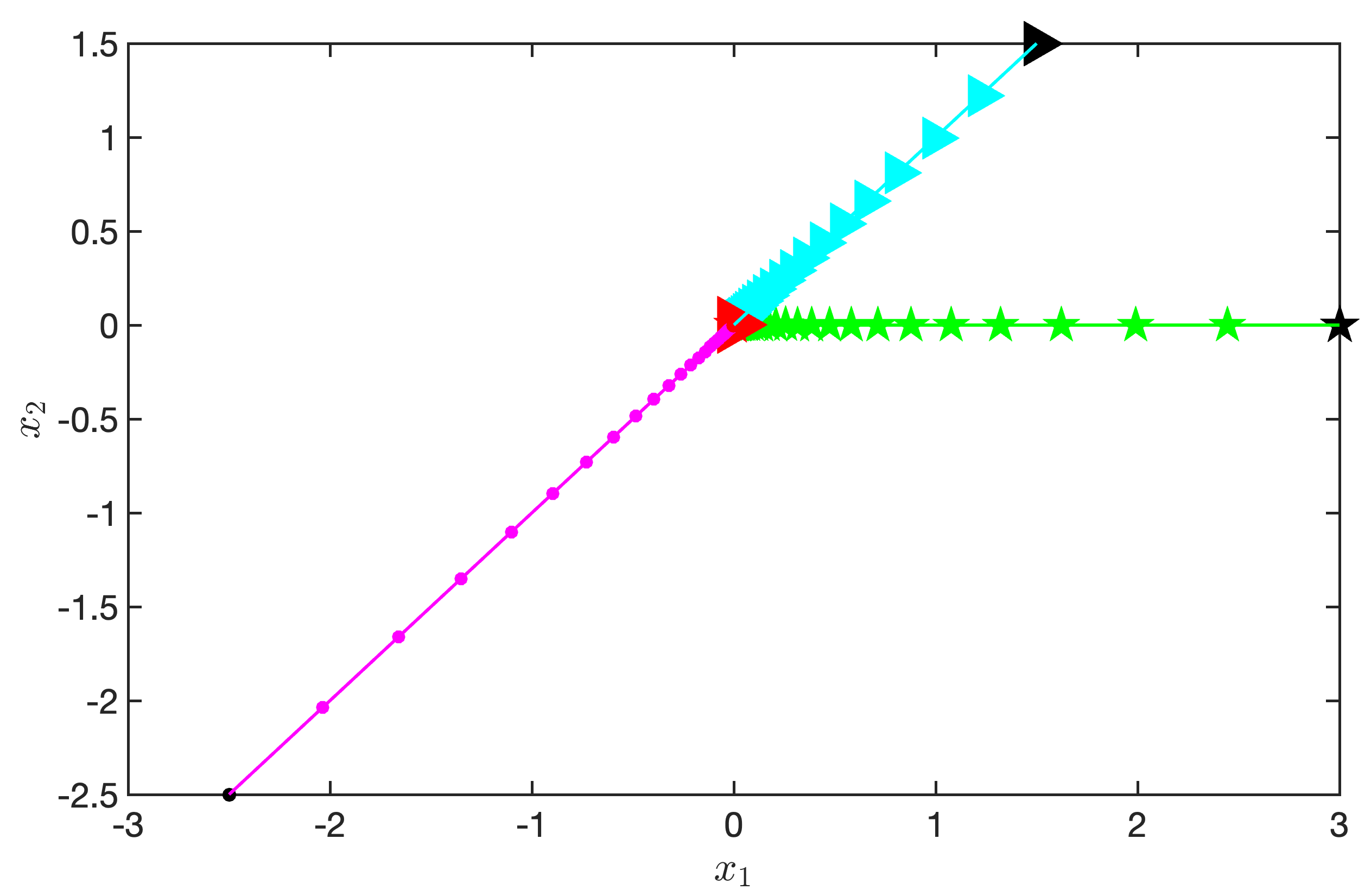}\label{figure1d}}
    \caption{Output of Algorithm \ref{algo1} for Example \ref{exam_1}}
    \label{figure1}
\end{figure}

\begin{table}[H]
\centering
\caption{Performance of Algorithm \ref{algo1} on Example \ref{exam_1}}
\centering 
\scalebox{0.71}{
\begin{tabular}{c c c c} 
\hline  
Number of &Algorithm  &Iterations &CPU time\\
initial points& & (Min, Max, Mean, Median, Mode, SD)  &(Min, Max, Mean, Median, $\lceil\text{ Mode }\rceil$, SD)   \\ 
\hline 
100 & NM ($t_k=1$) &($2,~2,~2,~2,~2,~0$) & ($2.0521,~3.1308,~2.1637,~2.1031,~2,~0.1601$) \\ 
& NM &($1,~14,~10.6800,~11,~12,~2.6963$) & ($1.0345,~25.6842,~18.9496,~20.4787,~1,~4.9776$) \\ 
& SD &($~1,~23,~11.4800,~11,~10,~5.7375$)&($1.0182,~41.4550,~19.9548,~,19.5869,~1,~10.6772$)\\
\hline 
\end{tabular}}
\label{table_1a}
\end{table}

\end{example}


\textcolor{black}{The next example, discussed below, is freshly introduced.}

\begin{example}\label{exam_2}
Consider the set-valued function $F:\mathbb{R}\rightrightarrows \mathbb{R}^2$ defined by
$$F(x)=\{f^1(x),f^2(x),\ldots,f^{50}(x)\},$$
where for each  $i\in[50],$ the function $f^i:\mathbb{R}\to\mathbb{R}^2$ is given by 
\[
f^i(x)=\begin{pmatrix}
0.35\sin\left(\tfrac{2\pi(i-1))}{50}\right)\cos\left(\tfrac{2\pi(i-1))}{50}\right)+x^2\\
0.35\cos\left(\tfrac{2\pi(i-1))}{50}\right)+\tfrac{1}{(1+e^{2x})}+\cos(2x)\\
\end{pmatrix}.
\]
\begin{figure}[ht]
\centering
\mbox{\subfloat[The value of $F$ at each iteration generated by Algorithm \ref{algo1} for Example \ref{exam_2} for the initial point $x_0=2.0000$]{ \includegraphics[width=0.48\textwidth]{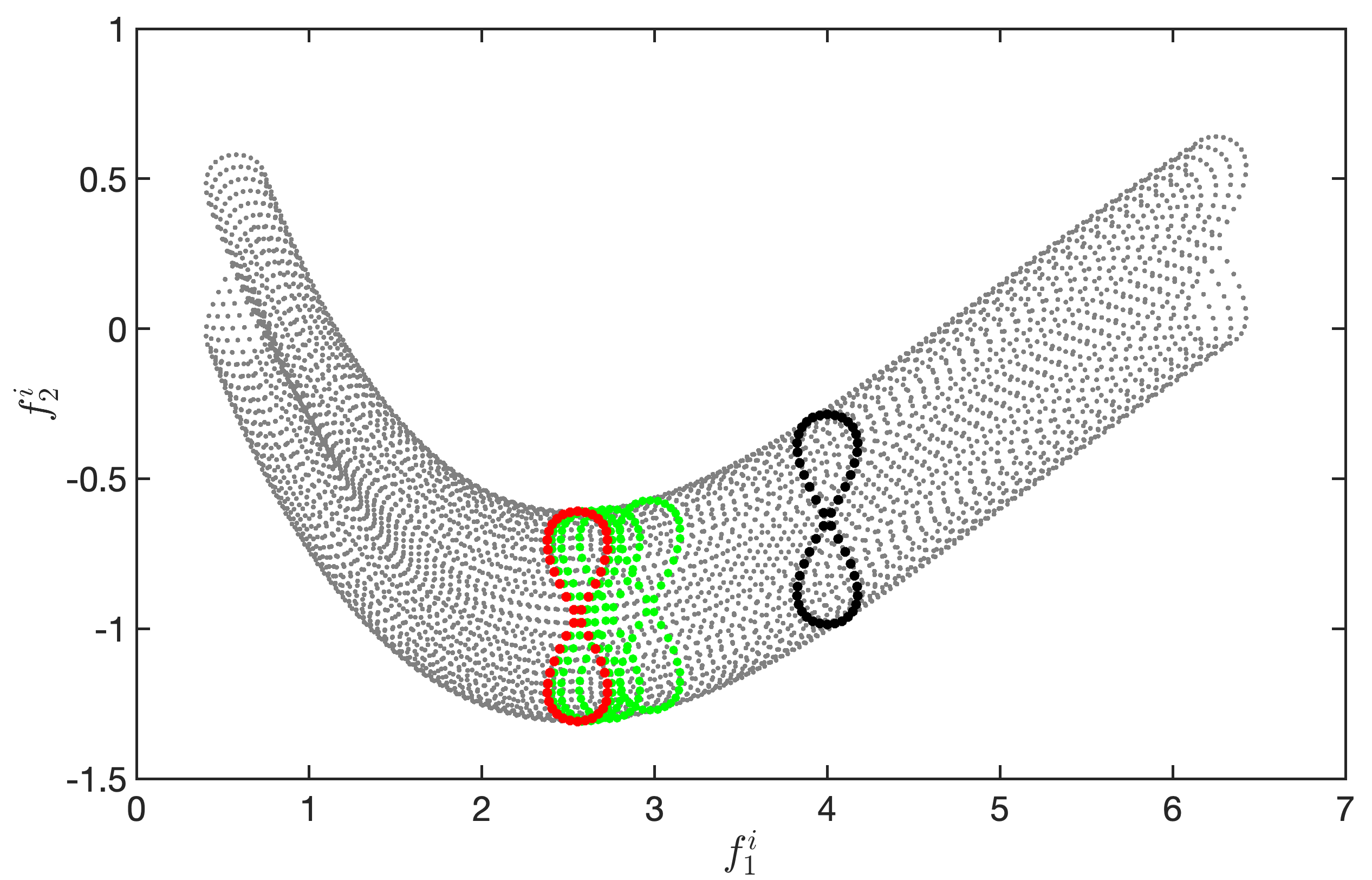}\label{figure2a}}\quad
\subfloat[The value of $F$ for three randomly chosen different initial points at each iteration generated by Algorithm \ref{algo1} of Example \ref{exam_2}]{\includegraphics[width=0.48\textwidth]{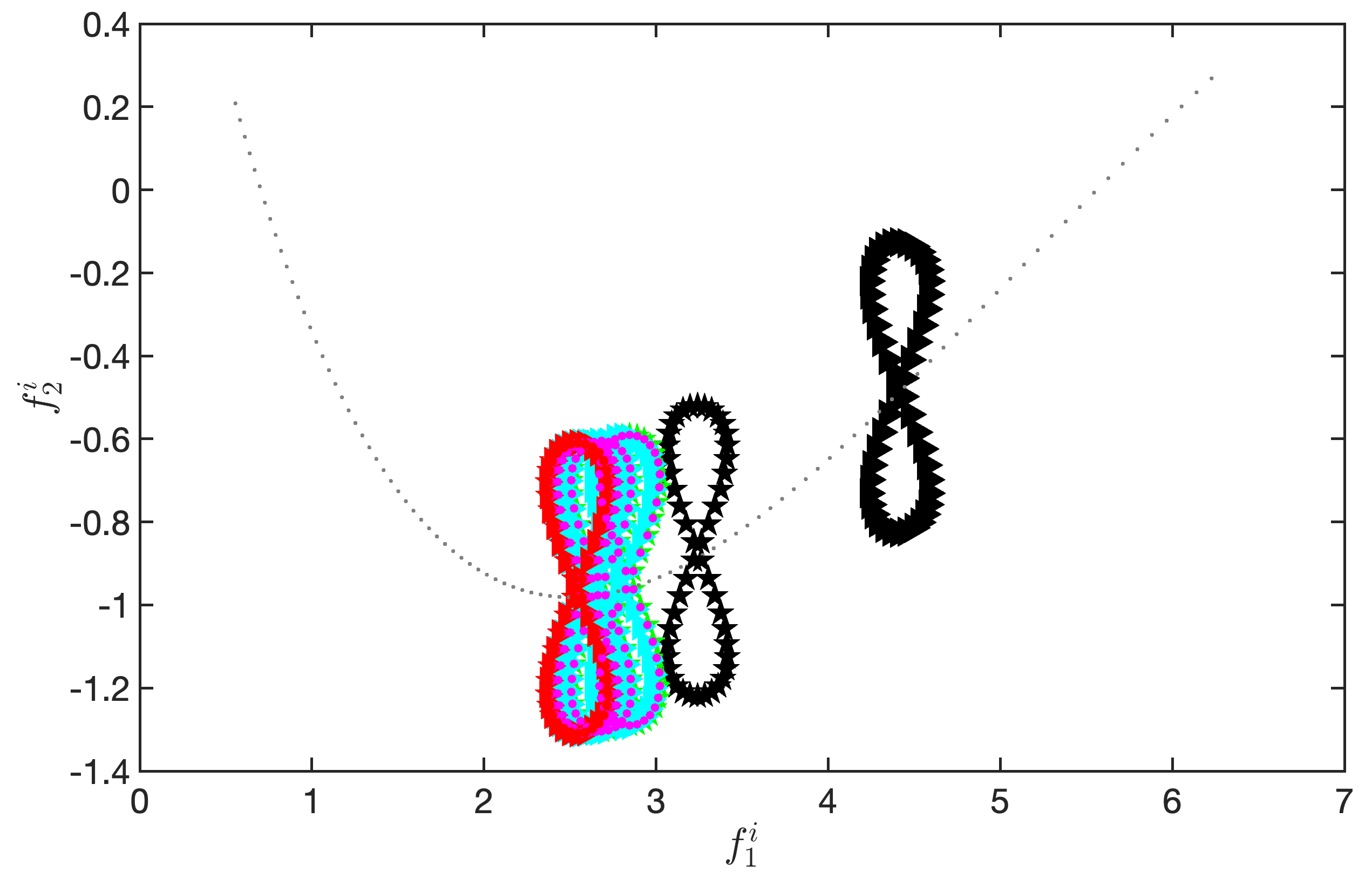}\label{figure2b}}\quad}
\caption{Obtained output of Algorithm \ref{algo1} for Example \ref{exam_2}}
\label{figure2} 
\end{figure}

The output of Algorithm \ref{algo1} for different initial points of Example \ref{exam_2} are depicted in Figure \ref{figure2}. The discretized $\infty$-shaped segments represent the objective values that transverse a curve within the interval $[0.7700,6.3000]$. Figure \ref{figure2a} depicts the sequence $\{F(x_k)\}$ generated by Algorithm \ref{algo1} for the starting point $x_0=2.0000$. \textcolor{black}{In Figure \ref{figure2b}, we exhibit the output sequence $\{F(x_k)\}$ generated by Algorithm \ref{algo1} for three randomly chosen initial points.}\\ 

The performance of Algorithm \ref{algo1} for Example \ref{exam_2} is shown in Table \ref{table_ex2}. Moreover, we have compared the results of the NM with the SD method for set optimization as presented in Table \ref{table_ex2}. The values in Table \ref{table_ex2} show that the proposed method performs better than the existing SD method.  

\begin{table}[!htb]
\centering
\caption{Performance of Algorithm \ref{algo1} on Example \ref{exam_2}}
\centering 
\scalebox{0.73}{
\begin{tabular}{c c c c} 
\hline 
Number of &Algorithm  &Iterations &CPU time\\
initial points& &(Min, Max, Mean, Median, Mode, SD)  &(Min, Max, Mean, Median, $\lceil\text{Mode}\rceil$, SD)   \\ 
\hline 
$100$ & NM ($t_k=1$) &($1,~3,~1.9900,~2,~2,~0.3013$) & ($0.2215,~1.6828,~0.8634,~0.8416,~0,~0.1625$) \\ 
& NM &($1,~5,~2.3600,~1,~1,~1.8175$) & ($0.3823,~3.9151,~1.6241,~0.9135,~0,~0.9606$) \\
 & SD &($1,~8,~6.0900,~7,~8,~2.4622$) & ($0.4170,~5.4045,~3.7906,~4.2701,~0,~1.3692$) \\
\hline 
\end{tabular}}
\label{table_ex2}
\end{table}

\end{example}

\textcolor{black}{Below, we consider Examples \ref{exam_3} and \ref{exam_4}, again motivated from \cite{kobis2016treatment}.}

\begin{example}\label{exam_3}
 Consider the function $F:\mathbb{R}^2\rightrightarrows \mathbb{R}^3$ defined as 
$$F(x)=\{f^1(x),f^2(x),\ldots,f^{14}(x)\},$$
where for each $i\in[14]$, $f^i:\mathbb{R}^2\to\mathbb{R}^3$ is given by
\[
f^i(x)=\begin{pmatrix}
x_1^2+x_2^2+0.25\sin\left(\tfrac{2\pi(i-1)}{14}\right)\\
4x_1^2+4x_2^2+0.25\cos\left(\tfrac{2\pi(i-1)}{14}\right)\\
x_1^2+x_2^2+i
\end{pmatrix}.
\]   

\textcolor{black}{
Figure \ref{figure3} shows the behaviour of Algorithm $\ref{algo1}$ for different initial points within the set $[-3,4]\times[-3,4]$. Firstly, Figure \ref{figure3a} depicts the sequence of iterates $\{F(x_k)\}$ generated by Algorithm \ref{algo1} for the chosen initial point $x_0=(3.2302,-0.5102)^\top$. Subsequently, the sequence of iterates $\{x_k\}$ corresponding to $\{f(x_k)\}$ is illustrated in Figure \ref{figure3b}}. 

\textcolor{black}{Further, in Figure \ref{figure3c}, the sequence of iterates $\{F(x_k)\}$ generated by Algorithm \ref{algo1} for three randomly selected initial points are depicted with cyan, magenta, and green colors. Additionally, the sequence of iterates $\{x_k\}$ corresponding to $\{F(x_k)\}$ generated by Algorithm \ref{algo1} are shown in Figure \ref{figure3d}}. 

The performance of Algorithm \ref{algo1} on Example \ref{exam_3} is shown in Table \ref{table_3a}. The comparison of the results of NM with the SD method for set optimization is given in Table \ref{table_3a}. The values in Table \ref{table_3a} show that the proposed method performs better than the existing SD method.  

\begin{figure}
    \centering
    \subfloat[~The value of $F$ at each iteration generated by Algorithm \ref{algo1} for initial point $x_0=(3.2302,-0.5102)^\top$ for Example \ref{exam_3}]{\includegraphics[width=0.47\textwidth]{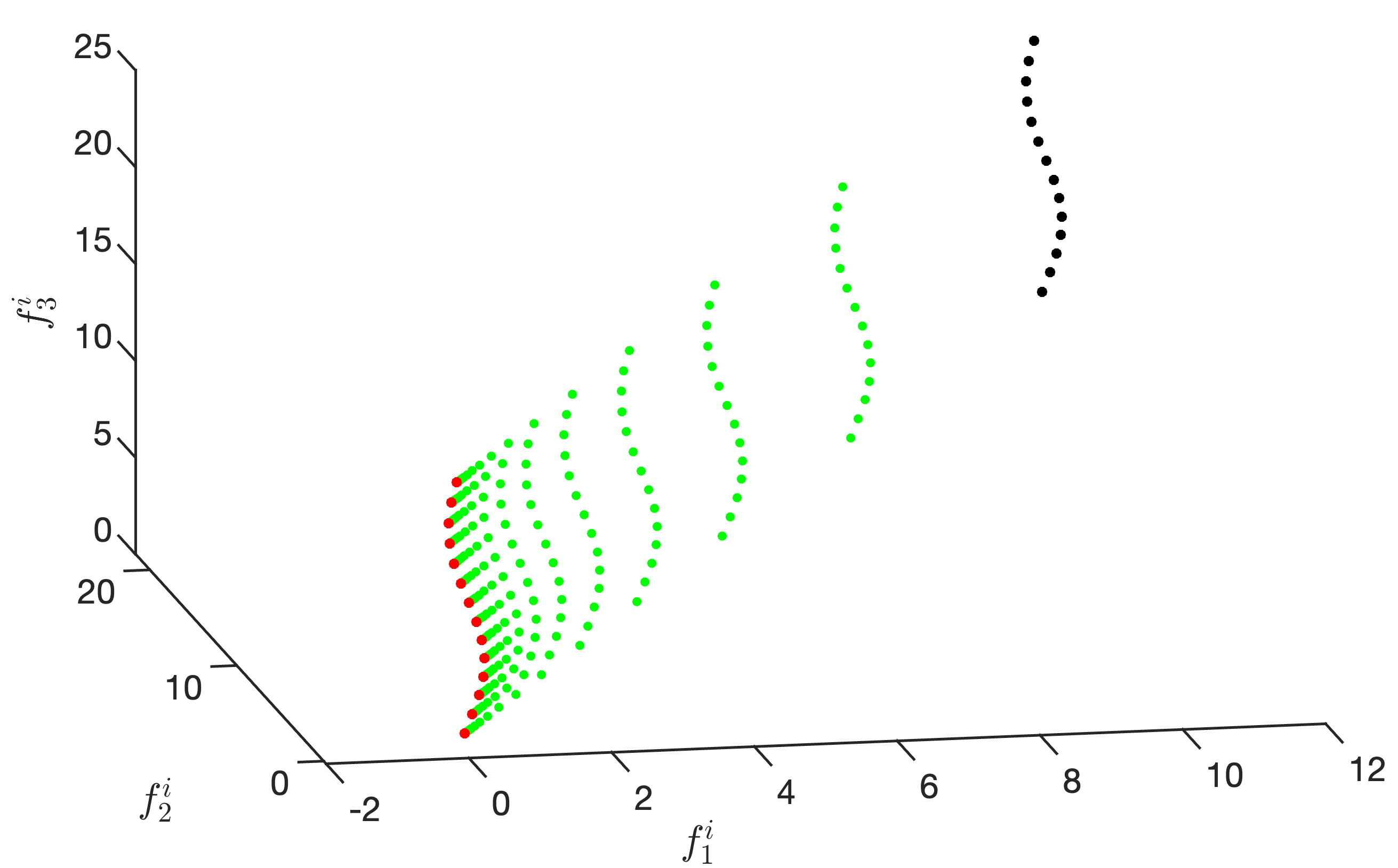}\label{figure3a}}
    \qquad
    \subfloat[~The value of $x_k$ at each iteration generated by Algorithm \ref{algo1} for initial point $x_0=(3.2302,-0.5102)^\top$ for Example \ref{exam_3}]{\includegraphics[width=0.47\textwidth]{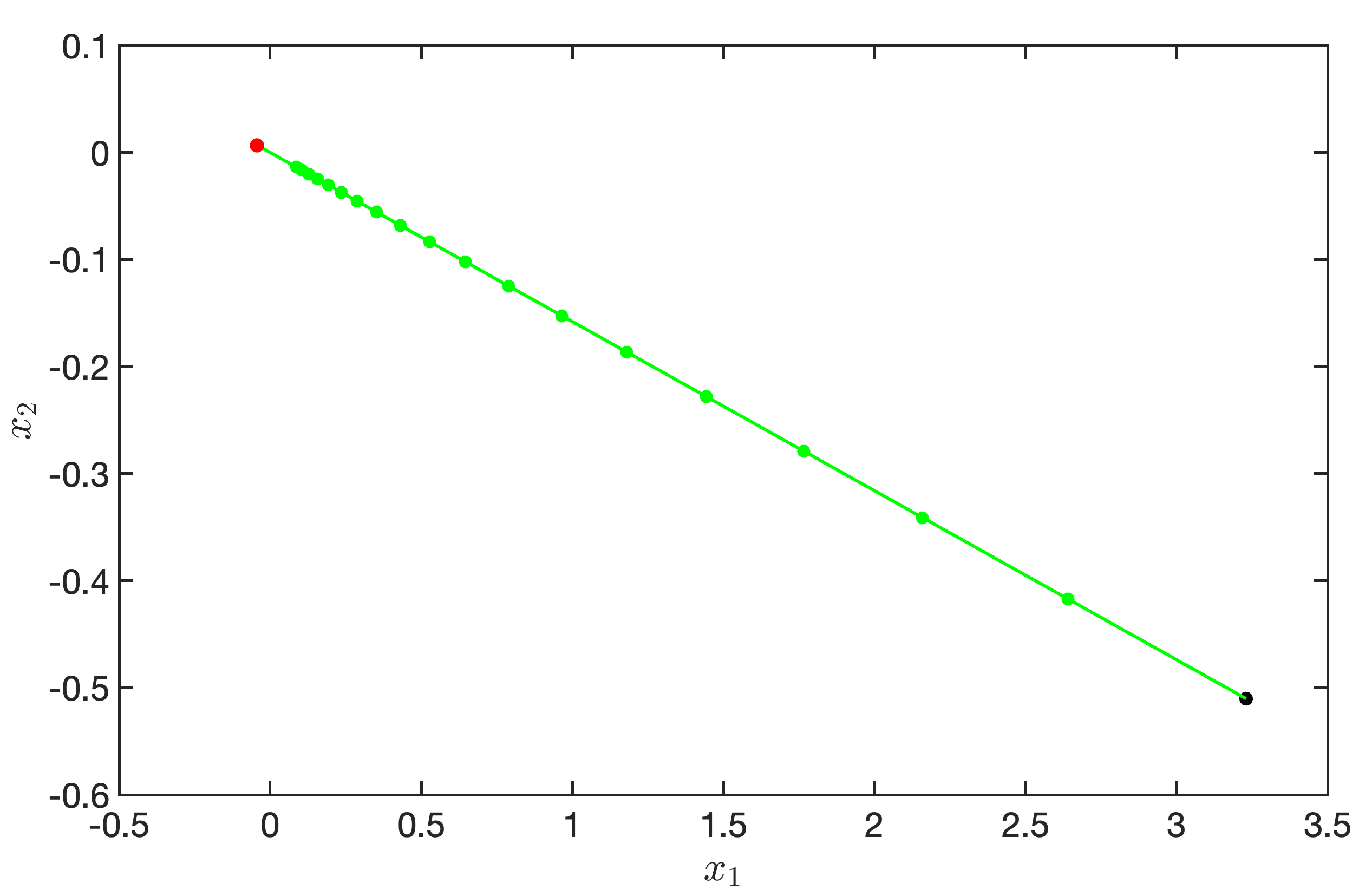}\label{figure3b}} 
    \qquad
    \subfloat[~The value of $F$ at each iteration generated by Algorithm \ref{algo1} for three different randomly chosen initial points for Example \ref{exam_3}]{\includegraphics[width=0.47\textwidth]{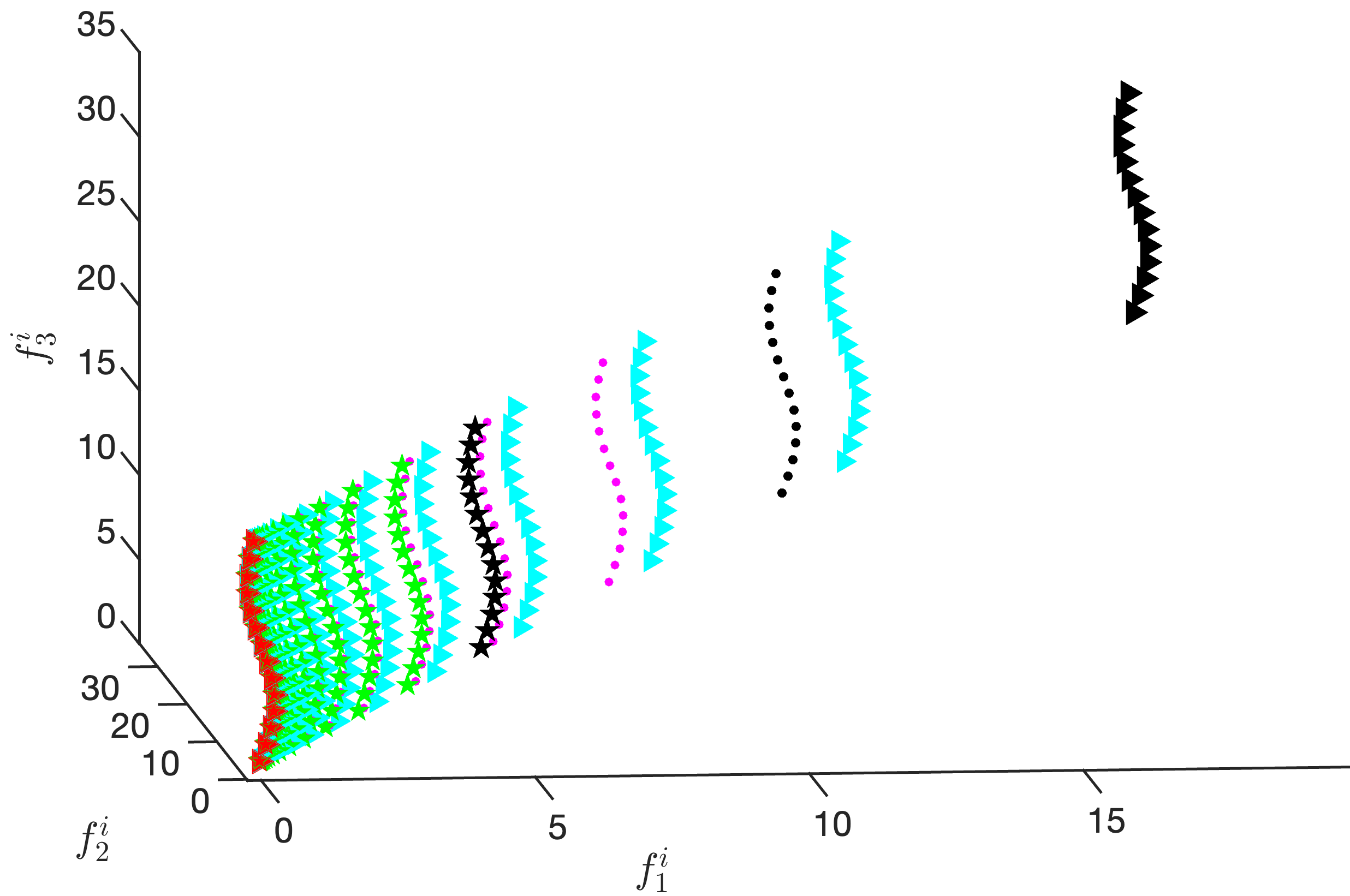}\label{figure3c}}
    \qquad
    \subfloat[~The value of $x_k$ at each iteration generated by Algorithm \ref{algo1} for three different randomly chosen initial points for Example \ref{exam_3}]{\includegraphics[width=0.47\textwidth]{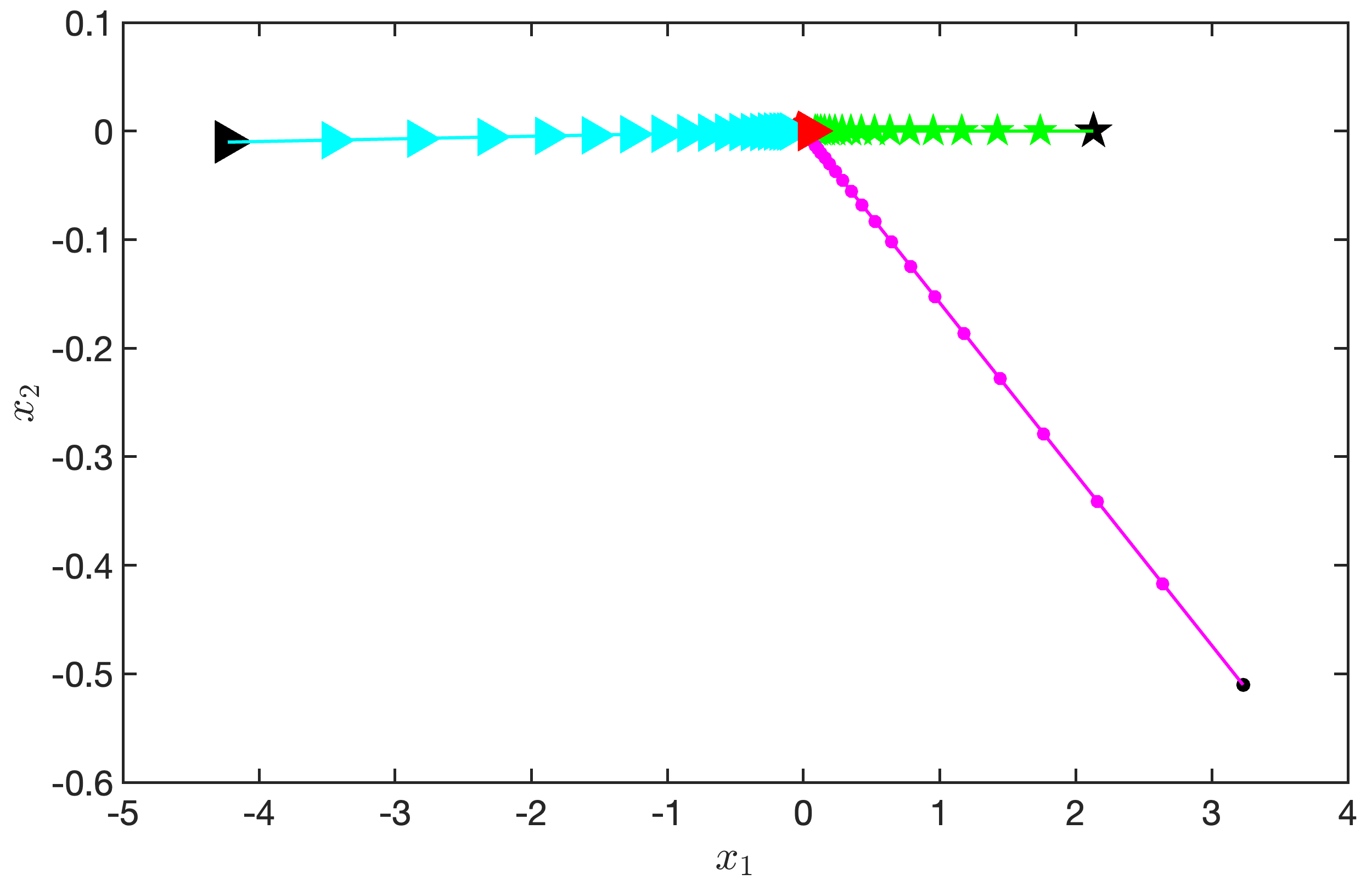}\label{figure3d}}
    \caption{Output of Algorithm \ref{algo1} for Example \ref{exam_3}}
    \label{figure3}
\end{figure}

\begin{table}[ht]
\centering
\caption{Performance of Algorithm \ref{algo1} on Example \ref{exam_3}}
\centering 
\scalebox{0.69}{
\begin{tabular}{c c c c} 
\hline 
Number of &Algorithm &Iterations &CPU time\\
initial points& & (Min, Max, Mean, Median, Mode, SD)  &(Min, Max, Mean, Median, $\lceil\text{Mode}\rceil$, SD)   \\ 
\hline 
$100$  & NM ($t_k=1$) &($2,~2,~2,~2,~2,~0$) & ($5.9813,~7.5494,~6.3048,~6.3082,~5,~0.1751$) \\ 
&NM &($1,~14,~10.8800,~12,~13,~2.8508$) & ($0.1190,~69.3281,~51.3830,~56.2429,~13,~13.9629$) \\
& SD &($1,~14,~11.0500,~12,~13,~2.2128$) & ($0.1202,~65.2457,~50.2756,~54.3065,~21,~10.3982$) \\
\hline 
\end{tabular}}
\label{table_3a}
\end{table}
\end{example}


\begin{example}\label{exam_4}
Consider the function $F:\mathbb{R}\rightrightarrows \mathbb{R}^3$ defined by 
$$F(x)=\{f^1(x),f^2(x),\ldots,f^{30}(x)\},$$
where for each $i\in[30],f^i:\mathbb{R}\to\mathbb{R}^3$ is given by 
\[
f^i(x)=\begin{pmatrix}
x^2+\tfrac{(i-1)}{30}\\
(x^2-4)(\sin(x^2-4))+\tfrac{(i-1)}{30}\\
\tfrac{(i-1)}{30} ~ x^2
\end{pmatrix}.
\]
  
\begin{figure}[ht]
\centering
\mbox{\subfloat[~The value of $F$ at each iteration generated by Algorithm \ref{algo1} for initial point $x_0=2.1300$ for Example \ref{exam_4}]{ \includegraphics[width=0.48\textwidth]{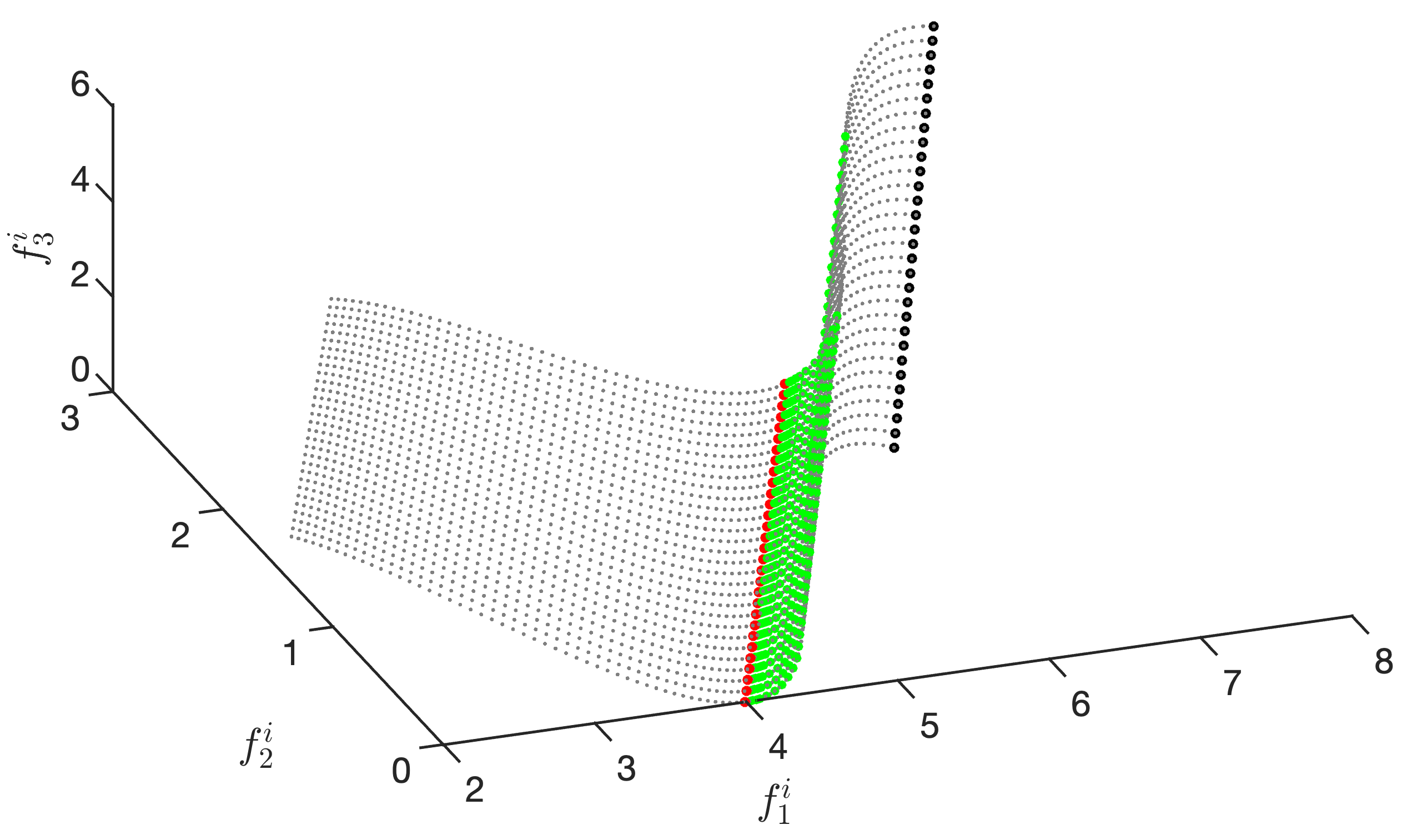}\label{figure4a}}\quad
\subfloat[~The value of $F$ at each iteration generated by Algorithm \ref{algo1} for three different randomly chosen initial points for Example \ref{exam_4}]{\includegraphics[width=0.48\textwidth]{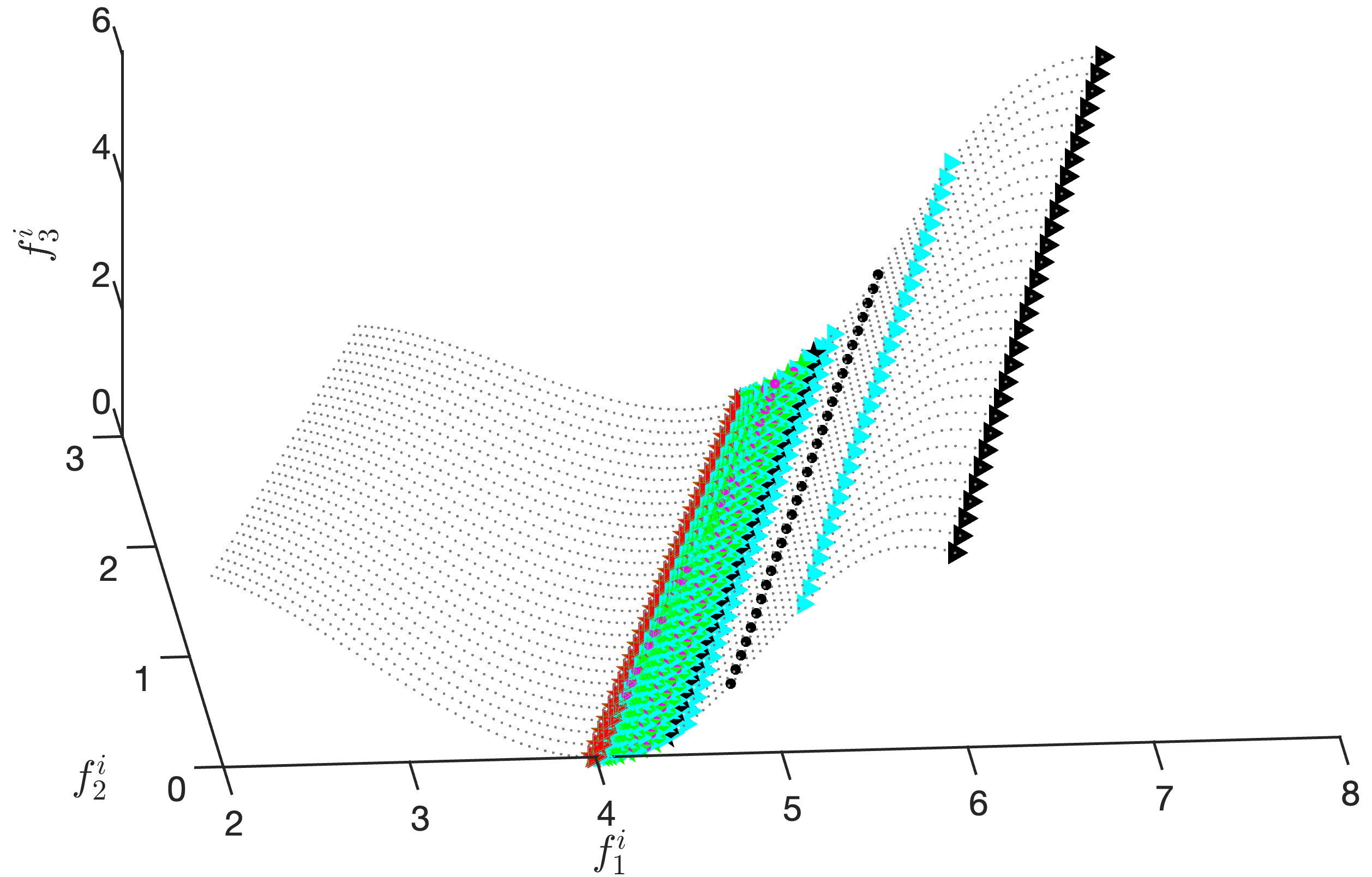} \label{figure4b}}\quad}
\caption{Obtained output of Algorithm \ref{algo1} for Example \ref{exam_4}}
\label{figure4} 
\end{figure}
The output of Algorithm \ref{algo1} for Example \ref{exam_4} is depicted in Figure \ref{figure4}. The collection of the objective values for all $x \in [1.5400, 2.1600]$ transverses a surface depicted by the dot-shaded region. Figure \ref{figure4a} depicts the sequence $\{F(x_k)\}$ generated by Algorithm \ref{algo1} for the starting point $x_0 = 2.1300$. We test Algorithm \ref{algo1} for three different randomly chosen initial points and observe their corresponding output $\{F(x_k)\}$ (marked with cyan, magenta, and green colors) as shown in Figure \ref{figure4b}. The performance of Algorithm \ref{algo1} for Example \ref{exam_4} is shown in Table \ref{table_4}. 

\begin{table}[ht]
\centering
\caption{Performance of Algorithm \ref{algo1} on Example \ref{exam_4}}
\centering 
\scalebox{0.64}{
\begin{tabular}{c c c c} 
\hline 
Number of& Algorithm  &Iterations &CPU time\\
initial points& &(Min, Max, Mean, Median, Mode, SD)  &(Min, Max, Mean, Median, $\lceil\text{Mode}\rceil$, SD)   \\ 
\hline 
$100$  & NM ($t_k=1$) &($2,~2,~2,~2,~2,~0$) & ($274.9991,~301.5326,~288.4099,~288.4619,~274,~7.6985$) \\ 
& NM &($0,~19,~7.4200,~9,~1,~5.1231$) & ($0.0018,~267.5772,~125.0167,~129.7198,~0,~ 81.4232$) \\  
\hline 
\end{tabular}}
\label{table_4}
\end{table}

For the initial point $x_0=2.1300$, the values of $F$ across the generated iterates for NM and SD methods are exhibited in \textcolor{black}{Table \ref{table_4b} and \ref{table_4c}, respectively}. We have observed that the SD method presented in \cite{bouza2021steepest} is not properly solving this problem due to the objective function's highly nonlinear nature and converging to the given initial point $x_0$ in zero number of iterations. On the other hand, the proposed NM method gives the minimal solution $\bar{x}=1.9957$ for the chosen initial point $x_0=2.1300$. Therefore, we can conclude that the proposed NM Algorithm \ref{algo1} is more efficient than the SD method for set optimization given in \cite{bouza2021steepest}.

\begin{table}[ht]
\centering
\caption{Performance of proposed NM (Algorithm \ref{algo1}) on Example \ref{exam_4} with $t_k\in(0,1)$}
\centering 
\scalebox{0.83}{
\begin{tabular}{c c c c c } 
\hline  
Iteration  ~&  ~& ~& ~&  \\
number $(k)$& $x_k$ &$f^{10}(x_k)$ & $f^{20}(x_k)$ & $f^{30}(x_k)$   \\ 
\hline 
 0&  $2.1300$& $(4.8369,~0.5746,~1.3611)$      &$(5.1702,~0.9079,~2.8734)$  
& $(5.5036,~1.2413,~4.3857)$\\  
 1&  $2.1115$& $( 4.7584,~0.5029,~1.3375)$       & $(5.0918,~0.8362,~2.8237)$
& $(5.4251,~1.1695,~4.3098)$\\
2&  $2.0929$    &$(4.6802,~0.4411,~1.3141)$      & $(5.0136,~0.7744,~2.7741)$ &$(5.3469,~1.1078,~4.2342)$\\
3&  $2.0744$& $(4.6031,~0.3905,~1.2909)$ & $(4.9365,~0.7238,~2.7253)$    
& $(5.2698,~1.0572,~4.1597)$ \\
4&  $2.0559$& $( 4.5267,~0.3510,~1.2680)$       & $( 4.8601,~0.6843,~2.6769)$   
& $(5.1934,~1.0176,~4.0858)$ \\
5& $2.0373$& $(4.4506,~0.3226,~1.2452)$ & $(4.7839,~0.6559,~2.6287)$ & $(5.1173,~0.9893,~4.0122)$ \\
6&  $2.0244$& $(4.3982,0.3096,1.2295)$ &$(4.7315,0.6430,2.5955)$    
& $(5.0649,0.9763,3.9616)$\\
7& $2.0179$ & $(4.3719,~0.3052,~1.2216)$       & $(4.7053,~0.6385,~2.5789)$    
& $(5.0386,~0.9718,~3.9362)$ \\
8&  $2.0114$& $(4.3457,~0.3021,~1.2137)$  &$(4.6791,~0.6354,~2.5623)$  
& $(5.0124,~0.9688,~3.9109)$\\
9&  $1.9957$& $(4.2828,~0.3003,~1.1948)$  
& $(4.6162,~0.6336,~2.5225)$
&$(4.9495,~0.9670,~3.8501)$ \\
\hline 
\end{tabular}}
\label{table_4b}
\end{table}

\begin{table}[ht]
\centering
\caption{Performance of SD method (see \cite{bouza2021steepest}) on Example \ref{exam_4} with $t_k\in(0,1)$}
\centering 
\scalebox{0.83}{
\begin{tabular}{c c c c c } 
\hline  
Iteration  ~&  ~& ~& ~&  \\
number $(k)$& $x_k$ &$f^{10}(x_k)$ & $f^{20}(x_k)$ & $f^{30}(x_k)$   \\ 
\hline 
 0&  $2.1300$& $(4.8369,~0.5746,~1.3611)$      &$(5.1702,~0.9079,~2.8734)$  
& $(5.5036,~1.2413,~4.3857)$\\  
 1&  $2.1300$& $(4.8369,~0.5746,~1.3611)$      &$(5.1702,~0.9079,~2.8734)$  
& $(5.5036,~1.2413,~4.3857)$\\  
\hline 
\end{tabular}}
\label{table_4c}
\end{table}

\end{example}


\textcolor{black}{In the next example, we consider the robust counterpart of a vector-valued facility location problem under uncertainty \cite{ide2014relationship}. A detailed discussion on this problem is given in \cite{bouza2021steepest}.}
\begin{example}\label{exam_5}
Consider the function $F:\mathbb{R}^2\rightrightarrows \mathbb{R}^3$ defined as 
$$F(x)=\{f^1(x),f^2(x),\ldots,f^{100}(x)\},$$
where for each $i\in[100],f^i:\mathbb{R}^2\to\mathbb{R}^3$ is given as
\[
f^i(x)=\tfrac{1}{2}
\begin{pmatrix}
\lVert x-l_1-u_i\rVert^2\\
\lVert x-l_2-u_i\rVert^2\\
\lVert x-l_3-u_i\rVert^2
\end{pmatrix},
\]  
where $l_1=\begin{pmatrix}
0\\
8
\end{pmatrix},~l_2=\begin{pmatrix}
0\\
0
\end{pmatrix} ~\text{ and }~l_3=\begin{pmatrix}
8\\
0
\end{pmatrix}$.  We consider a uniform partition set of $10$ points of the interval $[-1,1]$ given by 
$$
\mathcal{U}=\left\{-1,-1+\tfrac{1}{s},-1+\tfrac{2}{s},\ldots,-1+\tfrac{2(s-1)}{s},1 \right\}\text{ with }s=4.5.
$$
The set $\{u_i=(u_{1i},u_{2i})^\top:i\in[100]\}$ is an enumeration of the set $\mathcal{U}\times \mathcal{U}$.\\

In Figure \ref{figure_5}, the total of $100$ initial points were generated in the square $[-50,50]\times[-50,50]$. The grey points represent the set $(l_1+u_i)\cup(l_2+u_i)\cup(l_3+u_i)$ and the locations of $l_1,l_2,l_3$ are depicted in blue color. The values of 
$F(x_k)$ generated by Algorithm \ref{algo1} for three different randomly chosen initial points are given with cyan, magenta, and green colors as shown in Figure \ref{figure_5}. 

The performance of Algorithm \ref{algo1} on Example \ref{exam_5} is shown in Table \ref{table_5a}. A comparison of the results of NM with the existing SD method is presented in Table \ref{table_5a}.  The values in Table \ref{table_5a} show that the proposed method performs better than the existing SD method.   

\textcolor{black}{Next, the performance of the proposed Algorithm \ref{algo1} for the initial point $x_0 = (-2.0000, 10.0000)^\top$ with $\nu\in(0,1)$ and $t_k=1$ is shown in Table \ref{table_5b} and \ref{table_b_quadratic}. The decreasing behavior in the values of vector-valued functions at each iteration has been exhibited in Table \ref{table_5b}.}

\begin{figure}[ht]
\centering
\mbox{\subfloat[The value of $F$ in argument space at each iteration generated by Algorithm \ref{algo1} for three randomly chosen initial points for Example \ref{exam_5}]{ \includegraphics[width=0.6\textwidth]{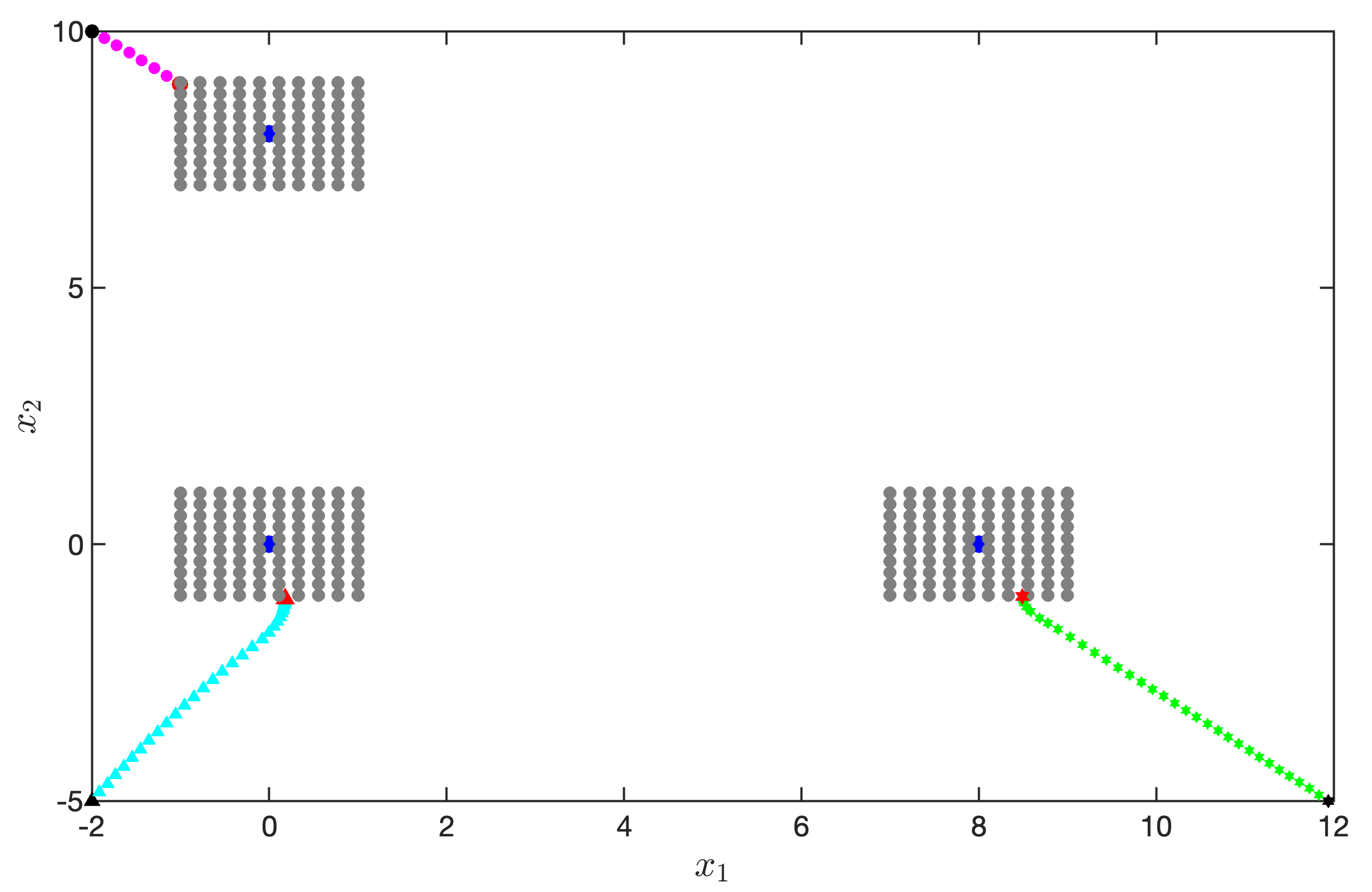}}\quad}
\caption{Obtained output of Algorithm \ref{algo1} for Example \ref{exam_5}}
\label{figure_5} 
\end{figure}

\begin{table}[ht]
\centering
\caption{Performance of Algorithm \ref{algo1} on Example \ref{exam_5}}
\centering 
\scalebox{0.73}{
\begin{tabular}{c c c c} 
\hline  
Number of &Algorithm &Iterations &CPU time\\
initial points& &(Min, Max, Mean, Median, Mode, SD)  & (Min, Max, Mean, Median, $\lceil\text{Mode}\rceil$, SD) \\ 
\hline
$100$  & NM ($t_k=1$) &($3,~2.1000,~2~,2,~0.3892$) & ($0.3677,~3.8689,~0.6442,~0.3840,~0,~0.8040$) \\ 
&NM &($1,~14,~6.4700,~7,~1,~3.9503$) & ($0.2631,~17.114,~5.0776,~4.2074,~0,~3.9056$) \\
 &SD &($1,~65,~14.4200,~5.5000,1,~16.6243$) & ($0.3871,~23.1876,~6.8118,~4.5792,~0,~5.8667$) \\
\hline 
\end{tabular}}
\label{table_5a}
\end{table}

\begin{table}[ht]
\centering
\caption{Performance of Algorithm \ref{algo1} on Example \ref{exam_5} with $\nu\in(0,1)$}
\centering 
\scalebox{0.63}{
\begin{tabular}{c c c c c c} 
\hline  
Iteration  ~&  ~& ~& ~& ~& \\
number $(k)$& $x_k^\top$ &$f^{25}(x_k)$ & $f^{50}(x_k)$ & $f^{75}(x_k)$ & $f^{100}(x_k)$   \\ 
\hline 
 0&  $(-2,10)$& $(0.6033,0.7958,0.8456)$& $(0.5820,0.7793,0.8399)$ &$(0.6354,0.7991,0.8550)$ & $(0.6295,0.7841,0.8502)$ \\  
 1&  $(-1.8614,9.8614)$& $(0.5940,0.7934,0.8432)$& $(0.5715,0.7765,0.8374)$ &$(0.6278,0.7965,0.8526)$ & $(0.6224,0.7812,0.8478)$ \\
2&  $(-1.7210,9.7210)$& $(0.5839,0.7909,0.8407)$& $(0.5601,0.7736,0.8348
)$ &$(0.6196,0.7939,0.8503)$ & $(0.6150,0.7782,0.8454)$\\
3&  $(-1.5794,9.5782)$& $(0.5729,0.7883,0.8382)$& $(0.5478,0.7707,0.8322)$ &$( 0.6109,0.7912,0.8479)$ & $(0.6073,0.7752,0.8430)$\\
4&  $(-1.4381,9.4317)$& $(0.5610,0.7857,0.8356)$& $(0.5345,0.7677,0.8296)$ &$( 0.6016,0.7884,0.8454)$ & $(0.5994,0.7720,0.8404)$\\
5& $(-1.2971,9.2813)$& ($0.5477,0.7831,0.8330$)& ($ 0.5202,0.7647,0.8269$) &($ 0.5915,0.7856,0.8429$) & $( 0.5914,0.7688,0.8379)$ \\
6&  $(-1.1556,9.1258)$& ($ 0.5329,0.7803,0.8302$)& ($0.5051,0.7615,0.8240$) &($0.5806,0.7826,0.8402$) & ($0.5832,0.7654,0.8352$)\\
7& $(-1.0149,8.9666)$ & ($0.5161,0.7775,0.8274$)  & ($0.4901,0.7582,0.8212$)& ($0.5688,0.7796,0.8376$) &($0.5752,0.7619,0.8325$) \\
8&  ($-1.0029,8.9639$)& ($0.5154,0.7774,0.8272$)& ($ 0.4885,0.7581,0.8210$) &($ 0.5680,0.7795,0.8374$) & ($0.5744,0.7618,0.8323$) \\
\hline 
\end{tabular}}
\label{table_5b}
\end{table}

\begin{table}[ht]
\centering
\caption{Performance of Algorithm 1 on Example 5.5 with $t_k=1$}
\centering 
\scalebox{0.57}{
\begin{tabular}{c c c c c c} 
\hline 
Iteration & & & & &\\
number $(k)$& $x_k^\top$& $f^{25}(x_k)$  & $f^{50}(x_k)$ & $f^{75}(x_k)$ & $f^{100}(x_k)$  \\ 
\hline 
$0$ &$(-5,~-5)$ &$(92.9382
,~21.8270,~89.3822)$ & ($109.9507,~29.9507,~101.0619
$) & $(98.4942,~27.3830,~103.8278)$ & $(116,~36,~116)$\\
1&$(-1,~-1)$ &$(39.6046,~0.4938,~36.0490)$ & ($50.3946,~2.3950,~41.5062
$) & $(40.7158,~1.6050,~46.0498)$ & $(51.9996,~3.9999,~51.9999)$\\
\hline 
\end{tabular}}
\label{table_b_quadratic}
\end{table}
\end{example}

\textcolor{black}{In the next two examples (Example \ref{exam_6} and Example \ref{exam_7}), we consider the ordering cone $K$ different from $\mathbb{R}^m_+$ and observe the performance of the proposed Algorithm \ref{algo1}. The Example \ref{exam_6} is a slight modification of Test instance 5.1 discussed in \cite{bouza2021steepest} with respect to a cone other than $\mathbb{R}^m$.}

\begin{example}\label{exam_6}
Consider the function $F:\mathbb{R}\rightrightarrows \mathbb{R}^2$ defined by  
$$F(x)=\{f^1(x),f^2(x),f^3(x),f^4(x)\},$$
where for each $i\in[4],f^i:\mathbb{R}\to\mathbb{R}^2$ is given as
\[
f^i(x)=\begin{pmatrix}
2x^2+ \frac{(i-3)}{2}+ 4x\\
\tfrac{x}{2}\cos(x) -\frac{(i-3)}{2} \sin x\end{pmatrix}.
\]
The cone is $K$ given by $K=\{(y_1,y_2)^\top\in\mathbb{R}^2:5y_1-y_2\geq 0,-9y_1+10y_2\geq 0\}$.

\begin{figure}[ht]
\centering
\mbox{\subfloat[The value of $F$ at each iteration generated by Algorithm \ref{algo1} for the initial point $x_0=4.3000$ for Example \ref{exam_6}]{ \includegraphics[width=0.48\textwidth]{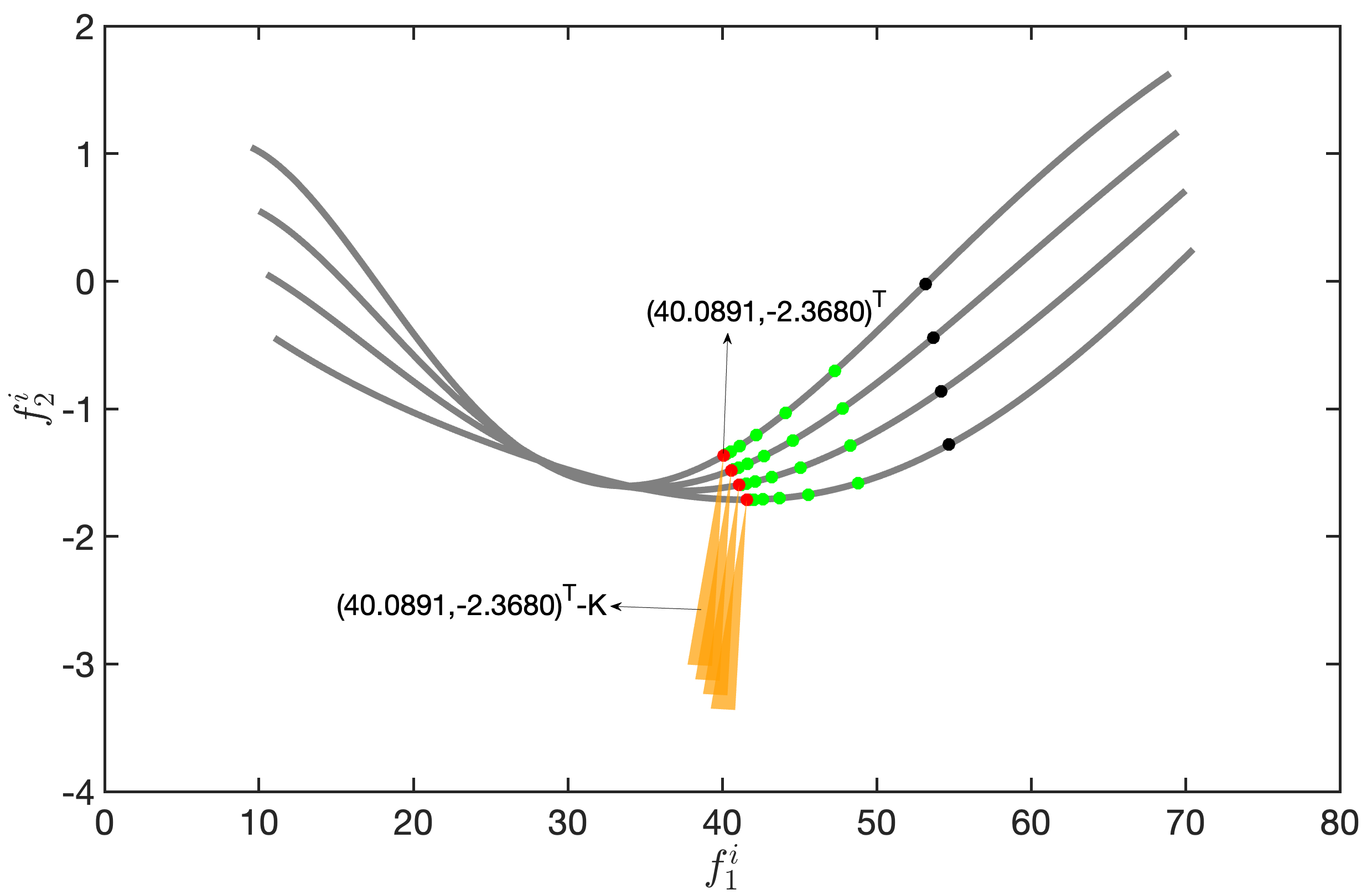}\label{figure6a}}\quad
\subfloat[The value of $F$ for three different initial points at each iteration generated by Algorithm \ref{algo1} for Example \ref{exam_6}]{\includegraphics[width=0.48\textwidth]{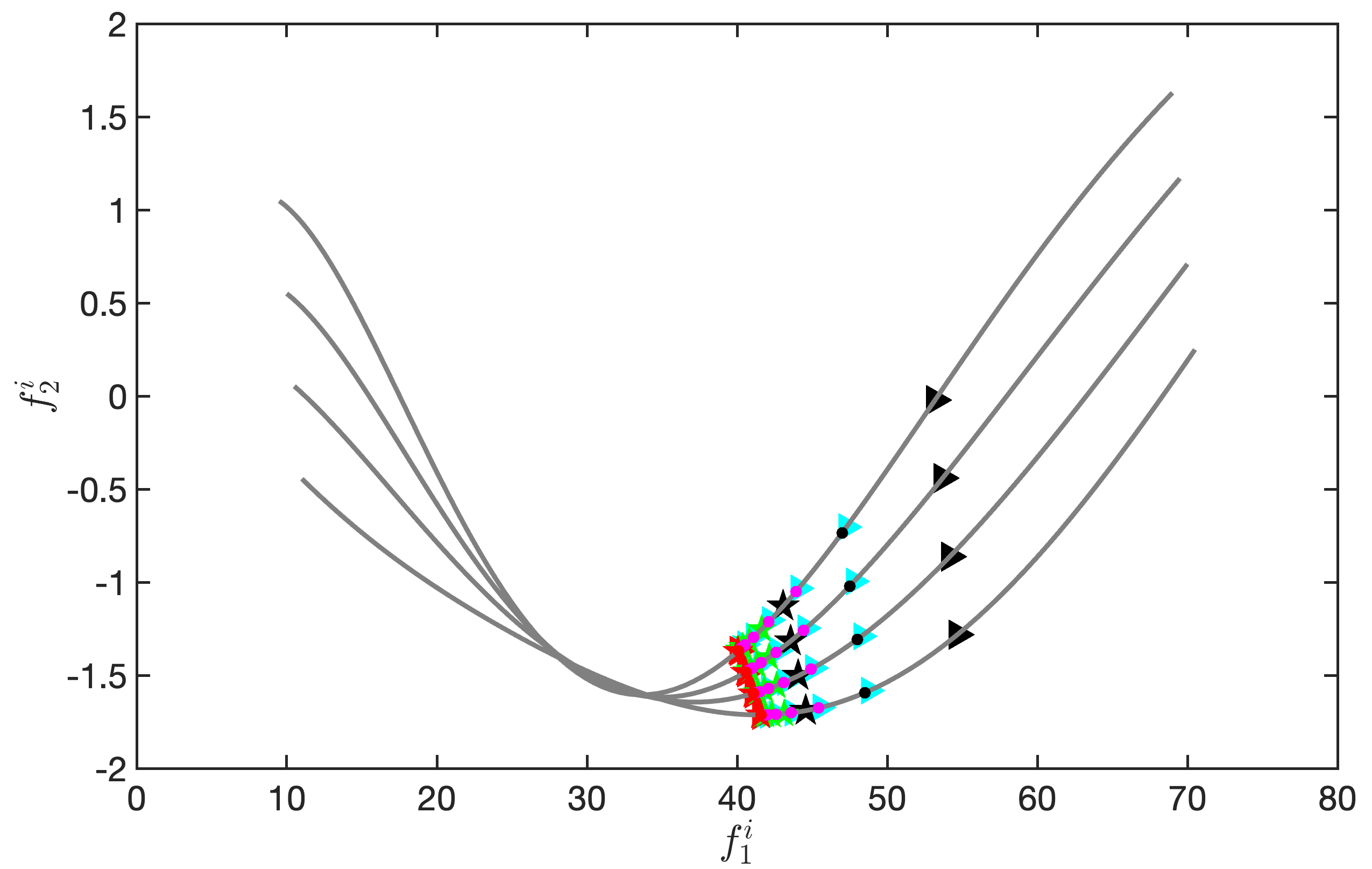} \label{figure6b}}\quad}
\caption{Obtained output of Algorithm \ref{algo1} for Example \ref{exam_6}}
\label{figure6} 
\end{figure}

\textcolor{black}{Figure \ref{figure6} shows the behaviour of Algorithm \ref{algo1} for different initial points within the set $[2.3350,4.4010]$. The collection of objective function values for all $x\in [2.3350,4.4010]$ transverse a surface as shown in Figure \ref{figure6}. Firstly, Figure \ref{figure6a} depicts the sequence of iterates $\{F(x_k)\}$ generated by Algorithm \ref{algo1} for the selected initial point $x_0=4.300$. In Figure \ref{figure6a}, the points depicted with red color collectively constitute a local weakly minimal point of $F$-values as the set $(40.0891,-2.3680)^\top-K$ does not contain any element of the image set of $F$ other than $(40.0891.-2.3680)$ for all $x\in[2.3350,4.4010]$.} 

\textcolor{black}{Subsequently, in Figure \ref{figure6b}, we exhibit the output sequence $\{F(x_k)\}$ generated by Algorithm \ref{algo1} for three randomly chosen initial points.} 

The performance of Algorithm \ref{algo1} on Example \ref{exam_6} is shown in Table \ref{table6}. Additionally, a comparison of the results of NM with the SD method for set optimization is presented in Table \ref{table6}. The values in Table \ref{table6} show that the proposed method performs better than the existing SD method.

\begin{table}[!htb]
\centering
\caption{Performance of Algorithm \ref{algo1} on Example \ref{exam_6}}
\centering 
\scalebox{0.66}{
\begin{tabular}{c c c c} 
\hline 
Number of& Algorithm  &Iterations &CPU time\\
initial points&& (Min, Max, Mean, Median, Mode, SD)  &(Min, Max, Mean, Median, $\lceil\text{Mode}\rceil$, SD)   \\ 
\hline 
$100$ & NM ($t_k=1$) &($1,~3,~1.5500,~1,~1,~0.8333$) & ($388.6205,~428.6387,~409.5273,~410.3945,~388,~11.6303$) \\ 
& NM& ($1,~8,~7.0400,~7,~8,~0.8980$) & ($0.1778,~3.6665,~2.3651,~ 2.3388,~1,~ 0.3064$) \\ 
 & SD& ($3,~10,~8.7800,~9,~9,~0.9813$) & ($1.0957,~3.8509,~2.6506,~2.6638,~2,~0.3428$) \\ 
\hline 
\end{tabular}}
\label{table6}
\end{table}
\end{example}

\begin{example}\label{exam_7}
 Consider the set-valued function $F:\mathbb{R}^2\rightrightarrows \mathbb{R}^2$ defined as
$$F(x)=\{f^1(x),f^2(x)\ldots,f^{100}(x)\},$$
where for each  $i\in[100],$ the function $f^i:\mathbb{R}^2\to\mathbb{R}^2$ is given as
\[
f^i(x)=\begin{pmatrix}
x_1^2+\sin(x_1)+x_1^2\cos(x_2)+0.25\cos\left(\tfrac{2\pi(i-1)}{100}\right)\sin^2\left(\tfrac{2\pi(i-1)}{100}\right)+e^{x_1+x_2}+x_2^2\\
2x_1^2+x_2^2\cos(x_1)+0.25\cos^2\left(\tfrac{2\pi(i-1)}{100}\right)\sin\left(\tfrac{2\pi(i-1)}{100}\right)+\cos(x_2)+e^{x_1+x_2}+2x_2^2\\
\end{pmatrix}.
\]  
The cone $K$ is given by $K =\{(z_1,z_2)^\top\in\mathbb{R}^2:2y_1-6y_2\geq 0, -6y_1+7y_2\geq 0\}$.
\begin{figure}[ht]
\centering
\subfloat[~The value of $F$ at each iteration generated by Algorithm \ref{algo1} for the initial point $x_0=(0.5000,-0.5000)^\top$ for Example \ref{exam_7}]{ \includegraphics[width=0.46\textwidth]{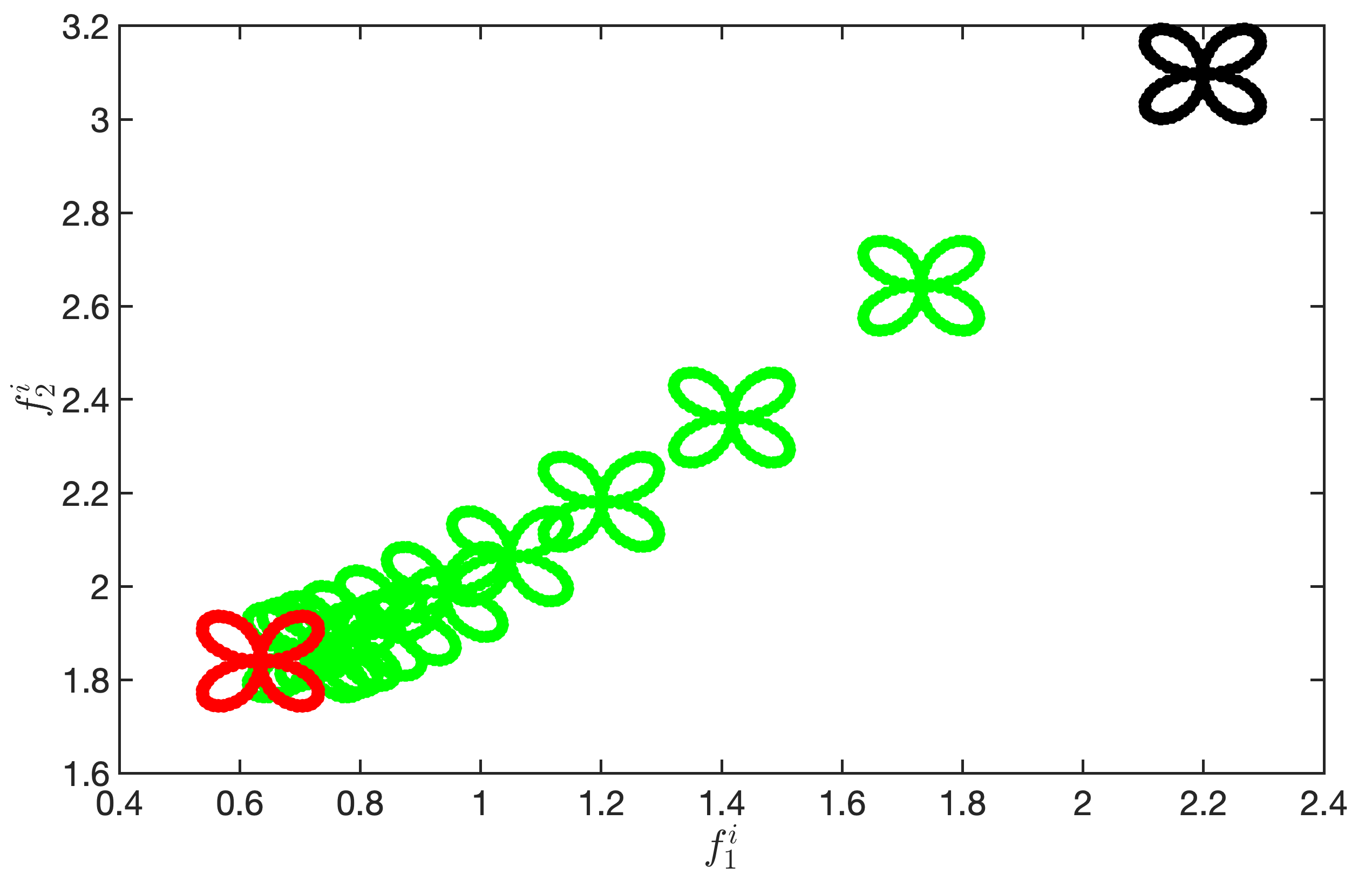}\label{figure7a}}
\qquad
\subfloat[~The value of $x_k$ at each iteration generated by Algorithm \ref{algo1} for for the initial point $x_0=(0.5000,-0.5000)^\top$ for Example \ref{exam_7}]{\includegraphics[width=0.46\textwidth]{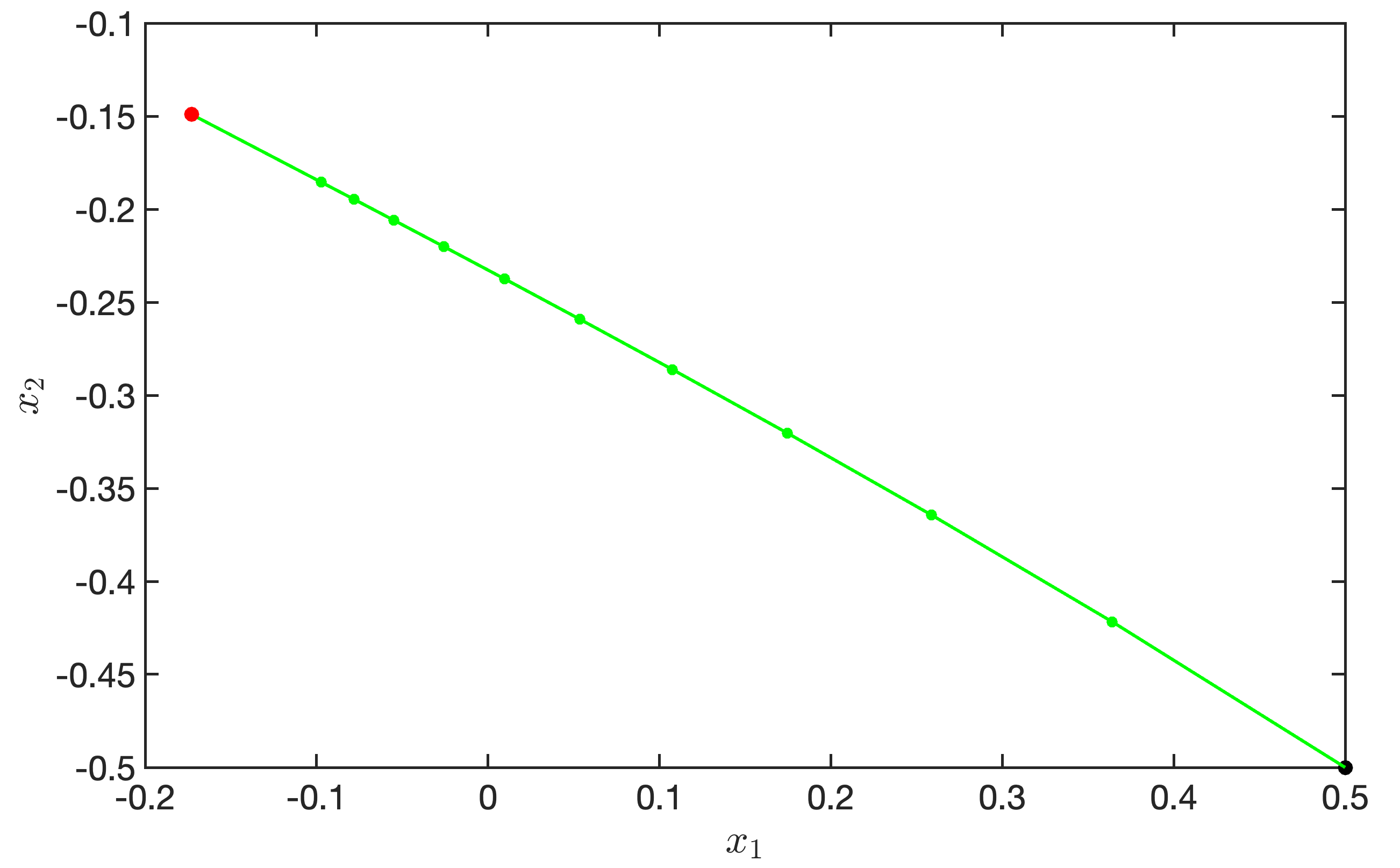} \label{figure7b}}
\qquad
\subfloat[~The value of $F$ at each iteration generated by Algorithm \ref{algo1} for three different randomly chosen initial points  for Example \ref{exam_7}]{\includegraphics[width=0.46\textwidth]{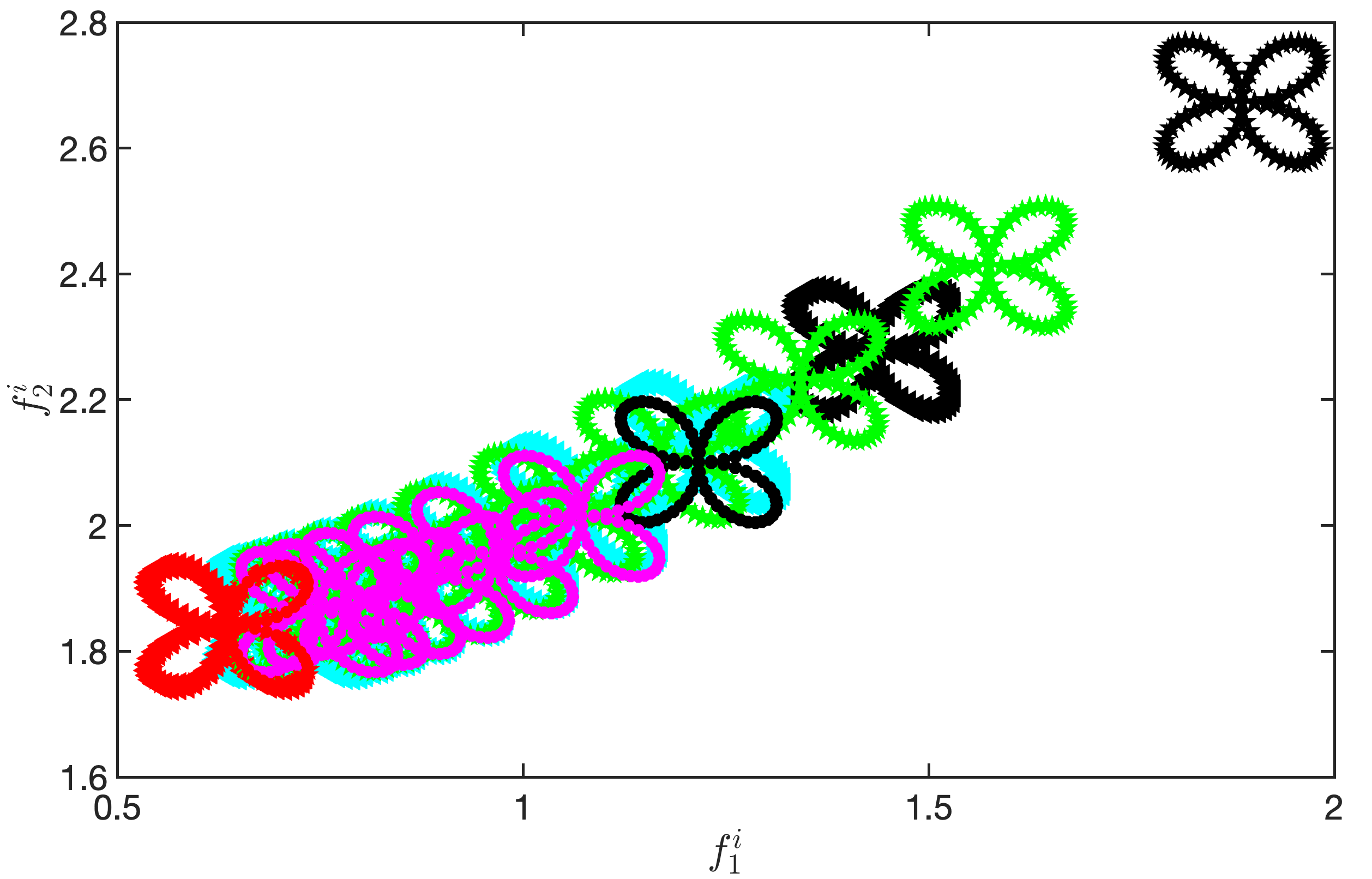} \label{figure7c}}
\qquad
\subfloat[~The value of $x_k$ at each iteration generated by Algorithm \ref{algo1} for three different randomly chosen initial points  for Example \ref{exam_7}]{\includegraphics[width=0.46\textwidth]{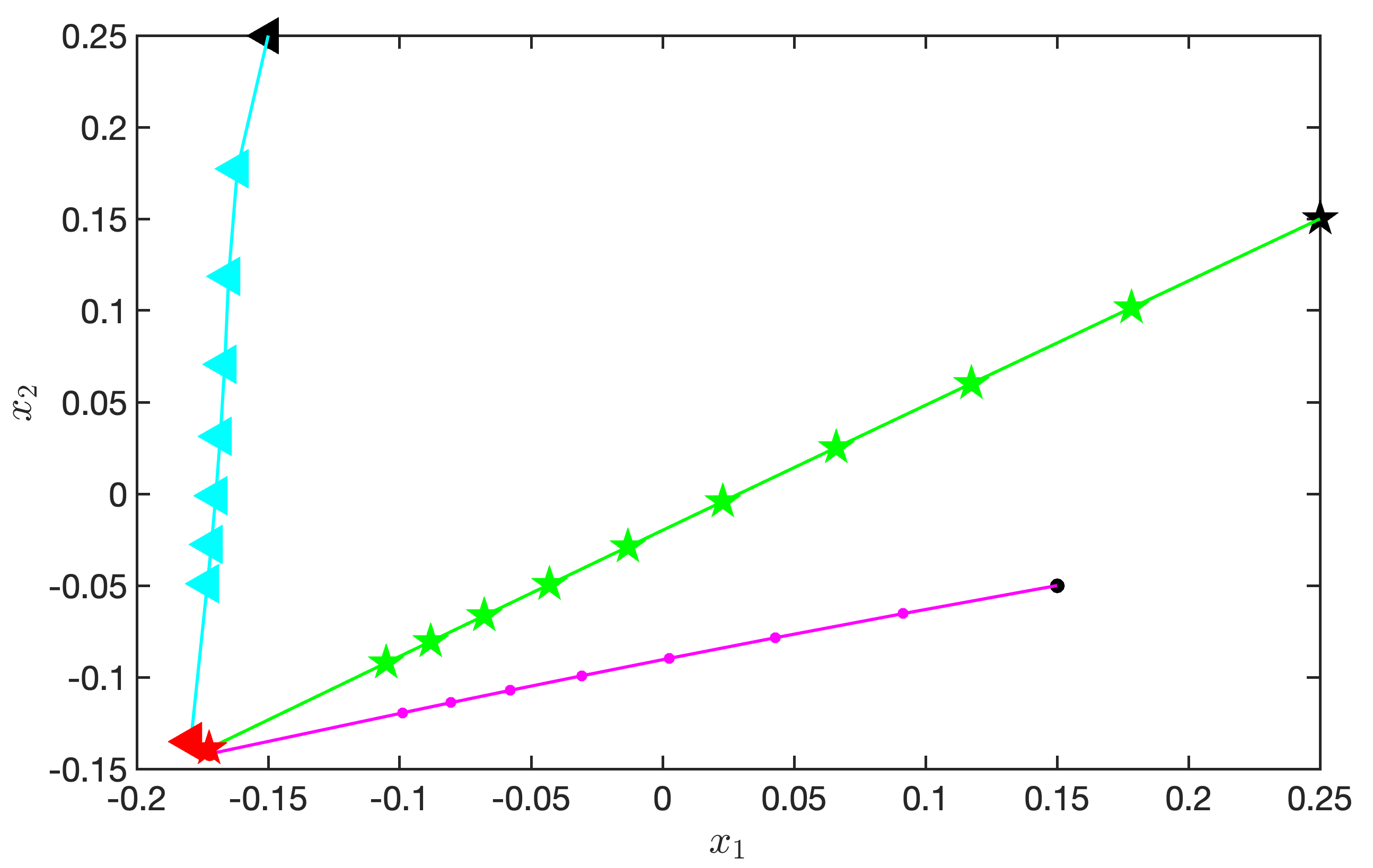} \label{figure7d}}
\caption{Obtained output of Algorithm \ref{algo1} for Example \ref{exam_7}}
\label{figure7} 
\end{figure}

The output of Algorithm \ref{algo1} on Example \ref{exam_7} is depicted in Figure \ref{figure7}. Figure \ref{figure7a} depicts the sequence $\{F(x_k)\}$ generated by Algorithm \ref{algo1} for the starting point $x_0 = (0.5000, -0.5000)^\top$. Figure \ref{figure7b} exhibits the sequence of iterates generated by Algorithm \ref{algo1} for three different randomly chosen starting points.

The performance of Algorithm \ref{algo1} for Example \ref{exam_7} is shown in Table \ref{table7}. A comparison of the results of NM with the SD method for set optimization is presented in Table \ref{table7}. The values in Table \ref{table7} show that the proposed method performs better than the existing SD method.  

\begin{table}[ht]
\centering
\caption{Performance of Algorithm \ref{algo1} on Example \ref{exam_7}}
\centering 
\scalebox{0.69}{
\begin{tabular}{c c c c} 
\hline  
Number of& Algorithm  &Iterations &CPU time\\
initial points&& (Min, Max, Mean, Median, Mode, SD)  &(Min, Max, Mean, Median, $\lceil\text{Mode}\rceil$, SD)   \\ 
\hline 
$100$  & NM ($t_k=1$) &($1,~2,~1.0100,~1,~1,~0.1000$) & ($5.4578,~12.5807,~5.7996,~5.7228,~5,~0.6878$) \\ 
& NM& ($5,~12,~10.6600,~10,~10,~ 1.4718$) & ($47.8361,~109.4561,~80.7287,~71.5817,~38,~  12.6247$) \\ 
 & SD& ($7,~14,~10.5900,~10,~10,~1.8592$) & ($61.5775,~113.1363,83.6014,~79.4974,~46,~15.1152$) \\ 
\hline 
\end{tabular}}
\label{table7}
\end{table}

\end{example}

\section{Conclusion}\label{section6}
In this paper, we studied set optimization problems with respect to the lower set less relation, where the set-valued objective mapping is given by finitely many twice continuously differentiable strongly convex functions. We have proposed a Newton method (Algorithm \ref{algo1}) to generate a sequence of iterates that converges to a weakly minimal solution of the problem. To generate the sequence for a tactfully chosen initial point, at each iteration $k$, we choose an element $a_k$ from the partition set $P_{x_k}$of the current iterate $x_k$, and then with the help of the \textcolor{black}{concepts in \cite{drummond2014quadratically,fliege2009newton}} we figured out the Newton direction $u_k$ (Step \ref{step3}) of the vector optimization problem \eqref{vp_equation} corresponding to $a_k$; this decent direction has been used to find the next iterate $x_{k+1}$. Algorithm \ref{algo1} kept generating iterates until the stopping condition (Step \ref{step4}) was met. The employed stopping condition is a necessary optimality condition (Proposition \ref{prop_minimal}) of weakly minimal or stationary points of the considered problem \eqref{sp_equation}. 

The well-definedness and convergence analysis (Subsection \ref{conv}) of the proposed Algorithm \ref{algo1} for the derived Newton method has been discussed in detail. Towards certifying the well-definedness, we have ensured the 
existence of $(a^k,u_k)$ in Step \ref{step3} and the existence of a step length $t_k$ in Step \ref{step5} (Proposition \ref{armijo}). In the convergence analysis of Algorithm \ref{algo1}, we derived  
\begin{enumerate}[(i)]
\item an equivalent condition of nonstationarity of a point (Proposition \ref{critical}), 
\item boundedness of the sequence of generated Newton direction (Proposition \ref{bounded}), 
\item convergence of the generated sequence of iterates (Theorem \ref{convergence}) under a regularity condition (Definition \ref{regular_condition}), 
\item superlinear convergence to a stationary point (Theorem \ref{superlinear}) of the generated sequence under a regularity condition and uniform continuity of the Hessian matrices, and 
\item local quadratic convergence of the generated sequence under a regularity condition and Lipschitz continuity of the Hessian matrices (Theorem \ref{quadratic}).
\end{enumerate}

Finally, we tested the performance of the proposed Newton method on some existing and freshly introduced numerical test problems in section \ref{section5}. It is found that the proposed Newton method outperforms the existing steepest descent method for strongly convex cases. 

In this paper, we have used the lower set less relation to compare images of the set-valued objective function of \eqref{sp_equation}. \textcolor{black}{Future research can be carried out by testing the proposed Algorithm \ref{algo1} with other set relations such as certainly less order relation, possibly less order relation, min-max less order relation, and min-max certainly less order relation (see \cite{jahn2011new}). A comparison between these order relations is also given in \cite{jahn2011new}. While working with min-max less order relation and min-max certainly less order relation, one may note that usual derivative concepts, like epiderivatives or coderivatives, may not be suitable for optimality concepts. Thus, new derivative concepts need to be developed that align with these order relations.}

\textcolor{black}{As a future work, one can try to find the set of Minimal elements for the proposed work using the approaches discussed in \cite{jahn2006graef,gunther2018new,gerth1990nonconvex,younes1993studies,eichfelder2013ordering} and references therein. A comparison of their performance with different sorting functions can be observed.} Moreover, in this paper, we have used Gerstewitz's scalarizing function for ordering vectors.  Further research can be performed on various other scalarizing functionals, such as separating functionals \cite{gerth1990nonconvex} with uniform level sets, Hiriart-Urruty functional \cite{hiriart1979tangent}, and Drummond-Svaiter functional \cite{drummond2005steepest}. Moreover, a comparison study can be made in the future, showing the performance profile of Newton’s method with different scalarizing functionals. 

\textcolor{black}{Furthermore, one can try to extend the proposed Newton method to different step-size conditions, such as strong Wolfe or Armijo-Wolfe conditions, and a comparison between the performance of the methods with different step-size conditions can be observed. Global convergence of the proposed Newton's method can be observed by analyzing the work done in the articles \cite{sun2002global,gonccalves2022globally,lai2020q}}.

\section*{Disclosure statement}
There were no conflicts of interest reported by the authors of this paper.

\section*{Acknowledgement}
Debdas Ghosh acknowledges the financial support of the research grants MATRICS (MTR/2021/000696) and Core Research Grant (CRG/2022/001347) by the Science and Engineering Research Board, India.

\bibliographystyle{tfnlm}
\bibliography{MyBIB}

\end{document}